\magnification=1200
\input amstex
\documentstyle{amsppt}
\def\bl {\bold\Lambda}
\def\wt#1{\widetilde {#1}}
\def\ssm {\smallsetminus}
\def\SO {\Cal O}
\def\SD {\Cal D}
\def\SE {\Cal E}
\def\SL {\Cal L}
\def\SF {\Cal F}
\def\SH {\Cal H}
\def\SP {\Cal P}
\def\SQ {\Cal Q}
\def\SR {\Cal R}
\def\SU {\Cal U}
\def\SW {\Cal W}
\def\SK {\Cal K}
\def\bwq {\bold{\widetilde Q}}
\def\bset {\text{\bf Set}}
\topmatter
\title Degeneration of moduli spaces and generalized theta functions  
\endtitle
\author Xiaotao Sun   \endauthor
\address Institute of Mathematics, Academia Sinica, 
Beijing 100080, China \endaddress
\email xsun$\@$math08.math.ac.cn\endemail
\thanks The work partly supported by a grant of Academia Sinica for
overseas Chinese scholar who just returned to China  
\endthanks
\endtopmatter
\document

\heading  Introduction \endheading

Let $C$ be a smooth projective curve of genus $g$ and $\SU_C$ the moduli
space of semistable vector bundles of rank $r$ and degree $d$ on $C$.
There is a natural ample line bundle $\Theta$ on $\SU_C$ that we call
theta line bundle of $\SU_C$, which generalises the line bundle on the
jacobian of $C$ defined by Riemann theta divisor [DN]. A section of
$\Theta^k$ over $\SU_C$ is called a generalised theta function of order
$k$. This definition of theta line bundle and generalised theta functions
can be generalised to the moduli spaces of semistable torsion free sheaves
of rank $r$ and degree $d$ on singular curves.
A natural problem suggested by the conformal field theory is to study
the space $H^0(\SU_C,\Theta^k)$ by relating it to space of generalised
theta functions associated with a smooth curve of genus $g-1$.

We consider a family of curves $f:\Cal X\to T$ of genus $g$, whose
singular fibre $\Cal X_0=X$ is irreducible, smooth except for a single
node, so that its normalisation $\wt X$ is a smooth curve of genus $g-1$.
There
exists a moduli scheme $\Cal M\to T$ such that $\Cal M_t$
for any $t\in T$ is the moduli space $\SU_{\Cal X_t}$ of semistable
torsion free sheaves
of rank $r$ and degree $d$. One can define a line bundle on $\Cal M$ such
that its restriction on $\Cal M_t=\SU_{\Cal X_t}$ is the theta line bundle
$\Theta_t$ on $\SU_{\Cal X_t}$. Moreover, if we have a vanishing theorem
$H^1(\Theta^k_t)=0$ for any $t\in T$, one would have that
$dim(H^0(\Theta^k_t))$ is constant. Thus we need to relate the space
$H^0(\SU_X,\Theta^k)$ with the spaces of generalised theta functions
associated with $\wt X$. Let $x_0$ be the node of $X$ and $\pi:\wt X
\to X$ the normalisation of $X$ with $\pi^{-1}(x_0)=\{x_1,x_2\}$.
The expected factorisation rule is
$$H^0(\SU_X,\Theta^k)=\bigoplus_{\mu}H^0(\SU^{\mu}_{\wt
X},\Theta_{\mu}),\tag A$$
where $\mu$ runs through a certain indexing set depending on $k$,
$\SU^{\mu}_{\wt X}$ is the moduli space of parabolic vector bundles of
rank $r$ and degree $d$ on $\wt X$ with parabolic structures at $x_1$ and
$x_2$ (with weights depending on $\mu$) and $\Theta_{\mu}$ is the
generalised theta line bundle. It is clear that to carry through the
induction on genus one has to start with moduli spaces of parabolic
torsion free sheaves of rank $r$ on a nodal curve $X$ with parabolic
structure at a finite number of smooth points and prove a factorisation
rule for generalised theta functions on them, as well as a vanishing
theorem for $H^1$. This was done in the case of rank two by [NR]. We
will treat the general case of any rank in this paper.

{\noindent Now we are going to state the main result. First, some preliminaries:}

{\noindent (1)} Let $X$ be an irreducible projective curve of genus $g$,
smooth but for one node $x_0$. Let $\pi:\wt X\to X$ be the normalization 
of $X$, and $\pi^{-1}(x_0)=\{x_1,x_2\}$.

{\noindent (2)} Let $I$ be a finite set of smooth points on $X$. Fix integers $d,$
$k,$ $r$ and   
$$\aligned
\vec a(x)&=(a_1(x),a_2(x),\cdots,a_{l_x+1}(x))\\
\vec
n(x)&=(n_1(x),n_2(x),\cdots,n_{l_x+1}(x))\endaligned$$
with $0\leq a_1(x)<a_2(x)<\cdots<a_{l_x+1}(x)\le k$ for each $x\in I$. 
Take $(\alpha_x)_{x\in I}\in \Bbb Z^I_{\geq 0}$ and $\ell >0$
satisfying 
$$ \sum_{x\in
I}\sum^{l_x}_{i=1}d_i(x)r_i(x)+r\sum_{x\in
I}
\alpha_x+r\ell=k(d+r(1-g)),\tag{$*$}$$
where $d_i(x)=a_{i+1}(x)-a_i(x)$ and $r_i(x)=n_1(x)+\cdots+n_i(x).$

{\noindent (3)} Let $\SU_X$ be the moduli space of ($s$-equivalence classes of) 
parabolic torsion free sheaves of rank $r$ and degree $d$ on $X$, with parabolic 
structures of type $\{\vec n(x)\}_{x\in I}$ at points $\{x\}_{x\in I},$ 
semistable with respect to the weights $\{\vec a(x)\}_{x\in I}$. The
definitions can be extended to cover the case that $l_x=0$ for $x\in Q\subset I$
(Remark 1.1).

{\noindent (4)} For $\mu=(\mu_1,\cdots,\mu_r)$  with $0\le\mu_r\le\cdots\le\mu_1\le k-1,$
let $$\{d_i=\mu_{r_i}-\mu_{r_i+1}\}_{1\le i\le l}$$ be the subset of nonzero
integers in $\{\mu_i-\mu_{i+1}\}_{i=1,\cdots,r-1},$ and for $j=1,2$ set
$$\aligned
\vec a(x_j)&=(\mu_r,\mu_r+d_1,\cdots,\mu_r+\sum^{l-1}_{i=1}d_i,\mu_r+
\sum^l_{i=1}d_i)\\
\vec n(x_j)&=(r_1,r_2-r_1,\cdots,r_l-r_{l-1}).\endaligned$$
Let $\SU^{\mu}_{\wt X}$ be the moduli space of semistable parabolic bundles
on $\wt X$ with parabolic structures of type $\{\vec n(x)\}_{x\in I\cup\{x_1,x_2\}}$ at points 
$\{x\}_{x\in I\cup\{x_1,x_2\}}$ and weights
$\{\vec a(x)\}_{x\in I\cup\{x_1,x_2\}}.$ We can extend the definition to cover
the case that $l=0$, namely, $\mu_1=\cdots=\mu_r.$

{\noindent (5)} For any data $\omega=(k,r,d,\ell,I,\{\vec a(x),\vec n(x),\alpha_x\}_{x\in I})$ 
satisfying the condition $(*)$, we will define a natural ample line bundle
$$\Theta_{\SU_X}=\Theta(k,r,d,\ell,I,\{\vec a(x),\vec n(x),\alpha_x\}_{x\in I})$$
on $\SU_X$, and $\Theta_{\SU^{\mu}_{\wt X}}$ is defined similarly with
$\alpha_{x_1}=\mu_r$ and $\alpha_{x_2}=k-\mu_1$.

\proclaim{Factorization theorem} There exists a (noncanonical) isomorphism
$$H^0(\SU_X,\Theta_{\SU_X})\cong\bigoplus_{\mu}H^0(\SU^{\mu}_{\wt
X},\Theta_{\SU^{\mu}_{\wt X}})$$
where $\mu=(\mu_1,\cdots,\mu_r)$ runs through the integers $0\le\mu_r\le\cdots\le
\mu_1\le k-1.$\endproclaim

\proclaim{Vanishing theorem} (1) Suppose that $C$ is a smooth projective curve of genus $g\ge 2$. 
Then $H^1(\SU_C,\Theta_{\SU_C})=0.$ (2) Assume that $g\ge 3$. Then $H^1(\SU_X,\Theta_{\SU_X})=0.$
\endproclaim

The Factorization theorem is proved in {\S 4} (Theorem 4.1) and the Vanishing theorems are proved
in {\S 5} (Theorem 5.1 and Theorem 5.3). Next we describe
briefly the main steps in the proof of the main theorems.

We adopt a variant of a concept in [B1], $GPS$, to relate $\SU_X$ with a 
suitable moduli space $\SP$ of $GPS$ on $\wt X$. Such a $GPS$ of rank $r$ is 
given by a pair $(E,Q)$ where $E$ is a sheaf, torsion free outside 
$\{x_1,x_2\}$, of rank $r$ on $\wt X$ and $Q$ a $r$-dimensional quotient of
$E_{x_1}\oplus E_{x_2}$ such that the torsion of $E$ injects to $Q$. Given
such a $GPS$, one defines a torsion free sheaf $F$ on $X$ by the exact sequence
$$0\to F\to \pi_*E\to\,_{x_0}Q\to 0$$
where $_{x_0}Q$ is the skyscraper sheaf on $X$ with support $\{x_0\}$ and 
fibre $Q$. One can define the notion of a semistable $GPS$, and prove that
$F$ is semistable iff $(E,Q)$ is semistable. All this goes through if there
are additional parabolic structures at $\{x\}_{x\in I}.$ There is therefore
a morphism $\phi:\SP\to\SU_X,$ which is actually the normalization 
of $\SU_X$ ({\S 2}).

Set $\Theta_{\SP}=\phi^*\Theta_{\SU_X},$ we will characterize the image
of $H^0(\SU_X,\Theta_{\SU_X})$ in $H^0(\SP,\Theta_{\SP}).$  Our strategy is
to consider the filtrations ($j=1,2$)
$$\SP\supset\SD_j:=\SD_j(r-1)\supset\cdots\supset\SD_j(a)\supset\SD_j(a-1)
\supset\cdots\supset\SD_j(0)$$
of subvarieties of $\SP$ (Notation 2.3 in {\S 2}), and the filtration
$$\SU_X\supset\SW_{r-1}\supset\cdots\supset\SW_a\supset\SW_{a-1}\supset
\cdots\supset\SW_0$$
of subvarieties of $\SU_X$ (Notation 2.4 in {\S 2}). We will prove in {\S 3}
and {\S 4} that $\SP$, $\SD_j(a)$ are reduced, irreducible, normal with only
rational singularities, and $\SU_X$, $\SW_a$  are seminormal 
(Proposition 3.2 and Theorem 4.2). Moreover, we will prove that the 
restriction $\phi_a$ of $\phi$ gives the normalization $\phi_a:\SD_1(a)\to\SW_a$ 
of $\SW_a$, and 
$\phi_a^{-1}(\SW_{a-1})=\SD_1(a)\cap\SD_2\cup\SD_1(a-1)$ (Proposition 2.1).
All of these properties are essentially used to prove that there exists a
(noncanonical) isomorphism 
$H^0(\SU_X,\Theta_{\SU_X})\cong H^0(\SP,\Theta_{\SP}(-\SD_2))$ in {\S 4} 
(Proposition 4.3). Note that
Proposition 2.1 is essential for the story. To prove it for the general
rank case, we have to clarify a fact: if $\SR_a$ (Notation 2.4 in {\S 2})
are saturated sets for the quotient map ? 
We prove that $\SR_a$ are indeed saturated sets for the quotient map 
(Lemma 2.6), which is not known in [NR] and [S2] (See Notation 3.1 of [NR] 
and the $\lq$Remarque' on page 172 of [S2]), thus we can even simplify 
the arguments of [NR] for the case of rank two by using our lemma.

Let $\wt\SR_F$ be the variety parametrizing a certain locally universal family 
of rank $r$ vector bundles $\SE$ on $\wt X$ with degree $d$ and
parabolic structures at $\{x\}_{x\in I},$  $\SU_{\wt X}$ is a geometric
invariant theory (GIT) quotient of the semistable points of $\wt\SR_F$
with respect to the action of a suitable reductive group and certain 
linearization by a line bundle $\hat\Theta.$  Let 
$\rho:\wt\SR'_F\to\wt\SR_F$ denote the grassmannian bundle of $r$-dimensional
quotients of $\SE_{x_1}\oplus\SE_{x_2}.$ One will see that
$$H^0(\SP,\Theta_{\SP}(-\SD_2))
=H^0(\wt\SR'_F,\rho^*\hat\Theta\otimes\bold L)^{inv}
=H^0(\wt\SR_F,\hat\Theta\otimes\rho_*\bold L)^{inv}$$
where $\bold L$ is essentially the line bundle $\SO(k-1)$ along the fibres of
the grassmannian bundle, and $\{\quad\}^{inv}$ denotes a space of invariants
for the group action. The computation of $\rho_*\bold L$ amounts to the 
following classical problem in representation theory: Let $Gr$ be the
grassmannian of $r$-dimensional subspaces of $\Bbb C^{2r}$ and $m$ a positive
integer, how to decompose the irreducible representation $H^0(Gr,\SO(m))$
of $GL(2r)$ into irreducible representations of 
$GL(r)\times GL(r)\subset GL(2r)$ (Lemma 4.5). The factorization theorem
follows from this.

We turn next to the vanishing theorem (1) for a smooth curve $C$. For the
given data $\omega$ satisfying the condition $(*)$, one has a line bundle
$\hat\Theta_{\omega}$ on $\wt\SR_F$, $\SU_C$ is the GIT quotient of 
semistable points 
$\wt\SR^{ss}_{\omega}\subset\wt\SR_F$ with respect to the
action of $SL(n)$ ($n=d+r(1-g)$) and the linearization by the line bundle
$\hat\Theta_{\omega}$, which descends to the ample line bundle 
$\Theta_{\SU_C}$ on $\SU_C$. We can write 
$$\hat\Theta_{\omega}=\omega_{\wt\SR_F}\otimes\hat\Theta_{\bar\omega}\otimes
Det^*\Theta_y^{-2}$$
on $\wt\SR_F$ (Proposition 2.2) for a new data $\bar\omega$ satisfying the
condition $(*)$. Let $\SU_{C,\bar\omega}$ be the GIT quotient of semistable 
points $\wt\SR^{ss}_{\bar\omega}\subset\wt\SR_F$ for the $SL(n)$ action under
the new linearization by $\hat\Theta_{\bar\omega}$, which descends to an ample 
line bundle $\Theta_{\bar\omega}$ on $\SU_{C,\bar\omega}.$ Use the fact that
the complements of $\wt\SR^{ss}_{\omega}$, $\wt\SR^{ss}_{\bar\omega}$ and
$\wt\SR^s_{\bar\omega}$ in $\wt\SR_F$ and $\wt\SR^{ss}_{\bar\omega}$ are
of high codimensions (One need here the restriction on genus, see Proposition 5.1), 
we have
$$H^1(\SU_C,\Theta_{\SU_C})=H^1(\wt\SR_F,\hat\Theta_{\omega})^{inv}
=H^1(\SU_{C,\bar\omega},\Theta_{\bar\omega}\otimes Det^*\Theta_y^{-2}\otimes
\omega_{\SU_{C,\bar\omega}})$$
where $Det$ denote the determinant map and $\Theta_y$ the theta bundle on the 
Jacobian $J^d_C$ of $C$. Then we prove that 
$\Theta_{\bar\omega}\otimes Det^*\Theta_y^{-2}$ is ample (Lemma 5.3) and thus
the vanishing of 
$H^1(\SU_{C,\bar\omega},\Theta_{\bar\omega}\otimes Det^*\Theta_y^{-2}\otimes
\omega_{\SU_{C,\bar\omega}})$ by applying a Kodaira-type vanishing theorem 
(Theorem 7.80(f) of [SS]).

The vanishing theorem (2) for the singular curve $X$ is reduced to prove the
vanishing of $H^1(\SP,\Theta_{\SP})$ (Lemma 5.5). There exists a flat morphism
$Det:\SP\to J^d_{\wt X}$ extending the determinant morphism on the open set of 
stable torsion free $GPS$ (Lemma 5.7), and a decomposition
$$(Det)_*\Theta_{\SP}=\bigoplus_{\mu}(Det_{\mu})_*\Theta_{\SU^{\mu}_{\wt X}}$$
where $Det_{\mu}:\SU^{\mu}_{\wt X}\to J^d_{\wt X}$ is the determinant
morphism. Thus $H^1(J^d_{\wt X},(Det)_*\Theta_{\SP})=0$ by using the 
vanishing theorem (1) for smooth curves, and we are left with the task of proving
$R^1Det_*\Theta_{\SP}=0.$ To prove that $H^1(\SP^L,\Theta_{\SP})=0$, where
$\SP^L$ denotes the fibre of $Det$ at any $L\in J^d_{\wt X}$, 
we follow the same line in the proof of the vanishing
theorem (1) except that $Det^*\Theta_y^{-2}$ disappears. 
We do need here the properties that $\SP$ is Gorenstein with only rational
singularities. It also takes more work to prove a formula for the dualizing
sheaf of $\SP$ (Proposition 3.4 and Lemma 5.6).

We introduce the moduli spaces and theta line bundles in {\S 1}. A detail
study of the morphism $\phi:\SP\to\SU_X$ is given in {\S 2}.  We prove
in {\S 3} that $\SP$ and its subvarieties $\SD_j(a)$ ($j=1,2,\,0\le a\le r-1$)
are normal with only rational singularities, and we also prove a formula 
expressing the canonical (dualizing) sheaf of $\SH$ 
(See {\S 2} for the definition) where we need to prove $\SH$ 
is Gorenstein (it is actually a
complete intersection by using a dimension formula for double determinant
varieties). The factorization theorem and the seminormality of $\SU_X$ and
its subvarieties $\SW_a$ ($0\le a\le r-1$) are proved in {\S 4}. {\S 5}
is devoted to the estimation of codimensions and the proof of vanishing
theorems.

\heading \S 1 Moduli spaces and theta bundles \endheading

We introduce the notation in this section by recalling the construction
of moduli spaces and theta bundles, whose proofs are contained in [NR],
where they deal with rank two, but the proof there goes through for any
rank. We also refer to [BR] and [Pa] about the theta
bundles on moduli spaces of parabolic bundles of any rank.

Let $X$ be an irreducible projective curve of genus $g$ over the complex
number field $\Bbb C$, which has at most one node $x_0$. Let $I$ be a
finite set of smooth points of $X$, and $E$ be a torsion free sheaf of
rank $r$ and degree $d$ on $X.$
\proclaim{Definition 1.1} By a
quasi-parabolic structure on E at a smooth
point $x\in X$, we mean a choice
of flag
$$E_x=F_0(E)_x\supset F_1(E)_x\supset\cdots\supset
F_{l_x}(E)_x\supset
F_{l_x+1}(E)_x=0$$
of the fibre $E_x$ of $E$ at $x$.
If, in addition, a sequence of integers
called the parabolic weights
$$0\leq a_1(x)<a_2(x)<\cdots
<a_{l_x+1}(x)\le k$$
are given, we call that
$E$ has a parabolic structure at $x$.\endproclaim

Let
$n_i(x)=dim(F_{i-1}(E)_x/F_i(E)_x)$ and $r_i(x)=dim(E_x/F_i(E)_x).$ Write
$$\aligned
\vec a(x):&=(a_1(x),a_2(x),\cdots,a_{l_x+1}(x))\\
\vec
n(x):&=(n_1(x),n_2(x),\cdots,n_{l_x+1}(x)).\endaligned$$
We use $\vec a$
(resp., $\vec n$) to denote the map $x\mapsto \vec a(x)$ 
(resp.,
$x\mapsto\vec n(x)$) from $I$ to a suitable set. Let $E'$ be a
subsheaf
of
$E$ such that $E/E'$ is torsion free, then the induced
parabolic
structure
on $E'$ is defined as follows: the quasi-parabolic
structure is defined by
 $F_i(E')_x:=F_i(E)_x\cap E'_x$, and the weights by
$a'_j(x)=a_i(x)$
where $i$ is the biggest integer satisfy that
$F_j(E')_x\subset
F_i(E)_x.$
\proclaim{Definition 1.2} The parabolic degree
of a parabolic sheaf $E$
is
$$pardeg(E):=deg(E)+\frac{1}{k}\sum_{x\in
I}\sum^{l_x+1}_{i=1}n_i(x)a_i(x).$$

$E$ is called semistable (resp., stable) for $(k,\vec a)$ if for
any
subsheaf
$E'\subset E$ such that $E/E'$ is torsion free with the
induced parabolic
structure, one has
$$pardeg(E')\leq
\frac{pardeg(E)}{rk(E)}\cdot
rk(E')\,\,(\text{resp., }<).$$\endproclaim

By
a family of rank $r$ parabolic sheaves parametrised by a variety $T$,
we
mean a sheaf $\Cal F_T$ on $X\times T$, flat over $T$, and torsion
free
(with rank $r$ and degree $d$) on $X\times \{t\}$ for every point
$t\in T$,
together with, for each $x\in I$, a flag of subbundles of
$\Cal
F_T|_{\{x\}\times T}$. The following theorem was proved in the
Appendix
of [NR].

\proclaim{Theorem 1.1} There exists a (coarse) moduli
space $\Cal U_X^s(d,r, I, k, \vec a,\vec n)$ of stable parabolic sheaves
$F.$
We have an open immersion $$\Cal U_X^s(d,r,I,k,\vec a,\vec
n)\hookrightarrow
\Cal U_X(d,r,I,k,\vec a,\vec n)$$ where $\Cal U_X(d,r,I,k,\vec
a,\vec n)$ denotes
the space of $s$-equivalent classes of semistable
parabolic sheaves. The
latter is a seminormal projective variety. If $X$ is
smooth, then it is
normal, with
only rational
singularities.\endproclaim

Fix $I$, $k$, $\vec a$ and $\vec n$, we set
$\Cal U_X:=\Cal U_X(d,r,I,k,\vec
a, \vec n)$ and $\Cal U^s_X:=\Cal
U_X^s(d,r,I,k,\vec a,\vec n)$, let us recall
the
construction of $\Cal
U_X$.

Let $\bold Q$ be the Quot scheme of coherent sheaves (of rank $r$
and
degree $d$) over $X$ which are quotients of $\Cal O^n$,
where
$n=d+r(1-g)$. Thus there is on $X\times\bold Q$ a sheaf
$\Cal
F_{\bold Q}$, flat over $\bold Q$, and $\Cal O_{X\times\bold Q}^n
\to
\Cal F_{\bold Q}\to 0$. Let $\Cal F_x$ be the sheaf given by
restricting
$\Cal F_{\bold Q}$ to $\{x\}\times\bold Q$. Let 
$Flag_{\vec n(x)}(\Cal
F_x)$ be the relative flag scheme of type $\vec
n(x)$, and $\Cal R$ be the
fibre product over $\bold Q$:
$$\Cal R=\underset x\in I\to{\times_{\bold
Q}}Flag_{\vec n(x)}(\SF_x).$$

Let $\Cal R^s$ (resp., $\Cal R^{ss}$) be the
open subscheme of 
$\Cal R$ corresponding to stable (resp., semistable)
parabolic
sheaves, which is generated by global sections and its first
cohomology 
vanishes when $d$ is large enough. The variety $\Cal U_X$ is
the good
quotient
of $R^{ss}$ by $SL(n)$ acting through $PGL(n)$. We denote
the projection
by $$\psi :\Cal R^{ss}\to \Cal U_X.$$

Choose an ample line
bundle of degree $1$ on $X$, denoted by $\Cal O_X(1)$
from now on. For
large enough $m$, we have a $SL(n)$-equivariant embedding
$\Cal
R\hookrightarrow \bold G$, where $\bold G$ is defined to
be
$$Grass_{P(m)}(\Bbb C^n\otimes W)\times
\prod_{x\in I}\{Grass_r(\Bbb
C^n)\times Grass_{r_1(x)}(\Bbb C^n)
\times\cdots\times
Grass_{r_{l_x}(x)}(\Bbb C^n)\}$$
where $P(m)=n+rm$, and $W=H^0(\Cal
O_X(m))$. For any $(\alpha_x)_{x\in I}
\in \Bbb Z^I_{\geq 0}$ and $\ell >0$
satisfying 
$$ \sum_{x\in
I}\sum^{l_x}_{i=1}d_i(x)r_i(x)+r\sum_{x\in
I}
\alpha_x+r\ell=kn,\tag{$*$}$$
where $d_i(x)=a_{i+1}(x)-a_i(x)$, we give $\bold G$ the polarisation (using
the obvious notation):
$$\frac{\ell}{m}\times\prod_{x\in
I}\{\alpha_x, d_1(x),\cdots,
d_{l_x}(x)\}$$
and take
the induced polarisation on $\Cal R$. It was proved in [NR] 
that the set
of semistable points for the $SL(n)$ action on $\Cal R$ is
precisely $\Cal
R^{ss}$. One remarks that this fact is independent on 
the choice of
$\vec\alpha:=(\alpha_x)_{x\in I}$ satisfying the condition 
($*$). $\Cal
R^{ss}$ is reduced and irreducible, $\Cal U_X$ is its GIT
quotient. 

For
any family of parabolic sheaf $\Cal F$ of type $\vec n(x)$ at $x\in
I$
parametrised by $T$, we denote the quotients $\Cal
F_{\{x\}\times
T}/F_i(\Cal F_{\{x\}\times T})$ by $\SQ_{\{x\}\times
T,i}$, and we define
$$\Theta_T:=
(det\,R\pi_T\Cal
F)^k\otimes\bigotimes_{x\in I}
\lbrace(det\,\Cal F_{\{x\}\times
T})^{\alpha_x}\otimes\bigotimes^{l_x}_{i=1}
(det\,\SQ_{\{x\}\times
T,i})^{d_i(x)}\rbrace\otimes 
(det\,\Cal F_{\{y\}\times
T})^{\ell}$$
where $\pi_T$ is the projection $X\times T\to T$, and
$det\,R\pi_T\Cal F$
is the determinant bundle defined as
$$\{det\,R\pi_T\Cal
F\}_t:=\{det\,H^0(X,\Cal F_t)\}^{-1}\otimes\{
det\,H^1(X,\Cal F_t)\}.$$
If
we take $T=\Cal R^{ss}$ and $\vec a$, 
$\vec\alpha$,$k$, $\ell$ satisfying
the condition $(*),$ it is easy to check
that $\Theta_{\Cal R^{ss}}$ is a
$PGL(n)$-bundle, which descends to $\Cal
U_X$. Moreover, we have the
following theorem for whose proof we refer to
[NR], [Pa] and
[BR].
\proclaim{Theorem 1.2} There is an unique ample line bundle
$\Theta_{\Cal U_X}=\Theta(k,\ell,\vec a,\vec{\alpha},I)$ 
on $\Cal U_X$
such that for any given family of semistable parabolic sheaf
$\Cal F$
parametrised by $T$, we have
$\phi_T^*\Theta_{\Cal U_X}=\Theta_T$, where
$\phi_T$ is the induced map
$T\to\Cal U_X$\endproclaim

\remark{Remark 1.1}
$(1)$ It is known that the analytic
local ring of $\Cal R^{ss}$ is
determined (up to smooth morphisms) by $\Bbb
C[X,Y]/(XY,YX),$
where $X$,
$Y$ are $r\times r$ matrices (See [Fa] and [S2]). Thus,
by Lemma 3.8 and
Lemma 3.13 of [NR], the seminormality of $\Cal U_X$ 
is equivalent to that
of $\Bbb C[X,Y]/(XY,YX)$, which is known to be
seminormal (See [Tr]).

$(2)$ If we replace, in the construction of theorem 1.2, 
$(det\,\Cal F_{\{y\}\times \SR^{ss}})^{\ell}$ by
$$\bigotimes_{q\in Q}(det\,\Cal F_{\{q\}\times\SR^{ss}})^{\beta_q}
\otimes (det\,\Cal F_{\{y\}\times\SR^{ss}})^{\ell+\ell_0},$$
where $Q$ is a set of smooth points of $X$, and 
$\sum_{q\in Q}\beta_q=-\ell_0,$  we get ample line bundles 
on $\SU_X$, which are all algebraically equivalent to $\Theta_{\SU_X}$.

$(3)$ We can extend the above definitions to cover the case that $l_q=0$
for $q\in Q\subset I$. In this case, $\SU_X$ denotes the moduli space 
of semistable parabolic sheaves with parabolic structures at $\{x\}_{x\in
I\ssm Q}$ and parabolic weights $\{\vec a(x)\}_{x\in I\ssm Q}$.
When $Q=I$, $\Cal U_X$ is the ordinary moduli space
of semistable torsion free sheave (i.e., no quasi-parabolic structure
is considered), the
definition of $\Theta _{\Cal U_X}$ in Theorem 1.2
gives ample line bundles
$\Theta (\vec\alpha,I)$ on $\Cal U_X$, all
of them are algebraically
equivalent to the descendant of 
$$(det\,R\pi_{\Cal R^{ss}}\Cal F)^k\otimes
(det\,\Cal F_{\{y\}\times
\Cal R^{ss}})^{\frac{kn}{r}}.$$
These $\Theta
(\vec\alpha,I)$ will appear naturally in the decomposition
theorems,
induced by the $1$-dimensional repsentations of $GL(r)$.
\endremark

Now we
are going to recall the notion of $\lq$generalised parabolic
sheaf'
(GPS)
and the construction of its moduli space ([B1], [B2] and
[NR]).
We do not define the general notation (as in [B1] and [B2]), but
we
have to consider the sheaves with torsion as in [NR].
Let
$\pi:\widetilde X\to X$ be the normalisation of $X$ and
$\pi^{-1}
(x_0)=\{x_1,x_2\},$ then we have 
\proclaim{Definition 1.3} Let
$E$ be a sheaf on $\widetilde X$, torsion
free of rank $r$ outside
$\{x_1,x_2\}.$ A generalised parabolic structure
on $E$ over the divisor
$x_1+x_2$ is a $r$-dimensional quotient $Q$ 
$$E_{x_1}\oplus
E_{x_2}@>q>>Q\to 0.$$  $(E,Q)$ is said to be a
generalised parabolic sheaf,
namely GPS.\endproclaim

We will consider generalised parabolic sheaves $E$
with, in addition,
parabolic structures at the points of $\pi^{-1}(I)$ (we
will identify
$I$ with $\pi^{-1}(I)$). Furthermore, by a family of GPS over
$T$, we
mean the following
\roster
\item a rank $r$ sheaf $\Cal E$ on
$\widetilde X\times T$ flat over
$T$ and locally free outside
$\{x_1,x_2\}\times T$.
\item a locally free rank $r$ quotient $\Cal Q$ of
$\Cal E_{x_1}\oplus
\Cal E_{x_2}$ on $T$.
\item a flag bundle $Flag(\Cal
E_x)$ on $T$ with given weights for each
$x\in
I$.\endroster
\proclaim{Definition 1.4} A GPS $(E,Q)$ is called
semistable
(resp., stable), if
for every nontrivial subsheaf $E'\subset E$
such that $E/E'$ is torsion
free outside $\{x_1,x_2\},$ we
have
$$pardeg(E')-dim(Q^{E'})\leq
rk(E')\cdot\frac{pardeg(E)-dim(Q)}{rk(E)}
\,\quad (\text{resp.,
$<$}),$$
where $Q^{E'}=q(E'_{x_1}\oplus E'_{x_2})\subset
Q.$\endproclaim

Set $\tilde n=d+r(1-\tilde g)$, where $\tilde g(=g-1)$ is
the genus of
$\widetilde X$, and let $\bold {\widetilde Q}$ be the Quot
scheme of 
coherent sheaves (of degree $d$ and rank $r$) over $\widetilde
X$ 
which are quotients of $\Cal O^{\tilde n}_{\widetilde X}.$ Taking $d$
to be
large enough, we can assume that for any semistable generalsed
parabolic
sheaf $E$ of rank $r$ and degree $d$ we have
$H^1(E(-x_1-x_2-x))=0,\,
x\in \widetilde X,$ which means that $\Bbb
C^{\tilde n}\to H^0(E)$ is an 
isomorphism, $E$ is generated by global
sections and $H^0(E)\to E_{x_1}
\oplus E_{x_2}$ is onto, $E(-x_1-x_2)$ is
generated by global sections.
Let $\Cal F$ be the universal quotient
$\widetilde {\Cal O}^{\tilde n}:=
\Cal O^{\tilde n}_{\widetilde
X\times\bold {\widetilde Q}}\to\Cal F\to 0$
on $\widetilde X\times\bold
{\widetilde Q}$ and
$$\widetilde {\Cal R}':=Grass_r(\Cal F_{x_1}\oplus\Cal
F_{x_2})\times
_{\bold {\widetilde Q}}\biggl\{\underset {x\in
I}\to{\times
_{\bold{\widetilde Q}}}Flag_{\vec n(x)}(\Cal
F_x)\biggr\}.$$
There is a locally universal family of GPS parametrised by
$\widetilde {\Cal R}'$ that we denote by $\Cal E$, which is actually
the
pull back of $\Cal F$ by the natural projection. Let $\tilde
P(m)=\tilde n +rm$ and 
$$\bold {\tilde G}':=
Grass_{\tilde P(m)}(\Bbb
C^{\tilde n}\otimes W)\times Grass_r(\Bbb C
^{\tilde n}\otimes\Bbb
C^2)\times \bold {Flag},$$
where $\bold{Flag}$ denotes the
variety
$$\prod_{x\in I}
\{Grass_r(\Bbb C^{\tilde n})\times
Grass_{r_1(x)}(\Bbb
C^{\tilde n})\times\cdots\times
Grass_{r_{l_x}(x)}(\Bbb
C^{\tilde n})\}.$$
Then we have a $SL(\tilde n)$-equivariant embedding
$\widetilde {\Cal R}'
\hookrightarrow \bold{\tilde G}'.$ Take the
polarisation
$$\frac{(\tilde{\ell}-k)}{m}\times k\times\prod_{x\in
I}
\{\alpha_x,d_1(x),\cdots,d_{l_x}(x)\}$$
such
that
$$ \sum_{x\in I}\sum^{l_x}_{i=1}d_i(x)r_i(x)+r\sum_{x\in
I}
\alpha_x+r\tilde{\ell}=k\tilde n,$$
which is nothing but $(*)$ with
$\tilde n=n+r$ and $\tilde{\ell}=\ell+k.$
Then one proves that the
GIT-semistable (stable) points of $\widetilde
{\Cal R}'$ are precisely the
semistable (stable) generalised parabolic
sheaves, namely $\widetilde{\Cal
R}^{\prime ss}$. Let $\Cal P:=\Cal
P_{\wt X}$ be the GIT quotient of
$\widetilde{\Cal R}^{\prime ss}$
by $SL(\tilde n)$ with the
projection
$$\tilde{\psi}':\widetilde{\Cal R}^{\prime ss}\to \Cal P.$$
One
defines an $s$-equivalence of GPS such that
$$E\sim E'\Longleftrightarrow\text{there exist
$E_1=E,\cdots,E_{s+1}=E'$ with $\overline{o(E_i)}\cap\overline
{o(E_{i+1})}\neq\emptyset,$}\tag1.1$$
where $\overline{o(E_i)}$ denotes
the schematic closure of the orbit of
$E_i$ under $SL(\tilde n).$ It is
clear that if $E_1$ and $E_2$ are stable 
then $E_1\sim E_2$ iff $E_1\cong
E_2$. Then we have 
\proclaim{Theorem 1.3} There exists a (coarse) moduli
space $\Cal P^s$
of stable GPS on $\widetilde X$, which is a smooth
variety. We have an
open immersion $\Cal P^s\hookrightarrow\Cal P$, where
$\Cal P$ is the
space of $s$-equivalence class of semistable GPS on
$\widetilde X$, which
is a reduced, irreducible and normal projective
variety with
rational singularities.\endproclaim

The existence of $\Cal
P$ is known as we have shown above. We will prove
in $\S 3$ that it is
reduced, irreducible and normal with rational
singularities. In fact $\SP$
is the normalisation of $\Cal U_X$ as we
will see in next section. We
complete this section by introducing a sheaf
theoretic description
of
$s$-equivalence of GPS, which was given in Appendix B of [NR] in the
case
of rank $2$. We
enlarge the category of GPS by adopting the following more
general
definition, and assume that $|I|=0$ for
simplicity.

\proclaim{Definition 1.5} A generalised $m$-parabolic
structure on a sheaf
$E$ over
the divisor $x_1+x_2$ is a choice of an
$m$-dimensional quotient $Q$ of 
$E_{x_1}\oplus E_{x_2}$. A sheaf with a
generalised $m$-parabolic
structure will be called an $m$-GPS, or GPS for
short. A GPS $E$ is said
to be semistable (resp., stable) if $E$ is torsion
free outside
$\{x_1,x_2\}$ and
\roster\item if $rank(E)>0$, then for every
proper subsheaf $E'$ such
that $E/E'$ is torsion free outside
$\{x_1,x_2\}$, we have
$$rank(E)(deg(E')-dim(Q^{E'}))\leq
rank(E')(deg(E)-m)\quad (\text{resp.,
$<$})$$
\item if $rank(E)=0$, then
$E_{x_1}\oplus E_{x_2}=Q$ (resp.,
$E_{x_1}\oplus E_{x_2}=Q$ and
$dim(Q)=1$).\endroster\endproclaim

\proclaim{Definition 1.6} If $(E,Q)$ is
a GPS and $rank(E)>0$, we
set
$$\mu_G[(E,Q)]=\frac{deg(E)-dim(Q)}{rank(E)}.$$\endproclaim

It is
useful to think of an $m$-GPS as a sheaf $E$ on $\wt X$ together
with a map
$\pi_*E\to\,_{x_0}Q\to 0$ and $h^0(_{x_0}Q)=m$. Let $K_E$
denote the kernel
of $\pi_*E\to Q$.

\proclaim{Definition 1.7} A morphism of GPS $(E,Q)\to
(E',Q')$
is a sheaf map $E\to E'$ which maps $K_E$ to $K_{E'}$ (and
therefore
induces a map $Q\to Q'$).\endproclaim

\proclaim{Definition 1.8}
Given an exact sequence
$$0\to E'\to E\to E''\to 0$$
of sheaves on $\wt X$,
and $\pi_*E\to Q\to 0$ a generalised parabolic
structure on $E$, we define
the generalised parabolic structures on $E'$
and $E''$ via the
diagram
$$\CD
0@>>> \pi_*E'@>>>\pi_*E@>>> \pi_*E''@>>>  0  \\
@.    @VVV
@VVV        @VVV        @. \\
0@>>> Q'    @>>>   Q   @>>>    Q''  @>>>
0
\endCD$$    
The first horizontal sequence is exact because $\pi$ is
finite, $Q'$
is defined as the image in $Q$ of $\pi_*E'$ so that the first
vertical
arrow is onto, $Q''$ is defined by demanding that the second
horizontal
sequence is exact, and finally the third vertical arrow is onto
by the
snake lemma. We will write
$$0\to (E',Q')\to (E,Q)\to (E'',Q'')\to
0$$
whose meaning is clear.
\endproclaim

\proclaim{Proposition 1.1} Fix a
rational number $\mu$. Then the category
$\Cal C_{\mu}$ of semistable GPS
$(E,Q)$ such that $rank(E)=0$ or,
$rank(E)>0$ with $\mu_G[(E,Q)]=\mu$, is
an abelian, artinian, noetherian
category whose simple objects are the
stable GPS in the
category.\endproclaim

One can conclude, as usual, that
given a semistable GPS $(E,Q)$ it has a
Jordan-Holder filtration, and the
associated graded GPS $gr(E,Q)$ is
uniquely determined by $(E,Q)$. Thus we
have

\proclaim{Definition 1.9} Two semistable GPS $(E_1,Q_1)$ and
$(E_2,Q_2)$
are said to be $s$-equivalent if they have the same associated
graded GPS,
namely,
$$(E_1,Q_1)\sim  (E_2,Q_2)\quad\Longleftrightarrow\quad
gr(E_1,Q_1)\cong
gr(E_2,Q_2).$$\endproclaim

\remark{Remark 1.2} Any stable
GPS $(E,Q)$ with $rank(E)>0$ must
be a GPB (i.e., E is a vector bundle) such that $E_{x_j}\to Q$ ($j=1,2$)
are isomorphisms, and two stable GPS are s-equivalent
iff they are isomorphic. In fact, let $Q_j$ be the image of $E_{x_j}\to Q$ 
and $\bar q:E_{x_1}\oplus E_{x_2}@>q>>Q@>>>Q/Q_1=\bar Q$. Then we define
$E'$ by the exact sequence
$$0\to E'\to E\to _{x_2}\bar Q_2\to 0,$$
where $\bar Q_2=Q_2/Q_1\cap Q_2$ is the image of $E_{x_2}\to Q\to\bar Q$. 
Thus $Q^{E'}
=Q_1$ and $$\mu_G[E',Q_1)]=\mu_G[(E,Q)]+\frac{dim(Q)-dim(Q_1)-dim( Q_2)+
dim(Q_1\cap Q_2)}{rank(E)}.$$
If $(E,Q)$ is stable, we must have $Q_1=Q_2=Q$. One can
imitate the proof of Lemma 4.7 and
Theorem 4.8 in [Gi] to show that this
s-equivalence satisfies the
requirement $(1.1)$.\endremark

\heading \S 2 The normalization of parabolic moduli spaces on a nodal
curve \endheading

Let $X$ be an irreducible
projective curve of genus $g$ and smooth except
for one node $x_0$, and
$\pi:\wt X\to X$ the normalisation,
$\pi^{-1}(x_0)=
\{x_1,x_2\}.$ It is
clear that we have the canonical exact sequence
$$0\to\SO_X\to\pi_*\SO_{\wt
X}\to\,_{x_0}k(x_0)\to 0,$$
where $k(x_0)$ denotes the residue field of
$x_0$, and we will use $_xW$
to denote the ``skyscraper sheaf" supported at
$\{x\}$, with fibre $W$.

Given a GPS $(E,Q)$ on $\wt X$, we have the
exact sequence
$$0\to F\to\pi_*E\to\,_{x_0}Q\to 0.$$
It is clear that $\phi
(E,Q):=F$ (which has the natural parabolic
structures at points of $I$) is
a torsion free sheaf on $X$ of rank
$r$ if and only
if
$$(\text{Tor}E)_{x_1}\oplus(\text{Tor}E)_{x_2}\overset
q
\to\hookrightarrow Q.\tag T$$
Note that, for any sheaf $E$ on $\wt X$, we have 
$deg(\pi_*E)=deg(E)+rank(E)$,  thus $deg(F)=deg(E).$ 

\proclaim{Lemma 2.1} Let $(E,Q)$ satisfy
condition (T), and $F=\phi(E,Q)$
the associated torsion free sheaf on $X$.
We have
\roster\item If $E$ is a vector bundle and the maps $E_{x_i}\to
Q$
are isomorphisms, then $F$ is a vector bundle.
\item If $F$ is a vector
bundle on $X$, then there is an unique $(E,Q)$
such that $\phi(E,Q)=F.$ In
fact, $E=\pi^*F.$
\item If $F$ is a torsion free sheaf, then there is a
$(E,Q)$, with
$E$ a vector bundle on $\wt X,$ such that $\phi(E,Q)=F$ and
$E_{x_2}\to Q$
is an isomorphism. The rank of the map $E_{x_1}\to Q$ is $a$
iff 
$F\otimes\hat\SO_{x_0}\cong \hat\SO_{x_0}^{\oplus a}\oplus
m_{x_0}^{\oplus
(r-a)}.$ 
The roles of $x_1$ and $x_2$ can be
reversed.
\item Every torsion free rank $r$ sheaf $F$ on $X$ comes from a
$(E,Q)$
such that $E$ is a vector
bundle.\endroster\endproclaim

\demo{Proof} Here we only check (3) since we
will need the construction later.
The proof of Lemma 4.6 in [NR] easily
extends to the other statements for
any rank.
Let
$F\otimes\hat\SO_{x_0}\cong \hat\SO_{x_0}^{\oplus
a}\oplus
m_{x_0}^{\oplus(r-a)},$ 
and define a vector bundle $\wt E$ on
$\wt X$ by
$$0\to\text{Tor}(\pi^*F)\to\pi^*F\to\wt E\to 0.$$  
By the
canonical exact sequence
$$0\to\SO_X\to\pi_*\SO_{\wt X}\to\,_{x_0}k(x_0)\to
0,$$
we get (note that $\pi_*\pi^*F=F\otimes\pi_*\SO_{\wt X}$ and $F$
is
torsion free)
$$0\to F\to \pi_*\pi^*F\to\,_{x_0}Q_F\to 0,$$
where
$Q_F:=k(x_0)\otimes_{\SO_X}F$ is a vector space of dimension
$2r-a$.
Consider the diagram
$$\CD
0  @>>> F  @>>> \pi_*\pi^*F @>>>
_{x_0}Q_F @>>>    0\\
@.      @|       @VVV               @VVV
@.\\
0  @>>> F @>d>> \pi_*\wt E  @>>> _{x_0}\wt Q @>>>  0
\endCD$$
where
the vertical arrows are surjections and $\wt
Q=Q_F/\pi_*\text{Tor}(\pi^*F)$
is of dimension $a$. 
Note that $\pi_*\pi^*F@>>>\, _{x_0}Q_F$ induces two
surjective maps
$(\pi^*F)_{x_i}\to Q_F$, so is $\pi_*\wt E\to \wt Q$.
Denote their kernel
by
$K_i$, we have exact sequences 
$$0\to K_i\to\wt
E_{x_i}\to \wt Q\to 0.\tag2.1$$
Let $h:\wt E \to E$ be the Hecke
modification at $x_2$ (See Remark 1.4 of
[NS]) such that
$ker(h_{x_2})=K_2$, where $h_{x_2}:\wt E_{x_2}\to E_{x_2}$
is the induced
map of $h$ between the fibres at $x_2,$ then one has the 
exact sequence
for some $Q_2$ of dimension $r-a$
$$0\to\wt E@>h>>E @>>>\, _{x_2}Q_2@>>>
0,$$
namely, $pardeg(E)=pardeg(F)$ and $E_{x_1}\cong\wt E_{x_1}$. We
define
$Q$ by the exact sequence
$$0\to F @>(\pi_*h)\cdot d>>\pi_*E@>>>\,
_{x_0}Q @>>>0\tag Q$$
which is clearly of dimension $r$ and $\phi(E,Q)=F$.
To check that
the induced $E_{x_1}\oplus E_{x_2}@>q>> Q\to 0$ (by $(Q)$)
satisfies the
requirement, we consider the restriction of $(Q)$ at
$x_0$
$$F_{x_0}@>d_0>>\wt E_{x_1}\oplus\wt E_{x_2}@>h_{x_1}\oplus
h_{x_2}>>
E_{x_1}\oplus E_{x_2}@>q>> Q\to 0$$
which implies that
$Im(d_0)\cap\wt E_{x_i}=K_i$ and
$ker(q_i)=h_{x_i}(K_i),$ where
$q_i:E_{x_i}\to Q.$ Thus $q_1$ has rank
$a$ (Since $h_{x_1}$ is an
isomorphism), and $q_2$ is an isomorphism
(Since
$ker(h_{x_2})=K_2)$).\enddemo

\proclaim{Lemma 2.2} Let $F=\phi(E,Q)$,
then $F$ is semistable if and
only if $(E,Q)$ is semistable. Moreover, one
has
\roster\item If $(E,Q)$ is stable, then $F$ is stable.
\item If $F$ is
a stable vector bundle, then $(E,Q)$ is
stable.\endroster
\endproclaim
\demo{Proof}  See the proof of Proposition
4.7 of [NR], or [B2].\enddemo

Given a family $\SE_T$ of GPS parametrised
by $T$, we can define a family
$\SF_T$ of sheaves on $X$ by the exact
sequence
$$0\to \SF_T\to (\pi\times I_T)_*\SE_T\to\, _{x_0}\SQ_T\to
0.\tag2.2$$
Since $\SE_T$ is flat on $T$ and $\SQ_T$ locally free on $T$,
$\SF_T$
is flat on $T$, namely a flat family, thus we have

$$\hat\phi:T\to\SR.\tag2.3$$
If $\SE_T$ is a semistable family,
we get a
morphism $\phi_T: T@>\hat\phi>>\SR^{ss}@>\psi>> \Cal U_X$ by
Lemma 2.1.
Thus take
$T=\wt{\SR}^{\prime ss}$, the morphism $\phi_{\wt{\SR}^{\prime
ss}}$
induces a morphism $$\phi:\SP\to\SU_X.$$
We use the notation in $\S
1$, and let 
$\wt \SR=\underset {x\in I}\to{\times_{\bold{\widetilde
Q}}}
Flag_{\vec n(x)}(\Cal F_x).$ From now on, we will understand that
sheaves
in $\wt{\bold Q}$ have torsions only at $\{x_1,x_2\}$, so $\wt\SR
\to
\wt{\bold Q}$ is a flag bundle. Let
$$\rho:{\wt\SR}'\to\wt\SR$$ be the
natural projection, then we defined in
$\S 1$ that $$\Theta_{\wt\SR}
=(det\,R\pi_{\wt\SR}\Cal F)^k\otimes\bigotimes_{x\in I}\lbrace
(det\,\Cal
F_x)^{\alpha_x}\otimes\bigotimes^{l_x}_{i=1}
(det\,\SQ_{x,i})^{d_i(x)}\rbrace\otimes 
(det\,\Cal
F_y)^{\tilde\ell},$$
where $\tilde\ell=\ell+k$. Let $(\SE,\SQ)$ be the
universal family of 
GPS on $\wt\SR'$, we define on $\wt\SR^{\prime ss}$
that
$$\hat\Theta'=\rho^*\Theta_{\wt\SR}\otimes(det\,\SQ)^k\otimes
(det\,\SE_y)^
{-k}.$$
It is easy to check that $\hat\Theta'$ is the (restriction of)
ample
line bundle on $\wt\SR^{\prime ss}$ used to linearise the action
of
$SL(\tilde n)$ (note that $\SE$ is the pullback of $\SF$ by $\rho$),
and
descends to an ample line bundle $\Theta_{\SP}$ on $\SP$.

\proclaim{Lemma
2.3} Let $\eta_x:=(det\,\SQ)(det\,\SE_x)^{-1}$ for a point
$x\in\wt X$ and
denote $\Theta_{\wt\SR}$ by $\hat\Theta$. Then
\roster\item
$\hat\Theta'=\rho^*\hat\Theta\otimes\eta^k_y$
\item
$\Theta_{\SP}=\phi^*\Theta_{\SU_X}.$\endroster\endproclaim

\demo{Proof}
(1) is the definition of $\hat\Theta'$. To check (2),
consider the morphism
$\phi_{\wt\SR^{\prime ss}}:\wt\SR^{\prime ss}\to
\SU_X,$ it is enough to
prove that
$$\phi_{\wt\SR^{\prime
ss}}^*(\Theta_{\SU_X})=\hat\Theta'.$$
From the exact sequence (2.2), we
have
$$detR\pi_T\SF_T=(detR\pi_T(\pi_*\SE_T))\otimes(det\SQ_T)=(detR\pi_T\SE_T)
\otimes (det\SQ_T),$$
which and Theorem 1.2 in $\S 1$ imply
(2).\enddemo

\proclaim{Notation 2.1} Define $\SH$ to be the subscheme of
$\wt\SR'$
parametrising the generalised parabolic sheaves $(\SO^{\tilde
n}\to E\to
0,Q)$ satisfying
\roster\item $\Bbb C^{\tilde n}\cong H^0(E),$
and $H^1(E(-x_1-x_2-x))=0$
for any $x\in\wt X$
\item $\text{Tor}E$ is
supported on $\{x_1,x_2\}$ and
$(\text{Tor}E)_{x_1}\oplus(\text{Tor}E)_{x_2}\hookrightarrow
Q.$
\endroster\endproclaim

\proclaim{Notation 2.2} Define $\wt{\bold Q}_F$
to be the open subscheme
of $\bold{\wt Q}$
consisting of locally free
quotients $(\SO^{\tilde n}\to E\to 0)$ such
that
\roster\item $\Bbb
C^{\tilde n}\to H^0(E)$ is an isomorphism, and
\item $H^1(E(-x_1-x_2-x))=0$
for any $x\in\wt X.$
\endroster\endproclaim

\remark{Remark 2.1} The
assumption $H^1(E(-x_1-x_2-x))=0$ implies
that $H^1(E)=0$, $E$ is generated
by global sections, $H^0(E)\to E_{x_1}
\oplus E_{x_2}$ is onto, and
$E(-x_1-x_2)$ is generated by global
sections.
It is not difficult to see
([NR]) that $\SH$ is an irreducible open
subscheme of $\wt\SR'$ and
$$\wt\SR^{\prime
ss}\overset\text{open}\to\hookrightarrow
\SH\overset\text{open}\to\hookrightarrow
\wt\SR'.$$\endremark

\proclaim{Notation 2.3} Let $\wt\SR_F=\underset
{x\in I}\to
{\times_{\bold{\widetilde Q}_F}}Flag_{\vec n(x)}(\Cal F_x)$
and
$\wt\SR'_F={\rho}^{-1}(\wt\SR_F),$ then
$$\rho:\wt\SR'_F\to\wt\SR_F$$
is a grassmannian bundle over $\wt\SR_F$,
and $\wt\SR'_F\subset\SH.$ We
define
$$R^1_{F,a}:=\{(E,Q)\in\wt\SR'_F|\text{$E_{x_1}\to Q$ has rank
$a$}\},$$
and $\hat\SD_{F,1}(i):=R^1_{F,0}\cup\cdots\cup R^1_{F,i},$ which
have the
natural scheme structures. The subschemes $R^2_{F,a}$
and
$\hat\SD_{F,2}(i)$ are defined similarly. Let $\hat\SD_1(i)$ and
$\hat\SD_2(i)$ be the zariski closure of $\hat\SD_{F,1}(i)$
and
$\hat\SD_{F,2}(i)$ in $\wt\SR'$. Then they are reduced, irreducible
and
$SL(\tilde n)$-invariant closed subschemes of $\wt\SR'$, thus
inducing
closed subschemes $\SD_1(i)$, $\SD_2(i)$ of $\SP$. Clearly, we
have 
(for $j=1,2$)
that
$$\wt\SR'\supset\hat\SD_j(r-1)\supset\hat\SD_j(r-2)\cdots\supset
\SD_j(1)\supset\hat\SD_j(0)$$
$$\SP\supset\SD_j(r-1)\supset\SD_j(r-2)\supset\cdots\SD_j(1)
\supset\SD_j(0).$$
\endproclaim

\proclaim{Notation 2.4} Let
$\SR_a=\{F\in\SR|\, F\otimes\hat\SO_{x_0}=
\hat\SO_{x_0}^{\oplus a}\oplus
m_{x_0}^{\oplus(r-a)}\}$,
and
$$\hat\SW_i=\SR_0\cup\SR_1\cup\cdots\cup\SR_i$$
(which are closed in
$\SR$) endowed with their reduced scheme structures.
The subschemes
$\hat\SW_i$ are $SL(n)$-invariant, and yield closed reduced
subschemes of
$\SU_X$. It is clear
that
$$\SR\supset\hat\SW_{r-1}\supset\hat\SW_{r-2}\supset\cdots\supset\hat\SW_1
\supset\hat\SW_0=\SR_0$$
$$\SU_X\supset\SW_{r-1}\supset\SW_{r-2}\supset\cdots\supset\SW_1
\supset\SW_0.$$\endproclaim

\proclaim{Lemma 2.4} With the above
notation and $\hat\phi:\wt\SR'_F\to
\wt\SR$ defined as (2.3), we
have
\roster\item $\hat\phi(R^1_{F,a}\cap R^2_{F,b})=\SR_{a+b-r}$
\item
$\hat\phi(R^1_{F,a})=\hat\phi(R^2_{F,a})=\hat\SW_a$
\item
$\hat\phi(\hat\SD_{F,1}(i))=\hat\phi(\hat\SD_{F,2}(i))=\hat\SW_i.$
\endroster\endproclaim

\demo{Proof}  (1) follows from the Proposition 4.2 and Proposition 4.7 (1)
of [B2]. To check (2), we note that (1) implies that
$R^1_{F,a}\cap R^2_{F,j}=\emptyset$ if $j<r-a.$ Thus
$$\hat\phi(R^1_{F,a})=\bigcup^r_{j=r-a}
\SR_{a-(r-j)}=\SR_0\cup\SR_1\cup\cdots\cup\SR_a=\hat\SW_a.$$
(3) follows (2)
immediately.\enddemo

\proclaim{Proposition 2.1} With the above notation
and denoting
$\SD_1(r-1)$,
$\SD_2(r-1)$, $\SW_{r-1}$ by $\SD_1$, $\SD_2$
and $\SW$, we have
\roster\item $\phi: \SP\to\SU_X$ is finite and
surjective, and
$\phi(\SD_1(a))=\phi(\SD_2(a))=\SW_a$,
\item $\phi(\SP\ssm
\{\SD_1\cup\SD_2\})=\SU_X\ssm\SW$, and induces an
isomorphism on
$\SP\ssm\{\SD_1\cup\SD_2\},$
\item $\phi |_{\SD_1(a)}:\SD_1(a)\to\SW_a$ is
finite and surjective,
\item
$\phi(\SD_1(a)\ssm\{\SD_1(a)\cap\SD_2\cup\SD_1(a-1)\})
=\SW_a\ssm\SW_{a-1}$, and
induces an isomorphism on
$\SD_1(a)\ssm\{
\SD_1(a)\cap\SD_2\cup\SD_1(a-1)\},$
\item
$\phi:\SP\to\SU_X$ is the normalisation of $\SU_X$,
\item $\phi|_{\SD_1(a)}:\SD_1(a)\to\SW_a$ is the normalisation
of $\SW_a$,
\item $\phi(\SD_1(a)\cap\SD_2)=\SW_{a-1}$, and $\SW_{a-1}$ is the non-normal
locus of $\SW_a$.
\endroster\endproclaim

\demo{Proof}  $(5)$ and $(6)$ are
corollaries of the above (1)-(4),
and Proposition 3.2 of $\S 3$.

$(1)$ and $(3)$ follow Lemmas (2.1)-(2.4). In fact, the surjectivity
follows from
Lemma 2.1 and
Lemma 2.4, the finiteness follows from Lemma 2.3 and the
ampleness of
$\Theta_{\SU_X}$ and $\Theta_{\SP}$. 

To prove (2) and (4),
we need the following Lemmas (2.5)-(2.7). We 
will check (4) here, (2)
follows similarly. 
For
any
$\tilde\psi'(E,Q)\in\SD_1(a)\ssm\{\SD_1(a)\cap\SD_2\cup\SD_1(a-1)\}$,
we
can assume that $E$ is a vector bundle by Lemma 2.5, and $E_{x_2}\to
Q$
is an isomorphism since $\tilde\psi'(E,Q)\notin \SD_2$. Thus
$\phi(E,Q)
\in \SW_a\ssm\SW_{a-1}$ by Lemma 2.6 and Lemma 2.4, so $\phi$
induces
a morphism
$$\phi:\SD_1(a)\ssm\{\SD_1(a)\cap\SD_2\cup\SD_1(a-1)\}
@>>>\SW_a\ssm\SW_{a-1}$$
whose surjectivity follows from Lemma 2.1(3) and
$\phi(\SD_1(a)\cap\SD_2\cup\SD_1(a-1))=\SW_{a-1}$ by Lemma 2.4(1).
Now
take $T=\SR_a$ in Lemma 2.7 and use the universal property of
$\wt\SR'$, we
get a section of $\phi$ on $\SW_a\ssm\SW_{a-1}$, which
proves (4).

(7) is easy to see. In fact, $\phi(\SD_1(a)\cap\SD_2)=\SW_{a-1}$ is clear by
Lemma 2.4(1), and the non-normal locus of $\SW_a$ is contained in $\SW_{a-1}$
by the above (4). On the other hand, $\SW_{a-1}$ is irreducible since 
$\SD_j(a-1)$ is so (Proposition 3.2) and $$\phi(\SD_j(a-1))=\SW_{a-1}.$$ Thus it
suffices to prove that $\phi_a:=\phi|_{\SD_1(a)}:\SD_1(a)\to\SW_a$ is not
an isomorphism unless $\SW_{a-1}$ is empty. If $\SW_{a-1}$ is not empty, so are
$\SD_1(a)\cap\SD_2$ and $\SD_1(a-1)$. One sees easily that the fibre of $\phi_a$
at the generic point of $\SW_{a-1}$ contains at least two points since 
$\SD_1(a-1)\nsubseteq\SD_1(a)\cap\SD_2$ clearly. Therefore $\phi_a$ is not an 
isomorphism at the generic point of $\SW_{a-1}$.
\enddemo

\proclaim{Lemma 2.5} Every semistable GPS $(E',Q')$ with
$rank(E')>0$ is
$s$-equivalent to a semistable GPS $(E,Q)$ with $E$
locally
free.\endproclaim

\demo{Proof}  For given $(E',Q')\in\Cal C_{\mu}$
with $rank(E')>0$,
if $\text{Tor}(E')=0,$ then we are done. Thus we assume
that one of 
$\text{Tor}(E')_{x_i}$, say $\text{Tor}(E')_{x_1}$, is
nontrivial and 
the lemma is true for all $(\wt E',\wt Q')\in \Cal C_{\mu}$
with
$dim(\text{Tor}(\wt E'))<dim(\text{Tor}( E'))$.

Let $q'_i:E'_{x_i}\to
Q'$ be the maps induced by the generalised parabolic
structure of
$(E',Q')$, and choose a projection $p:Q'\to\Bbb C$ such
that
$p(q'_1(\text{Tor}(E')_{x_1}))\neq 0$. Let $E'\to\,_{x_1}\Bbb C\to 0$
be
the morphism 
$$E'@>>>E'_{x_1}@>q'_1>>Q'@>p>>\Bbb C$$
and $\wt E'$ its
kernel, which has a smaller torsion than $E'$ by the
choice of $p$. Then we
have an exact sequence
$$0\to\wt E'\to E'\to\,_{x_1}\Bbb C\to 0$$
which
induces an exact sequence of GPS (See Definition 1.6) if we
set
$\tau=\,_{x_1}\Bbb C$
$$0\to (\wt E',\wt Q')\to(E',Q')\to(\tau,\Bbb
C)\to 0.$$
One can check that $(\wt E',\wt Q')\in\Cal C_{\mu}$, thus there
exists
a $(\wt E,\wt Q)\in\Cal C_{\mu}$ with $\wt E$ locally free such
that
$gr(\wt E,\wt Q)=gr(\wt E',\wt Q')$. Since $(\tau,\Bbb C)$ is stable,
we
have $$gr(E',Q')=gr(\wt E,\wt Q)\oplus(\tau,\Bbb C).$$
Let $\tilde q:\wt
E_{x_1}\oplus\wt E_{x_2}\to\wt Q\to 0$ and
$K_1=ker(\tilde q_1:\wt
E_{x_1}\to\wt Q)$, choosing a Hecke modification
$h:\wt E\to E$ at $x_1$
(See Remark 1.4 of [NS]) such that
$\wt K_1:=ker(h_{x_1})\subset K_1$ and
$dim(\wt K_1)=1$, we get the
extension
$$0@>>>\wt
E@>h>>E@>\gamma>>\tau@>>>0.$$
Let $Q=\wt Q\oplus\Bbb C$ and
$E_{x_1}=h_{x_1}(\wt E_{x_1})\oplus V_1$
for a subspace $V_1$. We define a
morphism $f:E_{x_1}\oplus E_{x_2}\to Q$
such that $E_{x_2}\to Q$ to be
$$E_{x_2}@>h^{-1}_{x_2}>>\wt E_{x_2}@>\tilde q_2>>\wt Q\hookrightarrow
Q$$
and $E_{x_1}\to Q$ to be
$$h_{x_1}(\wt E_{x_1})\oplus V_1@>({\bar
h}^{-1}_{x_1},\gamma_{x_1})>>
\frac{\wt E_{x_1}}{\wt K_1}\oplus\Bbb
C@>(\tilde q_1,id)>>\wt Q\oplus\Bbb
C=Q$$
where $\bar h_{x_1}:\wt
E_{x_1}/\wt K_1\cong h_{x_1}(\wt E_{x_1})$ and
$\tilde q_1:\wt E_{x_1}/\wt
K_1\to\wt Q$ (Note that $\wt K_1\subset K_1$).
Thus the following diagram
is commutative
$$\CD
@.\wt E_{x_1}\oplus\wt
E_{x_2}@>(h_{x_1},h_{x_2})>>E_{x_1}\oplus E_{x_2}
@>(\gamma_{x_1},0)>>\Bbb
C@>>>0\\
@.      @VqVV       @VfVV               @|           @.\\
0
@>>>\wt Q @>>> \wt Q\oplus\Bbb C  @>>> \Bbb C @>>>  0
\endCD$$
One checks
that $f$ is surjective by this diagram, and thus
$$0\to(\wt E,\wt
Q)\to(E,Q)\to(\tau,\Bbb C)\to 0$$
It is easy to see that $(E,Q)\in\Cal
C_{\mu}$ and $s$-equivalent to
$(E',Q')$.\enddemo

It is well known that
for any $F\in \SR^{ss}$ there is an integer $a_F$
such
that
$F\otimes\hat\SO_{x_0}\cong \hat\SO_{x_0}^{\oplus
a_F}\oplus
m_{x_0}^{\oplus(rk(F)-a_F)},$ 
thus defining a function
$$\bold
a:\SR^{ss}\to\Bbb Z_{\geq 0}$$
with $\bold a (F)=a_F$. It is not clear that
if this function induces
a function $$\bold a:\SU_X\to\Bbb Z_{\geq
0}$$
(See the $\lq$Remarque' on page 172 of [S2]). However, the
following
lemma implies that $\bold a$ is invariant under $s$-equivalence,
in
particular, descends to $\SU_X$.

\proclaim{Lemma 2.6} Let $0\to F_1\to
F\to F_2\to 0$ be an exact
sequence of torsion free sheaves. Then
$$\bold
a(F)=\bold a(F_1)+\bold a(F_2).$$
In particular, if $F$ is $s$-equivalent
to $F'$, then 
$\bold a(F)=\bold a(F')$.\endproclaim

\demo{Proof} For any
torsion free sheaf $F_i$, we define a vector
bundle on $\wt X$ to be
$E_i=\pi^*F_i/\text{Tor}(\pi^*F_i).$ Note that
we have a diagram
$$\CD
0
@>>> F_1         @>>> F       @>>> F_2     @>>> 0\\
@.      @VVV
@VVV         @VVV         @.  \\
@.
\pi_*\pi^*F_1@>>>\pi_*\pi^*F@>>>\pi_*\pi^*F_2@>>>0
\endCD$$  
which induces
$$\CD
@.      0  @.           0   @.        0  @.           \\
@.
@VVV            @VVV          @VVV         @. \\
0  @>>> F_1         @>>> F
@>>> F_2     @>>> 0  \\
@.      @VVV             @VVV         @VVV
@.  \\
0 @>>>\pi_*E_1  @>>>   \pi_*E@>>>\pi_*E_2      @>>>0   \\
@.
@VVV             @VVV          @VVV        @.
\\
0@>>>\,_{x_0}Q_1@>>>\,_{x_0}Q@>>>\,_{x_0}Q_2@>>> 0    \\
@.      @VVV
@VVV         @VVV        @.\\
@.      0    @.           0   @.        0
@.
\endCD$$  
where $Q_1$, $Q$, $Q_2$ are defined such that each vertical
sequence
is exact, the third horizontal sequence is defined such that
the
diagram is commutative, which must be exact. 

For a torsion free sheaf
$F$, if we define $Q_F$ by
$$0\to F\to\pi_*\pi^*F\to\,_{x_0}Q_F\to 0$$
then
one can see that 
$Q_i=Q_{F_i}/\pi_*\text{Tor}(\pi^*F_i)$ and
$Q=Q_F/\pi_*\text{Tor}(\pi^*F)$ in the above diagram (See the proof
of
Lemma 2.1). Thus $dim(Q_i)=\bold a(F_i)$ and $dim(Q)=\bold a(F)$,
which
proves the lemma.\enddemo

\proclaim{Lemma 2.7} Let $T$ be a reduced
scheme, $\SF$ a sheaf on
$X\times T$, flat over $T$, such that for $t\in T$
the sheaf $\SF_t$
on $X$ is torsion free of rank $r$ and $\bold a(\SF_t)=a$
is constant.
Then there exists a vector bundle $\SE$ of rank $r$ on $\wt
X\times T$ and
a locally free rank $r$ quotient
$q:\SE_{x_1}\oplus\SE_{x_2}\to\SQ\to 0$
on $T$ such that
$q_2:\SE_{x_2}\to\SQ$ is an isomorphism,
$q_1:\SE_{x_1}\to\SQ$ has rank $a$
at each fibre and
$$0\to\SF\to (\pi\times I)_*\SE\to\,_{x_0\times T}\SQ\to
0.$$
\endproclaim

\demo{Proof} We can assume that $\SF$ is torsion free
(See Lemma 4.13 of
[NR]), and define $\SQ_{\SF}$
by
$$0\to\SF\to\pi_*\pi^*\SF\to\,_{x_0\times T}\SQ_{\SF}\to 0\tag2.4$$
(we
will write $\pi$ for $\pi\times I$), which gives for any $t\in T$
an exact
sequence
$$0\to \SF_t\to\pi_*\pi^*\SF_t\to\,_{x_0}(\SQ_{\SF})_t\to
0$$
since $\SF_t$ is torsion free. This shows that
$(\SQ_{\SF})_t\cong\SQ_{\SF_t}$ has
constant dimension $2r-a$, hence
$\SQ_{\SF}$ is flat (in fact, a vector
bundle) on $T$,
which and (2.4)
imply that $\pi^*\SF$ is flat on $T$.

Take a resolution $0\to\Cal
K\to\SL^{-m}\to\pi^*\SF\to 0$ of $\pi^*\SF$,
where $\SL$ is a line bundle,
and dualise it, we get
$$0\to (\pi^*\SF)^{\vee}\to\SL^{m}\to\Cal
K^{\vee}\to\Cal Ext^1(\pi^*\SF,
\SO_{\wt X\times T})\to 0.$$
Note that
$\Cal K$ is a vector bundle on $\wt X\times T$ (Since $\pi^*\SF$
is flat
over $T$ and $\wt X$ is smooth) and $(\Cal K^{\vee})_t\cong
(\Cal
K_t)^{\vee}$, we have the following
diagram
$$\CD
@.((\pi^*\SF)^{\vee})_t@>>>(\SL^m)_t
@>>>(\Cal
K^{\vee})_t@>>>\Cal Ext^1(\pi^*\SF,\SO_{\wt X\times
T})_t@>>>0\\
@.      @VVV       @|               @|       @VVV
@.\\
@.(\pi^*\SF_t)^{\vee}@>>>(\SL_t)^m @>>>(\Cal
K_t)^{\vee}@>>>\Cal
Ext^1(\pi^*\SF_t,\SO_{\wt X})@>>>0
\endCD$$
which implies that $\Cal
Ext^1(\pi^*\SF,\SO_{\wt X\times T})_t
\cong\Cal Ext^1(\pi^*\SF_t,\SO_{\wt
X})$, 
thus $\Cal Ext^1(\pi^*\SF,\SO_{\wt X\times T})$ is flat over $T$.
This
shows that $(\pi^*\SF)^{\vee}$ is locally free and
$((\pi^*\SF)^{\vee})_t\cong(\pi^*\SF_t)^{\vee}$.
Let
$\wt\SE:=(\pi^*\SF)^{\vee\vee}$ be the double dual of $\pi^*\SF$, then
we 
have
$$0\to\Cal T\to\pi^*\SF\to\wt\SE\to 0\tag2.5$$
which specialises
for any $t\in T$ to 
$$0\to\text{Tor}(\pi^*\SF_t)\to\pi^*\SF_t\to
\wt\SE_t\to 0.$$
By (2.4) and (2.5), we
get
$$0\to\SF\to\pi_*\wt\SE\to\,_{x_0\times T}\wt\SQ\to 0,$$
and $\SQ$ is a
vector bundle of rank $a$ on $T$. Now the same 
construction in the proof
of Lemma 2.1 proves our lemma.
\enddemo

\proclaim{Lemma 2.8} Let $E'$ be a rank $r$ (stable) semistable parabolic
bundle of degree $d-r$ on $\wt X$. Then its direct image $F=\pi_*E'$ is
a (stable) semistable parabolic sheaf of degree $d$ on $X$, such that
$F\otimes\hat\SO_{x_0}\cong m_{x_0}^{\oplus r}$. This construction gives
a morphism $$\SU_{\wt X}(d-r)\to\SW_0.$$\endproclaim

\demo{Proof} The proof of Lemma 2.1 clearly shows that $E'\to F=\pi_*E'$
and $F\to E'=\pi^*F/\text{Tor($\pi^*F$)}$ gives a bijection between the set
of isomorphism classes of rank $r$ bundles $E'$ of degree $d-r$ on $\wt X$
and the set of torsion free sheaves $F$ of degree $d$ on $X$ with $F\otimes
\hat\SO_{x_0}\cong m_{x_0}^{\oplus r}.$

We check now that the (stablity) semistablity of $E'$ implies that of $F$.
For any subsheaf $F_1\subset F$ of rank $r_1$ such that $F/F_1=F_2$ is a torsion
free sheaf of rank $r_2$, then $\bold a(F_1)=\bold a(F_2)=0$ since $\bold a(F)=0$
and $\bold a(F)=\bold a(F_1)+\bold a(F_2)$ by Lemma 2.6.
Thus we have $\pi_*E'_i=F_i$, where $E'_i=\pi^*F_i/\text{Tor($\pi^*F_i$)}$, 
and an exact sequence
$0\to E'_1\to E'\to E_2'\to 0$. One computes that 
$$\frac{pardeg(E'_i)}{rank(E_i')}-\frac{pardeg(E')}{rank(E')}=\frac{pardeg(F_i)}
{rank(F_i)}-\frac{pardeg(F)}{rank(F)},$$
which proves that $E'$ is (stable) semistable if and only if $F=\pi_*E'$ is so.
Since $\SW_0$ has reduced scheme structure, the above construction gives a 
morphism $$\SU_{\wt X}(d-r)\to\SW_0.$$
\enddemo

\proclaim{Corollary 2.1} Suppose that $g>1$. Then $\SW_0$ is nonempty, and 
contains stable parabolic sheaves if $|I|>0.$  In particular, for any $0\le 
a<r$, $\SW_a\neq \SW_{a+1}.$\endproclaim

\demo{Proof} We will prove in \S5 that $codim(\wt\SR_F\ssm\wt\SR^{ss})\ge
(r-1)(\tilde g-1)+1$ and $codim(\wt\SR^{ss}\ssm\wt\SR^s)\ge (r-1)(\tilde g-1)+1$
if $|I|>0$ (See Proposition 5.1). Thus $\SU_{\wt X}(d-r)$ is nonempty if $\tilde g=g-1>0$, 
and there exist stable parabolic bundles of degree $d-r$ on $\wt X$ if
moreover $|I|>0.$ 
Now using Lemma 2.8, we conclude that $\SW_0$
is nonempty, and contains stable parabolic sheaves if $|I|>0.$

Since semistability is an open condition and $\hat\SW_{a+1}\ssm\hat\SW_a$ is a 
nonempty open set of $\hat\SW_{a+1}$, there is a semistable sheaf 
$F\in\hat\SW_{a+1}\ssm\hat\SW_a$ (Because we have shown that 
$\hat\SW_{a+1}^{ss}$ is nonempty), and we can see that 
$\psi(F)\in\SW_{a+1}\ssm\SW_a$ by Lemma 2.6.
\enddemo

\remark{Remark 2.2} When $X$ is a nodal curve of $g=1$ and $|I|=0$, 
it is possible that $\SW_0$ is empty (I am not saying that every $\SW_a$ is empty).
In fact, if $\SW_0$ is nonempty in this case, then there exists a semistable
bundle of degree $d-r$ on $\wt X=\Bbb P^1$ by Lemma 2.8, which implies that
$r|d$.\endremark

We will finish this section by computing the canonical sheaf of $\wt\SR_F$.
Recall that $\pi:\wt\SR_F=\underset{x\in I}\to{\times_{\bold{\widetilde Q}_F}}
Flag_{\vec n(x)}(\Cal F_x)\to \bold{\wt Q}_F$, let $\SE=\pi^*\SF$ and
$$\SE_x=F_0(\SE)_x\supset F_1(\SE)_x\supset\cdots
\supset F_{l_x}(\SE)_x\supset F_{l_x+1}(\SE)_x=0$$
the universal flag on $\wt\SR_F$. Write $\SQ_{x,i}=\SE_x/F_i(\SE)_x$ and 
$\pi_{\wt\SR_F}:\wt X\times\wt\SR_F\to\wt\SR_F$ the projection. Then we have

\proclaim{Proposition 2.2} Let $\omega_{\wt\SR_F}$ be the canonical sheaf of
$\wt\SR_F$, and $\omega_{\wt X}=\SO_{\wt X}(\sum_qq)$ the canonical sheaf of
$\wt X$. Then
$$\aligned\omega^{-1}_{\wt\SR_F}=&(det\,R\pi_{\wt\SR_F}\SE)^{2r}\otimes
\bigotimes_{x\in I}\left\{(det\,\SE_x)^{n_{l_x+1}-r}
\otimes\bigotimes^{l_x}_{i=1}(det\,\SQ_{x,i})
^{n_i(x)+n_{i+1}(x)}\right\}\\
&\otimes\bigotimes_q(det\,\SE_q)^{1-r}\otimes(det\,R\pi_{\wt\SR_F}det\SE)^{-2}.
\endaligned$$
\endproclaim
 
\demo{Proof}  Note that $\omega_{\wt\SR_F}=
\omega_{\wt\SR_F/\bold{\wt Q}_F}\otimes\pi^*\omega_{\bold{\wt Q_F}}$, 
the proposition is clearly
a corollary of the following two lemmas.\enddemo

\proclaim{Lemma 2.9} Let $E$ be a vector bundle of rank $r$ on $M$, and
$F(l,E)=Flag_{\vec n}(E)$ of type $\vec n=(n_1,\cdots,n_{l+1})$, with the
universal flag 
$$E=F_0(E)\supset F_1(E)\supset\cdots\supset F_l(E)\supset F_{l+1}(E)=0$$
on $F(l,E)$ and $\SQ_i=E/F_i(E)$. Then
$$\omega_{F(l,E)/M}=(det\,E)^{r-n_{l+1}}\otimes\bigotimes^l_{i=1}(det\SQ_i)
^{-(n_i+n_{i+1})}.$$\endproclaim

\demo{Proof} One considers $F(l,E)$ as the grassmannian bundle 
$$p:Grass_{rk(F_l(E))}(F_{l-1}(E))\to F(l-1,E)$$ 
over $F(l-1, E)$. Then $\omega_{F(l,E)/F(l-1,E)}=det(F_l(E)\otimes(F_{l-1}(E)/
F_l(E))^{\vee})$, and one has
$$\omega_{F(l,E)/M}=\bigotimes^l_{i=1}det(F_i(E)\otimes(F_{i-1}(E)/F_i(E))^{\vee}).$$
Thus one can compute that
$$\omega_{F(l,E)/M}=(det\,E)^{r-n_{l+1}}\otimes\bigotimes^l_{i=1}(det\SQ_i)
^{-(n_i+n_{i+1})}.$$
\enddemo

\proclaim{Lemma 2.10} Let 
$\SO_{\wt X\times\bold{\wt Q_F}}^{\tilde n}\to\SF\to 0$
be the universal quotient on $\wt X\times\bold{\wt Q_F}$. Then
$$\omega^{-1}_{\bold{\wt Q_F}}=(det\,R\pi_{\bold{\wt Q_F}}\SF)^{2r}\otimes
\bigotimes_q(det\,\SF_q)^{1-r}\otimes(det\,R\pi_{\bold{\wt Q_F}}det\,\SF)^{-2}.$$
\endproclaim

\demo{Proof} We have, on $\wt X\times\bold{\wt Q_F}$, the exact sequence
$0\to \Cal K\to\SO^{\tilde n}\to\SF\to 0,$
the tangent space of $\bold{\wt Q_F}$ at a point $(0\to K\to\SO^{\tilde n}\to
E\to 0)$ is $H^0(\wt X,K^{\vee}\otimes E)$. From the properties of $\bold{\wt Q_F}$ 
(the Notation 2.2), it follows that
$$\omega^{-1}_{\bold{\wt Q_F}}=det\,R\pi_{\bold{\wt Q_F}}(\SF\otimes\SF^{\vee}).$$
We will now use a variant of the method of [DN] to evaluate 
$det\,R\pi_{\bold{\wt Q_F}}(\SF\otimes\SF^{\vee}).$

Let $\Cal M=\pi_{\bold{\wt Q_F}*}(\SF)$ be the direct image sheaf of $\SF$, 
which is local free of rank $\tilde n=d+r(1-\tilde g)$. Let $Gr$ be the
grassmannian of rank $r-1$ subbundles of $\Cal M$, and 
$$p:Gr\to \bold{\wt Q_F}$$
the projection. We consider the canonical exact sequence on $Gr$
$$0\to U_{Gr}\to p^*\Cal M\to Q_{Gr}\to 0,$$
where $U_{Gr}$ is the relative universal subbundle of $p^*\Cal M$ on $Gr$, and
$Q_{Gr}$ the relative universal quotent. Let $\SO_{Gr}(-1)=det(U_{Gr})$ and consider
$$\CD
\wt X\times Gr  @>\pi_{Gr}>>     Gr          \\
@V1\times p VV                @VpVV          \\
\wt X\times\bold{\wt Q_F} @>\pi_{\bold{\wt Q_F}}>> \bold{\wt Q_F},
\endCD$$
we have the induced morphism
$$\pi^*_{Gr}U_{Gr}\hookrightarrow\pi^*_{Gr}p^*\Cal M=
(1\times p)^*\pi^*_{\bold{\wt Q_F}}\pi_{\bold{\wt Q_F}*}(\SF)
\to (1\times p)^*\SF.\tag2.6$$
Let $Gr_0$ be the open set of $Gr$ such that, on $\wt X\times Gr_0$, the above 
induced morphism $$\pi^*_{Gr}U_{Gr}\to (1\times p)^*\SF$$
in (2.6) is injective. Then, if we write $D=Gr\ssm Gr_0$, there is on 
$\wt X\times(Gr\ssm D)$ an exact sequence
$$ 0\to \pi^*_{Gr}U_{Gr}\to (1\times p)^*\SF\to (1\times p)^*det(\SF)
\otimes\pi^*_{Gr}\SO_{Gr}(1)\to 0.\tag2.7$$
We will denote the morphisms $Gr\ssm D\to \bold{\wt Q_F}$
and $\wt X\times(Gr\ssm D)\to Gr\ssm D$ by the same $p$ and $\pi_{Gr}$. 
By using (2.7), we compute that
$$\aligned det\,R\pi_{Gr}(1\times p)^*(\SF\otimes\SF^{\vee})
=&(det\,R\pi_{Gr}(1\times p)^*\SF^{\vee})^{r-1}\otimes\SO_{Gr}(-rd)\\
&\otimes det\,R\pi_{Gr}(1\times p)^*(\SF^{\vee}\otimes det\,\SF),
\endaligned\tag2.8$$
$$det\,R\pi_{Gr}(1\times p)^*\SF=det\,R\pi_{Gr}(1\times p)^*det(\SF)
\otimes\SO_{Gr}(-d).\tag2.9$$
Use $$0\to\pi^*_{Gr}\SO_{Gr}(-1)\to (1\times p)^*(\SF^{\vee}\otimes det\,\SF)\to
(1\times p)^*\otimes\pi^*_{Gr}U^{\vee}_{Gr}\to 0,$$ 
we get
$$det\,R\pi_{Gr}(1\times p)^*(\SF^{\vee}\otimes det\,\SF)= 
(det\,R\pi_{Gr}(1\times p)^*det\,\SF)^{r-1}\otimes \SO_{Gr}(-d).\tag2.10$$
Thus, by (2.8)--(2.10) and the base change theorem, we have
$$ \aligned p^*det\,R\pi_{\bold{\wt Q_F}}(\SF\otimes\SF^{\vee})
=p^*\{&(det\,R\pi_{\bold{\wt Q_F}}\SF)^{r+1}
\otimes(det\,R\pi_{\bold{\wt Q_F}}\SF^{\vee})^{r-1}\\
&\otimes(det\,R\pi_{\bold{\wt Q_F}}det\,\SF)^{-2}\}\endaligned.$$
By duality and the exact sequence
$$0\to\SF\to\SF\otimes\omega_{\wt X\times\bold{\wt Q_F}/\bold{\wt Q_F}}\to
\bigoplus_q\SF_q\to 0,$$
one has that
$$det\,R\pi_{\bold{\wt Q_F}}\SF^{\vee}=det\,R\pi_{\bold{\wt Q_F}}
(\SF\otimes\omega_{\wt X\times\bold{\wt Q_F}/\bold{\wt Q_F}})=
(det\,R\pi_{\bold{\wt Q_F}}\SF)\otimes\bigotimes_q(det\SF_q)^{-1}.$$
Thus the lemma follows if $p^*:Pic(\bold{\wt Q_F})\to Pic(Gr\ssm D)$
is injective, which will be proved in the next lemma.
\enddemo

\proclaim{Lemma 2.11} $p^*:Pic(\bold{\wt Q_F})\to Pic(Gr\ssm D)$ is injective.
\endproclaim

\demo{Proof} It is well known that 
$Pic(Gr)=Pic(\bold{\wt Q_F})\oplus\, \Bbb Z\SO_{Gr}(1)$. 
For each fibre $Gr(E)=p^{-1}(E)$ of $p:Gr\to \bold{\wt Q_F}$,
$D\cap Gr(E)$ is an irreducible hypersurface of $Gr(E)$ 
(See Lemma 7.3 of [DN]). Thus the ideal sheaf
$$\SO_{Gr}(-D)=p^*(\Lambda)\otimes\SO_{Gr}(a_0),\quad\text{for some $\Lambda\in
Pic(\bold{\wt Q_F})$},$$
with $a_0\neq 0.$ One has the exact sequence 
(See Charpter II, Proposition 6.5, of [Ha])
$$\Bbb Z@>i>> Pic(Gr)@>>> Pic(Gr\ssm D)@>>>0,$$
where $i(1)=\SO_{Gr}(-D).$ For any $\SL\in Pic(\bold{\wt Q_F})$, if
$p^*\SL|_{Gr\ssm D}$ is trivial, then there exists $m\in\Bbb Z$ such that
$$p^*\SL=i(m)=p^*(\Lambda^m)\otimes\SO_{Gr}(ma_0).$$
The $m$ has to be zero, namely, $p^*\SL=\SO_{Gr}$, which implies clearly that
$\SL$ is trivial. \enddemo

\heading \S 3 Geometry of moduli spaces of generalized
parabolic sheaves\endheading

We will prove in this section that $\SP$ and
its subvarieties $\SD_j(a)$
and $\SD_1(a)\cap\SD_2(b)$ are reduced,
normal with rational singularities. In
particular, we will
prove that $\SH$, $\hat\SD_j(a)$ and $\hat\SD_1(a)\cap\hat\SD_2(b)$
are
reduced, normal with rational singularities, and prove a formula to express
the canonical (dualizing) sheaf of $\SH.$ 
We will use the following device to analyse singularities of a variety $V$:
Find varieties $W$ and $V'$ and smooth morphisms $f:V\to W$ and
$f':V'\to W$, such that the singularities of $V'$ are easy to analyse. We will
call $V'$ (or its complete local ring at a point) the smooth model of $V$ (or
local smooth model at a point of $V$). For
simplicity, we assume that $|I|=0$, which will not affect
the generality. 
Let $Y$ be a scheme of finite type, $\SF$ a locally free $\SO_Y$-module
of rank $r$ and 
$$\Bbb Hom(\SO^r_Y,\SF):=Spec\,S(\Cal Hom(\SO^r_Y,\SF)^{\vee})\to Y,$$
which parametrizes homomorphisms from $\SO^r_Y$ to $\SF$. Let 
$$F_Y:=Isom(\SO^r_Y,\SF)\subset\Bbb Hom(\SO^r_Y,\SF)$$
be the open subscheme corresponding to isomorphisms. Then we call $F_Y\to Y$
the frame bundle associated to $\SF$. Moreover, if $\SE$ is a $\SO_Y$-module, 
then the functor
$$\Bbb Hom(\SE,\SF):T\longmapsto Hom_{\SO_T}(\SE_T,\SF_T)$$
from the category of $Y$-schemes to the category of sets is representable by
$$V(\SF^{\vee}\otimes_{\SO_Y}\SE):=Spec\,S(\SF^{\vee}\otimes_{\SO_Y}\SE)$$ 
(Proposition 9.6.1 of [EGA-I]), where
$\SE_T$ and $\SF_T$ denote the pulling backs of $\SE$ and $\SF$ to $T$.
We will also need that

\proclaim{Lemma 3.1} Let $f:W\to V$ be a smooth morphism. Then $W$ is reduced
(respectively, normal, Cohen-Macaulay, Gorenstein) with only rational 
singularities if and only if $V$ is.\endproclaim

\demo{Proof} See Proposition 4.19 of [NR] for the statment about the rational 
singularity, and the other statments are well-known commutative algebraic facts.
\enddemo  

By a point $(E,h)\in \SH,$ we mean that
$(\SO^{\tilde n}@>e>>E\to 0,E_{x_1}\oplus E_{x_2}@>h>>\Bbb C^r\to 0)$   
such that $E$ is
locally free outside $\{x_1,x_2\}$,
$(\text{Tor}E)_{x_1}
\oplus(\text{Tor}E)_{x_2}\overset
h\to\hookrightarrow\Bbb C^r,$
and $H^1(E(-x_1-x_2-x))=0$ for $x\in\wt X$,
$e$ induces isomorphism
$\Bbb C^{\tilde n}\cong H^0(E).$ Let $\bl$ be the
category of Artinian
local $\Bbb C$-algebras, and consider the
functor
$$\Phi_1:\bl\to \bset$$
defined by (write $S=Spec(A)$ for any $A\in
\bl$):
$$\Phi_1(A):=\left\{\aligned &\text{Equivalent classes of flat families $(\SE^S,h^S)$}\\ 
&\text{such that $(\SE^S,h^S)|_{Spec(A/m)\times\wt X}=(E,h)$}\endaligned 
\right\},$$
where
$(\SE^S,h^S)=(\SO_{S\times\wt X}^{\tilde n}@>e^S>>\SE^S\to 0,
\SE^S_{x_1}\oplus\SE^S_{x_2}@>h^S>>\SO_S^r\to 0).$

Given a point
$z=(E,h)\in\SH$, let $dim(h_1(E_{x_1}))=r_1$ and
$dim(h_2(E_{x_2}))=r_2$.
The map $e$ induces two maps $\Bbb C^{\tilde n}@>e_i>>E_{x_i}\to 0$. We
denote images of the canonical base of $\Bbb C^{\tilde n}$ (under $h_i$)
by 
$$(\alpha^i_{_j1},\cdots,\alpha^i_{jr})\in\Bbb C^r$$ 
where
$j=1,...,\tilde n$. Without loss of generality, we assume
that
$$rk\pmatrix
\alpha^i_{11}&\hdots&\alpha^i_{1r_i}\\
\vdots&\ddots&\vdots\\
\alpha^i_{\tilde n1}&\hdots&\alpha^i_{\tilde nr_i}
\endpmatrix=r_i$$
For
any algebra $A$, we use $A^r@>r_i>>A^{r_i}$ and $A^r@>[r-r_i]>>
A^{r-r_1}$
to denote the projections $(y_1,...,y_r)\mapsto
(y_1,...,y_{r_i})$ and
$(y_1,...,y_r)\mapsto (y_{r_i+1},...,y_r)$
respectively. Thus for any
$(\SE^S,h^S)\in\Phi_1(A)$
$$\SE^S_{x_i}@>h^S_i>>\SO^r_S@>r_i>>\SO_S^{r_i}$$
are
surjective since it is  so at the fibre then by Nakayama's lemma, 
from
which we get a surjection
$$\SE^S@>q^S>>\,_{x_1}\SO^{r_1}_S\oplus\,_{x_2}\SO_S^{r_2}\to 0.$$
Let
$\wt\SE^S=ker(q^S)$, we get
$$0\to\wt\SE^S\to\SE^S@>q^S>>\,_{x_1}\SO^{r_1}_S\oplus\,_{x_2}\SO_S^{r_2}\to 0,$$
$$\SE^S_{x_1}@>h^S_1>>\SO^r_S@>[r-r_1]>>\SO_S^{r-r_1},\quad
\SE^S_{x_2}@>h^S_2>>\SO^r_S@>[r-r_2]>>\SO_S^{r-r_2}$$
which we denote by
$(\wt\SE^S,q^S,[r-r_1]\cdot
h_1^S,[r-r_2]\cdot
h_2^S):=\varphi_A
(\SE^S,h^S)$. It is clear that the
restriction of $\varphi_A(\SE^S,h^S)$
at the fibre $Spec(A/m)\times\wt X$
is
$$0\to\wt E\to E@>q^0>>\,_{x_1}\Bbb C^{r_1}\oplus\,_{x_2}\Bbb
C^{r_2}\to 0,$$

$$E_{x_1}@>h_1>>\Bbb C^r@>[r-r_1]>>\Bbb
C^{r-r_1},\quad
E_{x_2}@>h_2>>\Bbb C^r@>[r-r_2]>>\Bbb C^{r-r_2}$$
which is
$\varphi_{\Bbb C}(E,h).$ Thus the above construction gives a
morphism
$$\varphi:\Phi_1\to\Phi_2$$ of functors, where $\Phi_2$ will be
defined later.

Let $\bwq^1$ be the Quot scheme of rank $r$, degree
$d-r_1-r_2$ quotients
$$\SO_{\wt X}^{\tilde n-r_1-r_2}\to\wt E\to 0$$ 
and
$\bwq^1_F$ the open subset of locally free quotients with vanishing
$H^1(\wt E)$ such that $\Bbb C^{\tilde n-r_1-r_2}\to H^0(\wt E)$ is an
isomorphism. It is known that $\bwq^1_F$ is smooth. Let
$f:\bwq^1_F\times\wt X\to\bwq^1_F$ be the projection, $\wt\SE$ the
universal quotient on
$\bwq^1_F\times\wt X$, then the sheaf (See [La] for
the definition)
$$\Cal G:=Ext^1_f(\,_{x_1}\SO^{r_1}\oplus\,_{x_2}\SO^{r_2},\,\wt\SE)$$
where
$\SO=\SO_{\bwq^1_F}$, is locally free. Write $V:=V(\Cal G^{\vee})\to\bwq^1_F$ and
$$p_V:V\times\wt X\to\bwq^1_F\times\wt X,$$
then
there exists a universal extension on $V\times\wt X$
$$0\to p_V^*\wt\SE\to
\SE\to\,_{x_1}\SO^{r_1}_V\oplus\,_{x_2}\SO^{r_2}_V
\to 0.$$
Let $W_i$ be
the total space of 
$\Cal Hom_{\SO_V}(\SE_{x_i},\SO_V^{r-r_i})$, namely,
the $V$-scheme
$$W_i=V(\SE_{x_i}^{\oplus (r-r_i)})\to V,$$
which represents
the functor $\Bbb Hom(\SE_{x_i},\SO_V^{r-r_i})$ (See
Proposition (9.6.1)
of [EGA-I])
and let $Y:=W_1\times_VW_i$, then the $S$-points of $Y$ can be
expressed
as
$$\left\{(\wt\SE^S,q,e^S,\xi_1,\xi_2):=\left(\aligned
&\quad 0\to\wt\SE^S\to\SE^S@>q>>\,_{x_1}\SO^{r_1}_S\oplus\,_{x_2}\SO_S^{r_2}\to
0,\\
&\SO^{\tilde n-r_1-r_2}@>e^S>>\wt\SE^S\to 0,\,\SE^S_{x_1}@>\xi_1>>
\SO_S^{r-r_1},\,\SE^S_{x_2}@>\xi_2>>\SO_S^{r-r_2}
\endaligned\right)\right\}$$
where $\SO^{\tilde n-r_1-r_2}@>e^S>>\wt\SE^S\to 0$ is induced from $\wt\SE$ by a
$S$-point of $\bwq^1_F$. Thus we can define the functor
$\Phi_2:\bl\to\bset $ as
$$\Phi_2(A)=\left\{\aligned
&(\wt\SE^S,q,\xi_1,\xi_2)\quad\text{such that
its restriction at
$Spec(A/m)\times\wt X$}\\ 
&\text{is $(\wt E,q^0,[r-r_1]\cdot
h_1,[r-r_2]\cdot
h_2)$, which is $y:=\varphi_{\Bbb
C}(E,h)$}\endaligned\right\}$$
where $\wt\SE^S$ is a flat family of bundles parametrized by $S$, which are
generated by global sections and with vanishing $H^1$. We have an obvious smooth 
morphism $Y(-)\to \Phi_2$ by forgetting the frames of ${\pi_S}_*\wt\SE^S$, 
where $Y(-)$ denotes the functor defined by the variety $Y$.

\proclaim{Lemma 3.2} The morphism
$\varphi:\Phi_1\to\Phi_2$ is formally
smooth.\endproclaim

\demo{Proof} For
any $A\to A_0\to 0$, where $A_0=A/I$, $I^2=0$,
we consider 
$$\CD
\Phi_1(A)
@>>>     \Phi_1(A_0)   \\
@V\varphi_A VV          @V\varphi_{A_0}VV
\\
\Phi_2(A)      @>>>          \Phi_2(A_0)
\endCD$$
For any given points
$(\SE^{S_0},h^{S_0})\in\Phi_1(A_0)$ 
and
$(\wt\SE^S,q,\xi_,\xi_2)\in\Phi_2(A)$ such that
$$\CD
0@>>>\wt\SE^S|_{S_0\times\wt X}@>>>\SE^S|_{S_0\times\wt
X}@>q|_{S_0\times\wt
X}>>\,_{x_1}\SO_{S_0}^{r_1}\oplus\,_{x_2}\SO^{r_2}_{S_0}@>
>> 0\\
@.      @|       @|               @|               @.
\\
0@>>>\wt\SE^{S_0}@>>>\SE^{S_0}@>q^{S_0}>>\,_{x_1}\SO^{r_1}_{S_0}\oplus\,
_{x_2}\SO^{r_2}_{S_0}@>>> 0
\endCD$$
and $\xi_1|_{S_0}=[r-r_1]\cdot
h_1^{S_0}$, $\xi_2|_{S_0}=[r-r_2]\cdot
h^{S_0}_2,$ we need to show that
there exists a point
$$(\SE',h')=(\SO^{\tilde n}@>e'>>\SE'\to
0,\SE'_{x_1}\oplus\SE'_{x_2}
@>h'>>\SO^r_S\to 0)\in \Phi_1(A)$$
such that
$\varphi_A(\SE',h')=(\wt\SE^S,q,\xi_,\xi_2)$ and
$(\SE',h')|
_{S_0}=(\SE^{S_0},h^{S_0}).$

We take $\SE'=\SE^S$ and
$e'=e^S$ to be a lifting of $e^{S_0}$, which 
always exists and is
surjective since $I^2=0$. To define $h'$,
let
$q_i:\SE^S_{x_i}\to\SO^{r_i}_S\to 0$ be the two induced surjections
by
$q$, and define $h'_i:\SE^S_{x_i}@>(q_i,\xi_i)>>\SO^r_S$, which gives
a
surjective morphism $h':\SE^S_{x_1}\oplus\SE^S_{x_2}\to\SO^r_S\to 0$
since
its restriction on $S_0$ is $h^{S_0}$. One can check that it is what
we
want.\enddemo

By the above lemma, we are reduced to consider the
singularities of $Y$.
To 
analyse the singularities of $Y$, we can fix a
$\wt E\in \bwq^1_F$ since
$$Y=W_1\times_VW_2\to \bwq^1_F$$
is locally
trivial, namely the singularities of $Y$ are the same with that
of any
fibre (Note that $\bwq^1_F$ is smooth).

\proclaim{Proposition 3.1} Let
$\wt E$ be a vector bundle of rank $r$ on
$\wt X$, $x_1,\,x_2\in\wt X$ and
$V=V(Ext^1(\,_{x_1}\Bbb C^{r_1}\oplus\,_{x_2}\Bbb C^{r_2},\,\wt
E)^{\vee})$, $p:V\times\wt X\to\wt X.$
Consider the universal
extension
$$0\to p^*\wt E\to
E\to\,_{x_1}\SO_V^{r_1}\oplus\,_{x_2}\SO_V^{r_2}\to 0$$
on $\wt X\times V$,
then the space
$\bold
E=V(E_{x_1}^{\oplus(r-r_1)})\times_V
V(E_{x_2}^{\oplus(r-r_2)})$ is
reduced, irreducible and normal with
rational
singularities.
\endproclaim

\demo{Proof} Replace $\wt X$ by an affine
neighbourhood of
$\{x_1,x_2\}$ where $\wt E$ is trivial. Furthermore, we
can assume that
$\wt X=\Bbb A^1$ such that $x_1=\{t=0\}$ and $x_2=\{t=1\}$,
namely,
$$\wt X=Spec\,\Bbb C[t]\supset\{(t),\,(1-t)\}\qquad\wt E=\Bbb
C[t]^r.$$
Let $F=Ext^1(\,_{x_1}\Bbb C^{r_1}\oplus\,_{x_2}\Bbb C^{r_2},\,\wt
E)$, and
$\{e_{ij}\}$ is a $\Bbb C$-basis of $F$, $\{x_{ij}\}=\{e_{ij}^*\}$
is the
dual basis of $\{e_{ij}\}$, then
$V=Spec(S(F^{\vee}))=Spec(\Bbb
C[\{x_{ij}\}])$
and
the
element
$$\gamma=\sum_{i,j} x_{ij}e_{ij}\in F^{\vee}\otimes
F\hookrightarrow
S(F^{\vee})
\otimes F=Ext^1(\,_{x_1}\Bbb
C[\{x_{ij}\}]^{r_1}\oplus\,_{x_2}\Bbb
C[\{x_{ij}\}]^{r_2}),p^*\wt E)$$
determines the universal extension
$$0\to p^*\wt E\to E\to\,
_{x_1}\Bbb
C[\{x_{ij}\}]^{r_1}\oplus\,_{x_2}\Bbb C[\{x_{ij}\}]^{r_2}\to
0.$$
To
construct the universal extension $E$, we need a resolution
of $_{x_1}\Bbb
C^{r_1}\oplus\,_{x_2}\Bbb C^{r_2}$, 
$$0\to\Bbb C[t]^{r_1}\oplus\Bbb
C[t]^{r_2}@>\alpha>>
\Bbb C[t]^{r_1}\oplus\Bbb C[t]^{r_2}
\to\,_{x_1}\Bbb
C^{r_1}\oplus\,_{x_2}\Bbb C^{r_2}\to 0\tag3.1$$
where $\alpha$ is defined
by 
$$\Bbb C[t]^{r_1}\oplus\Bbb C[t]^{r_2}
=\bigoplus^{r_1+r_2}_{i=1}
\Bbb
C[t]e_i,\qquad\alpha(e_i)=\cases te_i,&\text{for $i\leq
r_1$}\\
(1-t)e_i,&\text{for $i>r_1$}.\endcases$$
We have
$$F=\frac{Hom_{\Bbb C[t]}(\Bbb C[t]^{r_1}\oplus\Bbb C[t]^{r_2},\,
\Bbb
C[t]^r)}{\alpha^*Hom_{\Bbb C[t]}(\Bbb C[t]^{r_1}\oplus\Bbb
C[t]
^{r_2},\,\Bbb C[t]^r)}.$$
Define $e_{ij}:\Bbb C[t]^{r_1}\oplus\Bbb
C[t]^{r_2}\to\Bbb
C[t]^r\,1(\leq i\leq r_1+r_2,1\leq j\leq r)$ to
be
$$e_{ij}(e_i)= (\oversetbrace j\to{0,\cdots,0,1},0,\cdots,0),
\qquad
e_{ij}(e_k)=(0,\cdots,0)\quad\text{if $k\ne i$}.$$
Then $\{e_{ij}\}$ is a
basis of $F$.
Let $p^*\wt E=\Bbb C[t]^r\otimes_{\Bbb C}\Bbb
C[x_{ij}]=\Bbb
C[t,x_{ij}]^r=\SO^r_{V\times\wt X},$ and
$$0\to\Bbb
C[t,x_{ij}]^{r_1}\oplus\Bbb C[t,x_{ij}]^{r_2}@>\alpha>>
\Bbb
C[t,x_{ij}]^{r_1}\oplus\Bbb C[t,x_{ij}]^{r_2}
@>>>\,_{x_1}\Bbb
C[x_{ij}]^{r_1}\oplus\,_{x_2}\Bbb C[x_{ij}]^{r_2}\to 0$$
the pullback of
(3.1). The extension $E$ determined by 
$$\gamma:\Bbb
C[t,x_{ij}]^{r_1}\oplus\Bbb C[t,x_{ij}]^{r_2}\to 
p^*\wt E=\Bbb
C[x_{ij},t]^r$$ 
is
$$E=\frac{\Bbb C[x_{ij},t]^r\oplus\Bbb C[x_{ij},t]^{r_1+r_2}}
{W}=\frac{\sum^{r+r_1+r_2}_{k=1}\Bbb C[x_{ij},t]y_k}{W},$$
where $W=\left\{(\gamma(a),-\alpha(a))\,|\,a\in\Bbb
C[x_{ij},t]^{r_1+r_2}\right\}$. We can
describe $E$ by the following
exact
sequence
$$0\to\bigoplus^{r_1+r_2}_{k=1}\Bbb
C[x_{ij},t]e_k@>\beta>>
\bigoplus_{k=1}^{r+r_1+r_2}\Bbb C[x_{ij},t]y_k\to
E\to 0$$
where $\beta(e_k)=\gamma(e_k)\oplus(-\alpha(e_k))$ equals
to
$$\sum^r_{j=1}x_{kj}y_j-\cases ty_{r+k},&\qquad k\leq
r_1\\(1-t)y_{r+k}&\qquad k>r_1.\endcases$$
Thus we
get
$$0\to\bigoplus^{r_1+r_2}_{k=1}\Bbb
C[x_{ij}]e_k@>\beta_{x_i}>>
\bigoplus_{k=1}^{r+r_1+r_2}\Bbb C[x_{ij}]y_k\to
E_{x_i}\to 0\tag3.2$$
where 
$$\beta_{x_1}(e_k)=\cases
\sum^r_{j=1}x_{kj}y_j,&\text{if $k\leq
r_1$}\\
\sum^r_{j=1}x_{kj}y_j-y_{r+k}&\text{if
$k>r_1$}\endcases$$
$$\beta_{x_2}(e_k)=\cases\sum^r_{j=1}x_{kj}y_j-y_{r+k},&\text{if $k\leq
r_1$}\\ \sum^r_{j=1}x_{kj}y_j&\text{if
$k>r_1$}\endcases$$
Let
$$\bold
X_1=
\pmatrix
x_{11}&\hdots&x_{1r}\\
\vdots&\ddots&\vdots\\
x_{r_11}&\hdots&x_{r
_1r}
\endpmatrix\qquad\bold
X_2=
\pmatrix
x_{r_1+1,1}&\hdots&x_{r_1+1,r}\\
\vdots&\ddots&\vdots\\
x_{r_1+r_2
,1}&\hdots&x_{r_1+r_2,r}
\endpmatrix$$
We have 
$$S(E_{x_1})=\frac{\Bbb
C[\bold X_1,\bold X_2,\vec\bold y\,,y_{r+1},
\cdots,y_{r+r_1}]}{(\bold
X_1\cdot\vec\bold y)}\qquad\text{where
 $\vec\bold y=\pmatrix
y_1\\
\vdots\\
y_r\endpmatrix$}$$
and
$$S(E_{x_2})=\frac{\Bbb C[\bold
X_1,\bold X_2,\vec\bold y\,,y_{r+r_1+1},
\cdots,y_{r+r_1+r_2}]}{(\bold
X_2\cdot\vec\bold y)}.$$
Note that 
$$\SO_{\bold
E}=\overbrace{S(E_{x_1})\otimes_{\Bbb C[x_{ij}]}
\cdots\otimes_{\Bbb
C[x_{ij}]}S(E_{x_1})}^{r-r_1}\otimes_{\Bbb
C[x_{ij}]}
\overbrace{S(E_{x_2})\otimes_{\Bbb
C[x_{ij}]}\cdots\otimes_{\Bbb
C[x_{ij}]}S(E_{x_2})}^{r-r_2}$$
and set
$$y_{jl}=\overbrace{1\otimes\cdots\otimes
y_j}^l\otimes\cdots\otimes
1,\quad (1
\leq j\leq r+r_1+r_2,\,1\leq l\leq 2r-r_1-r_2)$$
we have
$$\SO_{\bold
E}=\frac{\Bbb C[\bold X_1,\bold X_2,\bold Y_1,\bold Y_2]}
{(\bold
X_1\cdot\bold Y_1,\,\bold X_2\cdot\bold Y_2)}\otimes_{\Bbb C}
\Bbb C[\bold
Y_3,\bold Y_4],$$
where $$\bold
Y_1=\pmatrix
y_{11}&\hdots&y_{1,r-r_1}\\
\vdots&\ddots&\vdots\\
y_{r1}&\hdots&y_
{r,r-r_1}
\endpmatrix\quad
\bold
Y_4=\pmatrix
y_{r+r_1+1,1}&\hdots&y_{r+r_1+1,2r-r_1-r_2}\\
\vdots&\ddots&\vdots\\
y_{r+r_1+r_2,1}&\hdots&y_{r+r_1+r_2,2r-r_1-r_2}
\endpmatrix$$
and $$\bold
Y_2=\pmatrix
y_{1,r-r_1+1}&\hdots&y_{1,2r-r_1-r_2}\\
\vdots&\ddots&\vdots\\
y_{r,r-r_1+1}&\hdots&y_{r,2r-r_1-r_2}
\endpmatrix\quad
\bold
Y_3=\pmatrix
y_{r+1,1}&\hdots&y_{r+1,2r-r_1-r_2}\\
\vdots&\ddots&\vdots\\
y_{r+r_1,1}&\hdots&y_{r+r_1,2r-r_1-r_2}
\endpmatrix$$
By taking $a=r_i$ and
$b=r-r_i$ in the following Lemma 3.3, note that
$I_{r_i}(\bold X_i)=\{0\}$
and $I_{r-r_i}(\bold Y_i)=\{0\}$, the
proposition is
proved.\enddemo

\proclaim{Lemma 3.3} Let $\bold X=(x_{ij})_{p\times r}$,
$\bold
Y=(y_{ij})_{r\times q}$ be two matrices and $I_a(\bold X)$
(resp.
$I_b(\bold Y)$) denote the set of rank $a+1$ (resp. $b+1$)
subdeterminants
of $\bold X$ (resp. $\bold Y$). Then
$$D_{a,b}=\text{Spec}\frac{\Bbb C[\bold X,\bold Y]}{(\bold X\cdot\bold
Y,I_a(\bold
X),I_b(\bold Y))}$$
is reduced, irreducible and normal with
rational singularities if $a+b\leq
r$.\endproclaim

\demo{Proof} The fact
that the variety is reduced, normal and Cohen-Macaulay
is a special case of
theorems in [CS]. Theorems in [He] imply that it has
only rational
singularities (See Example 6.5 of [He]).\enddemo

\remark{Remark 3.1} The above varieties $D_{a,b}$ were called double
determinantal varieties in [He], whose dimension formula is
$$dim(D_{a,b})=a(r+p)+b(r+q)-a^2-b^2-ab .$$
Take $a=p=r_i$, $b=q=r-r_i$ ($i=1,2$), we have $D_{r_i,r-r_i}=
\text{Spec}\,\Bbb C[\bold X_i,\bold Y_i]/(\bold X_i\cdot\bold Y_i)$ and
$dim(D_{r_i,r-r_i})=r^2+r_i^2-rr_i$. It is easy to see that the ideal 
$(\bold X_i\cdot\bold Y_i)$ has $r_i(r-r_i)$ generators at most and
$$dim(\text{Spec}\,\Bbb C[\bold X_i,\bold Y_i])-dim(D_{r_i,r-r_i})=r_i(r-r_i)
=height(I).$$
Thus $D_{r_i,r-r_i}$ are complete interesections. In particular, $\SH$ and $\SP$
are Gorenstein.\endremark 

\proclaim{Proposition
3.2} $\SH$, $\hat\SD_j(a)$ and $\hat\SD_1(a)\cap\hat\SD_2(b)$ are
reduced, normal with rational singularities. In particular,
$\SP$, $\SD_j(a)$ and $\SD_1(a)\cap\SD_2(b)$ are reduced, normal 
with rational singularities.\endproclaim

\demo{Proof}  By Lemma 3.2 and
Proposition 3.1, it is true for $\SH$.
We only need to show the proposition
for $\hat\SD_j(a)$ and $\hat\SD_1(a)\cap\hat\SD_2(b)$. Let us
rewrite
$(3.2)$ into
$$0\to\bigoplus^{r_1}_{k=1}\Bbb
C[x_{ij}]e_k@>\beta_{x_1}>>\bigoplus_{k=1}
^{r+r_1}\Bbb C[x_{ij}]y_k\to
E_{x_1}\to 0$$
$$0\to\bigoplus^{r_1+r_2}_{k=r_1+1}\Bbb
C[x_{ij}]e_k@>\beta_{x_2}>>
\bigoplus_{k=1}^{r+r_2}\Bbb C[x_{ij}]y_k\to
E_{x_2}\to 0$$
where $\beta_{x_i}(e_k)=\sum^r_{j=1}x_{kj}y_j.$ The
universal maps 
$E_{x_j}\to \SO^{r-r_j}$
are induced
by
$$(f_1,\cdots,f_{r+r_1})\overset u_1\to\longmapsto
(\sum^{r+r_1}_{i=1}
f_iy_{i1},\cdots,\sum^{r+r_1}_{i=1}f_iy_{i,r-r_1})$$
$$(f_1,
\cdots,f_{r+r_2})\overset u_2\to\longmapsto (\sum^{r+r_2}_{i=1}
f_iy_{i,r-r_1+1},\cdots,\sum^{r+r_2}_{i=1}f_iy_{i,2r-r_1-r_2}).$$
Let $E_{x_j}\to \SO^{r_j}$ be the induced projections by the projection
in the universal extension ($(f_1,\cdots,f_{r+r_j})\overset p_j\to
\longmapsto (f_{r+1},\cdots,f_{r+r_j})$). Then the maps $E_{x_j}\to\SO^r$
are induced by $\SO^{r+r_j}@>(p_j,u_j)>>\SO^r,$ matrices of which are
$$\pmatrix
\bold 0&\bold Y_1\\
\bold I_{r_1}&\bold Y'_3
\endpmatrix,\quad\pmatrix
\bold 0&\bold Y_2\\
\bold I_{r_2}&\bold Y^{''}_3
\endpmatrix$$
where $\bold I_{r_i}$ denote $r_i\times r_i$ unit matrice and
$(\bold Y'_3,\bold Y^{''}_3)=\bold Y_3$ (We use the notions in
Proposition 3.1). It is not difficult to see that the local smooth
models for
$\hat\SD_j(a)$ and $\hat\SD_1(a)\cap\hat\SD_2(b)$ at $z=(E,h)$ are
$$\text{Spec}\frac{\Bbb C[\bold X_j,\bold Y_j]}{(\bold X_j\cdot\bold
Y_j, I_{a-r_j}(\bold Y_j))}\times\text{Spec}\frac{\Bbb C
[\bold X_i,\bold Y_i]}{(\bold X_i\cdot\bold Y_i)}\quad (j=1,2,\,i\neq j)$$
and $$\text{Spec}\frac{\Bbb C[\bold X_1,\bold Y_1]}{(\bold X_1\cdot\bold
Y_1, I_{a-r_1}(\bold Y_1))}\times\text{Spec}\frac{\Bbb C
[\bold X_2,\bold Y_2]}{(\bold X_2\cdot\bold Y_2,I_{b-r_2}(\bold Y_2))}.$$
The proposition follows Lemma 3.3 (Note that $a\leq r$).
\enddemo

\remark{Remark 3.2} (1) It is easy to see from the proof that a point
$(E,h)\in\SH$ is a smooth point in the following cases: (i) $E$ is torsion
free at $x_1$ and $h_2:E_{x_2}\to\Bbb C^r$ is surjective, i.e., $r_2=r$
(The roles of $x_1$ and $x_2$ are reversed). (ii) Both of 
$h_j:E_{x_j}\to\Bbb C^r$ are surjective (i.e., $r_1=r_2=r$). In particular,
one can see that $\hat\SD_j(0)$ are smooth.

(2) The locus of non-locally-free extensions is 
$$Spec\frac{\Bbb C[\bold X_1,\bold Y_1,\bold X_2,\bold Y_2,\bold Y_3,\bold Y_4]}
{(\bold X_1\cdot\bold Y_1,\bold X_2\cdot\bold Y_2,I_{r_1-1}(\bold X_1))}\cup
Spec\frac{\Bbb C[\bold X_1,\bold Y_1,\bold X_2,\bold Y_2,\bold Y_3,\bold Y_4]}
{(\bold X_1\cdot\bold Y_1,\bold X_2\cdot\bold Y_2,I_{r_2-1}(\bold X_2))}.$$
More precisely, the non-locally-free locus $\SH\ssm\SH_F$ of $\SH$ has two
components $\hat\SD^t_j$ ($j=1,2$): $\hat\SD_1^t$ is the component of 
$\SH\ssm\SH_F$ parametrising sheaves with non-zero torsion at $x_2$ (The
$\hat\SD^t_2$ is defined similarly), whose local smooth models are
$$Spec\frac{\Bbb C[\bold X_1,\bold Y_1,\bold X_2,\bold Y_2,\bold Y_3,\bold Y_4]}
{(\bold X_1\cdot\bold Y_1,\bold X_2\cdot\bold Y_2,I_{r_1-1}(\bold X_1))}\quad
\text{and}\quad Spec\frac{\Bbb C[\bold X_1,\bold Y_1,\bold X_2,\bold Y_2,\bold Y_3,\bold Y_4]}
{(\bold X_1\cdot\bold Y_1,\bold X_2\cdot\bold Y_2,I_{r_2-1}(\bold X_2))}.$$
We will give more information about the subschemes $\hat\SD^t_j$ of $\SH$ 
in the following Proposition 3.3.\endremark

\proclaim{Proposition 3.3} Let $\tilde{\psi}':\widetilde{\Cal R}^{\prime ss}\to \Cal P$
be the projection. Then we have\roster\item 
$\tilde{\psi}'(\hat\SD^t_1\cap\wt\SR^{\prime ss})=\SD_1$, 
$\tilde{\psi}'(\hat\SD^t_2\cap\wt\SR^{\prime ss})=\SD_2$ and
\item the codimension one subschemes
$\hat\SD^t_j$ in $\SH$ are reduced, irreducible and normal.\endroster\endproclaim

\demo{Proof} (1) We will prove that $\tilde{\psi}'(\hat\SD^t_2\cap\wt\SR^{\prime ss})=\SD_2$, 
the other one is similar.
For any $(E',Q')\in \hat\SD^t_2\cap\wt\SR^{\prime ss}$, $\text{Tor$(E')_{x_1}\neq 0$}$ by
the definition, thus it is $s$-equivalent to a semistable GPS $(E,Q)$ with $E$ 
locally free by Lemma 2.5. Moreover, by checking the proof of Lemma 2.5, one find
that $E_{x_2}\to Q$ has rank $r-1$, so $\tilde{\psi}'((E,Q))\in\SD_2$. Each point of $\SD_2$ 
is the image of a GPS $(E,Q)$ with $E$ locally free, 
and $E_{x_2}\to Q$ is not surjective. Thus $E_{x_1}\to Q$ is nonzero, 
and we can choose a projection $Q\to\Bbb C$ such that $E_{x_1}\to Q\to \Bbb C$
is nonzero. Take $\wt E$ to be the kernel of $E\to \,_{x_1}\Bbb C\to 0$ and $\wt Q$ 
to be the image of $\wt E_{x_1}\oplus\wt E_{x_2}$ under $E_{x_1}\oplus E_{x_2}\to Q$, 
we get an
extension $$0\to(\wt E,\wt Q)\to (E,Q)\to (\, _{x_1}\Bbb C,\,\Bbb C)\to 0,$$
and one checks that $(E,Q)$ is $s$-equivalent to 
$(\wt E\oplus\, _{x_1}\Bbb C,\wt Q\oplus \,\Bbb C).$ 
Hence we proved that $\tilde{\psi}'(\hat\SD^t_2\cap\wt\SR^{\prime ss})=\SD_2.$

To prove (2), we only need to check the irreducibility, and the other facts 
follow the Remark 3.2 (2) and Lemma 3.3. On $\wt X\times\SH$, there is an exact
sequence $$0\to \Cal K\to\SO^{\tilde n}\to \SE\to 0$$
where $\Cal K$ is a vector bundle. It is easy to see that $\hat\SD^t_1$ is the
subscheme of $\SH$ defined by
$$\hat\SD^t_1=\{h\in\SH\,|\,\text{rank$(\Cal K_{(x_2,h)}\to 
\SO^{\tilde n}_{(x_2,h)})\le \tilde n-r-1$}\}.$$
Thus we only need to prove the open subset 
$$\hat\SD^{t,0}_1=\{h\in\hat\SD^t_1\,|\,rank(\Cal K_{(x_2,h)}\to\SO^{\tilde n}
_{(x_2,h)})=\tilde n-r-1\}$$
of $\hat\SD^t_1$ is irreducible. $\hat\SD^{t,0}_1$ is the open subset
of sheaves of the form $\wt E\oplus \,_{x_2}\Bbb C$ with $\wt E$ 
generated by global sections and having vanishing $H^1(\wt E)$. 
It is now straightforward
to imitate the proof of Remark 5.5 of [Ne].
  
\enddemo

We have shown (See Remark 3.1) that $\SH$ is Gorenstein, so it
has a canonical sheaf. Before closing this section, we will 
prove a formula to express the canonical sheaf of $\SH$. Let $$\SO^{\tilde n}\to
\SE\to 0,\quad \SE_{x_1}\oplus\SE_{x_2}\to\SQ\to 0$$
be the universal quotients on $\wt X\times\SH$ and $\SH$, we
write down an obvious lemma at first

\proclaim{Lemma 3.4} Let $\omega_{\wt X}=\SO(\sum_qq)$ be the
canonical sheaf of $\wt X$ and $\omega_{\wt\SR'_F}$ denote the
canonical bundle of $\wt\SR'_F$. Then
 $$\aligned&\omega^{-1}_{\wt\SR'_F}=(det\,R\pi_{\wt\SR'_F}\SE)^{2r}\otimes
\bigotimes_{x\in I}\left\{(det\,\SE_x)^{n_{l_x+1}-r}
\otimes\bigotimes^{l_x}_{i=1}(det\,\SQ_{x,i})
^{n_i(x)+n_{i+1}(x)}\right\}\\
&\otimes\bigotimes_q(det\,\SE_q)^{1-r}\otimes
(det\,R\pi_{\wt\SR'_F}det\SE)^{-2}\otimes(det\,\SQ)^{2r}\otimes
(det\,\SE_{x_1})^{-r}\otimes(det\,\SE_{x_2})^{-r}.
\endaligned$$\endproclaim

\demo{Proof} $\wt\SR'_F\to\wt\SR_F$ is a grassmannian bundle
over $\wt\SR_F$, then use Proposition 2.2.\enddemo 

We will give an extension of the right-hand side of the above
formula to $\SH$ as a line $PGL(\tilde n)$-bundle, 
then to prove that the extension gives the canonical bundle of 
$\SH$. Note that we have an exact squence
$$0\to\SK\to\SO^{\tilde n}\to\SE\to 0$$
on $\wt X\times\SH$, and $\SK$ is flat over $\SH$ since $\SE$
is so. One proves that $\SK$ is locally free on $\wt X\times\SH$
(By using Lemma 5.4 of [Ne]). For $x\in\wt X\ssm\{x_1,x_2\}$,
we have the identity $det\,\SK_x)^{-1}=det\,\SE_x$ on $\SH$. 
It is clear that  
$$\aligned\Omega^{-1}:&=(det\,R\pi_{\SH}\SE)^{2r}\otimes
\bigotimes_{x\in I}\left\{(det\,\SE_x)^{n_{l_x+1}-r}
\otimes\bigotimes^{l_x}_{i=1}(det\,\SQ_{x,i})
^{n_i(x)+n_{i+1}(x)}\right\}\otimes\\
&\bigotimes_q(det\,\SE_q)^{1-r}\otimes
(det\,R\pi_{\SH}det\,\SK^{-1})^{-2}\otimes(det\,\SQ)^{2r}\otimes
(det\,\SK_{x_1})^r\otimes(det\,\SK_{x_2})^r
\endaligned$$
is an extension of the line bundle in Lemma 3.4. We now prove 
that it is the dual of the canonical sheaf of $\SH$.

\proclaim{Proposition 3.4} Let $\SK$ be the kernel of the
universal surjection $\SO^{\tilde n}\to\SE$ on $\wt X\times\SH$,
and $\omega_{\SH}$ denote the canonical bundle of $\SH$. Then
$$\aligned&\omega_{\SH}^{-1}=\Omega^{-1}=(det\,R\pi_{\SH}\SE)^{2r}\otimes
\bigotimes_{x\in I}\left\{(det\,\SE_x)^{n_{l_x+1}-r}
\otimes\bigotimes^{l_x}_{i=1}(det\,\SQ_{x,i})
^{n_i(x)+n_{i+1}(x)}\right\}\\
&\otimes\bigotimes_q(det\,\SE_q)^{1-r}\otimes
(det\,R\pi_{\SH}det\,\SK^{-1})^{-2}\otimes(det\,\SQ)^{2r}\otimes
(det\,\SK_{x_1})^r\otimes(det\,\SK_{x_2})^r.
\endaligned$$\endproclaim

\demo{Proof} By Lemma 3.4, $\omega_{\SH}^{-1}=\Omega^{-1}$ holds outside
the $\hat\SD^t_j$, we will check that it extends to each $\hat\SD^t_j$.
For definiteness take $j=1$ and for simplicity of notation suppose there 
is no ordinary parabolic point. We will determine $\omega_{\SH}$ in a neighbourhood
of a suitable point of $\hat\SD^t_1$. Since $\hat\SD^t_1$ is irreducible, it
will be enough to show that $\omega_{\SH}^{-1}=\Omega^{-1}$ holds in one such 
neighbourhood.

We consider a point $(\SO^{\tilde n}\to E\to 0,Q)$ in $\SH$ satisfying
\roster\item $E$ has torsion at $x_2$ (i.e., the point lies on $\hat\SD^t_1)$,
\item $E$ is locally free at $x_1$, and
\item the maps $E_{x_j}\to Q$ are surjective for both $j=1,2$.\endroster
The conditions (2) and (3) will hold in a neighourhood $U$ of the point. On
$\wt X\times U$, we define a locally free sheaf $\wt\SE$ by the exact sequence
$$0\to\wt\SE\to\SE\to \,_{x_2}\SQ\to 0,\tag3.3$$
where $_{x_2}\SQ$ is the sheaf on $\wt X\times\SH$ got by pulling back $\SQ$ 
from $\SH$ and then restricting to $\{x_2\}\times\SH$. We can assume that for
any $u\in U$ the vector bundle $\wt\SE_u$ is generated by global sections and
$H^1(\wt\SE_u)=0$. Thus $\pi_{U*}\wt\SE$ is a locally free sheaf of rank
$\tilde n-r$ and commutes with any base change ($\pi_U$ denotes the projection
$\wt X\times U\to U$). Let $p:F_U\to U$ denote the frame bundle of   $\pi_{U*}\wt\SE$. 
We will use the same notation $\wt\SE$, $\SE$ and $\SQ$ to denote 
their pulling back to $\wt X\times F_U$ and $F_U$.

Let $\bwq$ be the Quot scheme of rank $r$, degree $d-r$ quotients 
$$\SO^{\tilde n-r}\to\wt E'\to 0$$
and $\bwq_F$ the open subset of locally free quotients generated by global 
sections with $H^1(\wt E')=0$ and $\SO^{\tilde n-r}\to\wt E'\to 0$
induces isomorphism $\Bbb C^{\tilde n-r}\cong H^0(\wt E').$ Let 
$$\SO^{\tilde n-r}_{\wt X\times\bwq_F}\to\wt\SE'\to 0$$
be the universal quotient, then there is a morphism $f_1:F_U\to\bwq_F$ such
that $(1\times f_1)^*\wt\SE'=\wt\SE$. Let $\pi_{\bwq_F}:\wt X\times\bwq_F\to
\bwq_F$ be the projection and 
$\Bbb E=Ext^1_{\pi_{\bwq_F}}(\,_{x_2}\wt\SE'_{x_1},\wt\SE')$ 
(See [La] for the 
definition of the sheaf), where $_{x_2}\wt\SE'_{x_1}$ is the sheaf on 
$\wt X\times\bwq_F$ got by pulling back $\wt\SE'_{x_1}$ from $\bwq_F$ and 
then restricting to $\{x_2\}\times\bwq_F$. Then there exists a universal
extension
$$0\to(1\times q_1)^*\wt\SE'\to\SE'\to\,_{x_2}q_1^*\wt\SE'_{x_1}\to 0$$
on $\wt X\times V(\Bbb E^{\vee}),$ where $q_1:V(\Bbb E^{\vee})\to\bwq_F$
is the projection. Note that $\wt\SE_{x_1}=\SE_{x_1}\cong\SQ$ on $U$ and
(3.3), we have a morphism $f_2:F_U\to V(\Bbb E^{\vee})$ such that
$f_1=q_1\cdot f_2$ and $(1\times f_2)^*\SE'=\SE$. Let $q:F_V\to V(\Bbb E^{\vee})$ 
be the frame bundle of $\pi_{V*}\SE'$, where $\pi_V:\wt X\times
V(\Bbb E^{\vee})\to V(\Bbb E^{\vee})$ is the projection, and 
$\SO^{\tilde n}_{F_V}@>u>>q^*(\pi_{V*}\SE')$ the universal frame. 
Then, if we denotes the pulling back of $\wt\SE'$ and $\SE'$ still by 
$\wt\SE'$ and $\SE'$,  there is a morphism $f:F_U\to F_V$ such that
$(1\times f)^*\wt\SE'=\wt\SE$, $(1\times f)^*\SE'=\SE$ and
$$\CD
\SO^{\tilde n}_{F_U}  @>f^*(u)>>     f^*(\pi_{F_V*}\SE')  \\
@|                                      @|            \\
\SO^{\tilde n}_{F_U}      @>>>    \pi_{F_U*}\SE
\endCD$$
is commutative, where $\SO^{\tilde n}_{F_U}\to \pi_{F_U*}\SE$ is the induced
isomorphism got by taking direct image of 
$\SO^{\tilde n}_{\wt X\times F_U}\to\SE\to 0$.

It is not difficult to check that $f:F_U\to F_V$ is unramified. In fact, the
universal frame $\SO^{\tilde n}_{F_V}@>u>>q^*(\pi_{V*}\SE')$ induces a quotient
$$\SO^{\tilde n}_{\wt X\times F_V}@>1\otimes u>>(1\times q)^*(\pi_V^*\pi_{V*}\SE')
\to(1\times q)^*\SE':=\SE'\to 0$$
on $\wt X\times F_V$, which gives a morphism $g_1:F_V\to \wt\SR$ such that
$\SE'=(1\times g_1)^*\SE$, the universal extension gives a quotient
$$\SE'_{x_1}\oplus\SE'_{x_2}\to\wt\SE'_{x_1}\to 0$$
on $F_V$, and thus a morphism $g_2:F_V\to U$ such that $g_2^*\SQ=\wt\SE'_{x_1}$
and $(1\times g_2)^*\wt\SE=\wt\SE'$. Finally, the universal quotient
$\SO^{\tilde n-r}\to\wt\SE'\to 0$ induces an isomorphism 
$\SO_{F_V}^{\tilde n-r}\cong\pi_{F_V*}\wt\SE'$, thus we have a morphism
$g:F_V\to F_U,$ which can be checked to be a section of $f:F_U\to F_V.$ Hence
$f$ is actually an isomorphism if $dim(F_U)=dim(F_V)$.

Now we check that $dim(F_U)=dim(F_V)$, it is easily to check that
$$dim(F_U)-dim(F_V)=r^2-rank(\Bbb E).$$
Thus we need to determine the locally free sheaf 
$\Bbb E=Ext^1_{\pi_{\bwq_F}}(\,_{x_2}\wt\SE'_{x_1},\wt\SE').$
Use the exact sequence
$$0\to\SO_{\wt X\times\bwq_F}(-\{x_2\}\times\bwq_F)
\otimes\pi_{\bwq_F}^*\wt\SE'_{x_1}\to\pi_{\bwq_F}^*\wt\SE'_{x_1}\to
\,_{x_2}\wt\SE'_{x_1}\to 0$$
and that $Ext^1_{\pi_{\bwq_F}}(\pi_{\bwq_F}^*\wt\SE'_{x_1},\wt\SE')
=R^1\pi_{\bwq_F*}(\wt\SE'\otimes\pi_{\bwq_F}^*\wt\SE^{'\vee}_{x_1})=0,$ we have
$$0\to\pi_{\bwq_F*}(\wt\SE')\otimes\wt\SE^{'\vee}_{x_1}\to
\pi_{\bwq_F*}(\wt\SE'\otimes\SO_{\wt X\times\bwq_F}(\{x_2\}\times\bwq_F))
\otimes\wt\SE^{'\vee}_{x_1}\to\Bbb E\to 0.\tag3.4$$
One can see easily that $rank(\Bbb E)=r^2$, and thus 
$\omega_{F_U}=f^*\omega_{F_V}=f^*q^*\omega_V$.

Since $\omega_V^{-1}=q_1^*\omega^{-1}_{\bwq_F}\otimes det(q_1^*\Bbb E)$ and
$det(\Bbb E)=(det\,\wt\SE'_{x_2})^r\otimes(det\,\wt\SE'_{x_1})^{-r}$ (By using
(3.4) and Riemann-Roch theorem), we have, using Lemma 2.10 and the pulling back
of (3.3),
$$\aligned\omega_{F_U}^{-1}=&(det\,R\pi_{F_U}\SE)^{2r}\otimes\bigotimes_q
(det\,\SE_q)^{1-r}\otimes(det\,R\pi_{F_U}det\,\wt\SE)^{-2}\\
&\otimes (det\,\SQ)^{2r}\otimes
(det\,\wt\SE_{x_2})^r\otimes(det\,\wt\SE_{x_1})^{-r}.\endaligned\tag3.5$$
On $\wt X\times F_U$, let $\SK'$ be the kernel of $\SO^{\tilde n}\to\SE\to
\,_{x_2}\SQ$, then we have the commutative diagram
$$\CD
@.        @.           0   @.        0  @.           \\
@.
@.                     @VVV          @VVV         @. \\
0  @>>> \SK         @>>> \SK' @>>> \wt\SE     @>>> 0  \\
@.      @|             @VVV         @VVV        @.  \\
0 @>>>\SK  @>>>  \SO^{\tilde n}@>>>\SE     @>>>0   \\
@.    @.             @VVV          @VVV        @. \\
@.       @.       \,_{x_2}\SQ@=\,_{x_2}\SQ       @.   \\
@.      @.              @VVV         @VVV        @.\\
@.         @.           0   @.        0   @.
\endCD\tag3.6$$
One sees easily that $\SK'$ is a vector bundle of rank $\tilde n$, and (Note that $F_U$ is smooth) 
$$det\,\SK'=\SO_{\wt X\times F_U}(-r\cdot\{x_2\}\times F_U).$$ 
Thus we can compute easily that
$$\aligned det\,R\pi_{F_U}det\,\wt\SE=&det\,R\pi_{F_U}(det\,\SK'\,\otimes\,
det\,\SK^{-1})\\&= det\,R\pi_{F_U}(det\,\SK^{-1})\otimes (det\,\SK_{x_2})^{-r}.
\endaligned$$
Note that $det\,\wt\SE_{x_2}=(det\,\SK_{x_2})^{-1}\otimes(det\,\SK'_{x_2})$ and
that $det\,\SK'_{x_2}$ is trivial, we have
$$\aligned\omega_{F_U}^{-1}=&(det\,R\pi_{F_U}\SE)^{2r}\otimes\bigotimes_q
(det\,\SE_q)^{1-r}\otimes(det\,R\pi_{F_U}det\,\SK^{-1})^{-2}\\
&\otimes (det\,\SQ)^{2r}\otimes
(det\,\SK_{x_2})^r\otimes(det\,\SK_{x_1})^r.\endaligned\tag3.7$$
Thus, by (3.7), we have $p^*\omega_U^{-1}=p^*\Omega^{-1}$, which shows
clearly that $\omega_{\SH}^{-1}=\Omega^{-1}$ holds on $U$ since $p:F_U\to U$
is locally trivial for the Zariski's topology, and we are done.
\enddemo

\heading \S 4 seminormality and the decomposition theorem \endheading

Let $I_Z$  denote the ideal sheaf of closed subscheme $Z$ in a scheme $X$. 
When $Z$ is of codimension one (not necessarily a Cartier divisor), 
we set $\SO_X(-Z):=I_Z$. If $\SL$ is a line bundle on $X$ and $Y$ is a closed
subscheme of $X$, we denote $\SL\otimes I_Z$ and the restriction $I_Z\otimes
\SO_Y$ of $I_Z$ on $Y$ by $\SL(-Z)$ and $\SO_Y(-Z)$. Now we collect 
some general facts at frist.

\proclaim{Lemma 4.1} Let $V$ be a projective scheme on which a reductive
group $G$ acts, $\wt\SL$ an ample line bundle linearising the $G$-action,
and $V^{ss}$ the open subscheme of semistable points. Let $V'$ be a
$G$-invariant closed subscheme of $V^{ss}$, $\bar V'$ its schematic closure
in $V$. Then
\roster\item $\bar V^{'ss}=V'$, and $V'\diagup\diagup G$ is a closed subscheme
of $V^{ss}\diagup\diagup G.$
\item $H^0(V^{ss},\wt\SL)^{inv}=H^0(W,\wt\SL)^{inv}$, where $W$ is an open 
$G$-invariant (irreducible) normal subscheme of $V$ containing $V^{ss}$ and
$(\quad)^{inv}$ denotes the invariant subspace for an action of $G$.
\endroster\endproclaim

\demo{Proof}  See Lemma 4.14 and Lemma 4.15 of [NR].\enddemo

\proclaim{Lemma 4.2} Let $V$ be a normal variety with a $G$-action, where $G$
is a reductive algebraic group. Suppose a good quotient $\pi:V\to U$ exists.
Let $\wt\SL$ be a $G$-line bundle on $V$, and suppose it descends as a line
bundle $\SL$ on $U$. Let $V^{\prime\prime}\subset V'\subset V$ be open $G$-invariant 
subvarieties of $V$, such that $V'$ maps onto $U$ and 
$V^{\prime\prime}=\pi^{-1}(U^{\prime\prime})$ for
some nonempty open subset $U^{\prime\prime}$ of $U$. Then any invariant section of
$\wt\SL$ on $V'$ extends to $V$.\endproclaim

\demo{Proof}  See Lemma 4.16 of [NR].\enddemo

\proclaim{Lemma 4.3} Suppose given a seminormal variety $V$, with
normalization $\sigma:\wt V\to V$. Let the non-normal locus be $W$,
endowed with its reduced structure. Let $\wt W$ be set-theoretic inverse
image of $W$ in $\wt V$, endowed with its reduced structure. Let
$N$ be a line bundle on $V$, and let $\wt N$ be its pull-back to
$\wt V$ ($\wt N=\sigma^*N$). Suppose $H^0(\wt V,\wt N)\to H^0(\wt W,\wt N)$
is surjective. Then
\roster
\item 
There is an exact sequence
$$0\to H^0(\wt V,\wt N\otimes I_{\wt W})\to H^0(V,N)\to H^0(W,N)\to 0.$$
\item If $H^1(W,N)\to H^1(\wt W,\wt N)$ is injective, so is
$H^1(V,N)\to H^1(\wt V,\wt N).$
\endroster\endproclaim

\demo{Proof} See Proposition 5.8 of [NR].\enddemo

\proclaim{Lemma 4.4} Let $G$ be a reductive group and $P$ a parabolic subgroup.
Let $\SL$ be an ample line bundle on $G/P$, and $X$ a union of Schubert
varieties with the reduced structure. Then \roster
\item the restriction map $H^0(G/P,\SL)\to H^0(X,\SL)$ is surjective, and
\item $H^i(G/P,\SL)$ and $H^i(X,\SL)$ vanish for $i>0$.\endroster\endproclaim

\demo{Proof} See Theorem 3 of [MRa].\enddemo

\remark{Remark 4.1} Let $F_1$ and $F_2$ be two vector spaces
of dimension $r$, and $Gr$ denote the grassmannian $Grass_r(F_1
\oplus F_2)$ of $r$-dimensional quotients. Let $E_j$ 
be the vector bundle on $Gr$ generated by $F_j$, and let
$E_1\oplus E_2\to Q$ be the universal quotient. Write 
$l_j=(det\,E_j)^{-1}\otimes det\,Q$ and 
$$D_j(a)=\{q\in Gr|\text{$rank(E_{jq}\to Q_q)\leq a$}\}.$$
The action of $GL(F_1)\times GL(F_2)$ on $Gr$ lifts to the line bundles $l_j$, 
and the subvarieties $D_j(a)$ are invariant under
the action of $GL(F_1)\times GL(F_2)$. Thus $H^0(l_j^k)$ and
$H^0(l_1^k|_{D_1(a)\cap D_2\cup D_1(a-1)})$ are $GL(F_1)\times
GL(F_2)$ modules. By Lemma 4.4, 
$$0\to H^0(l_1^k\otimes I_{D_1(a)\cap D_2\cup D_1(a-1)})
\to H^0(l_1^k)\to H^0(l_1^k|_{D_1(a)\cap D_2\cup D_1(a-1)})\to 0\tag4.1$$ 
is an exact sequence of $GL(F_1)\times GL(F_2)$ modules.
Thus it is splitting.\endremark

\proclaim{Proposition 4.1} The following maps are surjective for any
$1\le a\le r$
\roster \item $H^0(\SD_1(a),\Theta_{\SP})\to H^0(\SD_1(a)\cap\SD_2\cup\SD_1(a-1),\Theta_{\SP})$.
\item $H^0(\SD_1(a),\Theta_{\SP})\to H^0(\SD_1(a)\cap\SD_2,\Theta_{\SP}).$\endroster\endproclaim

\demo{Proof} It is clear that the following Proposition 4.2 implies
Proposition 4.1. \enddemo

\proclaim{Proposition 4.2} The following maps are surjective for any 
$1\le a\le r$
\roster
\item $H^0(\SP,\Theta_{\SP})\to H^0(\SD_1(a),\Theta_{\SP}).$ 
\item $H^0(\SP,\Theta_{\SP})\to H^0(\SD_1(a)\cap\SD_2\cup\SD_1(a-1),\Theta_{\SP}).$
\item $H^0(\SP,\Theta_{\SP})\to H^0(\SD_1(a)\cap\SD_2,\Theta_{\SP}).$
\endroster\endproclaim

\demo{Proof} We will deal with (2) in detail, the other statements will follow
the same proof. Let us consider the diagram
$$\CD
H^0(\wt\SR^{'ss},\hat\Theta')^{inv} @>a>>H^0(\hat\SD_1(a)^{ss}\cap\hat\SD_2^{ss}
\cup\hat\SD_1(a-1)^{ss},\hat\Theta')^{inv}   \\
 @AeAA                      @AfAA            \\ 
H^0(\SH,\hat\Theta')^{inv} @>>>H^0(\hat\SD_1(a)\cap\hat\SD_2
\cup\hat\SD_1(a-1),\hat\Theta')^{inv}   \\
 @VbVV                       @VdVV        \\
H^0(\wt\SR'_F,\hat\Theta')^{inv} @>c>>H^0(\hat\SD_{F,1}(a)\cap\hat\SD_{F,2}
\cup\hat\SD_{F,1}(a-1),\hat\Theta')^{inv}   \\
\endCD$$  
We need to prove that $a$ is surjective. The map $e$ is an isomorphism by Lemma 4.1. 
To prove that $b$ is an isomorphism,
it is enough to check that
$$H^0(\wt\SR^{'ss},\hat\Theta')^{inv}\to
H^0(\wt\SR^{'ss}\cap\wt\SR'_F,\hat\Theta')^{inv}$$
is an isomorphism. We use Lemma 4.2 with the identification
$V=\wt\SR^{'ss}$, $U=\SP$, $\pi=\tilde{\psi}'$, $V'=\wt\SR^{'ss}
\cap\wt\SR'_F$ and $U^{\prime\prime}=\SP\ssm(\SD_1\cup\SD_2)$ 
(One can show that $U^{\prime\prime}=\SP\ssm(\SD_1\cup\SD_2)$
is nonempty, for example, by Corollary 2.1). Lemma 2.5 shows 
that $V'=\wt\SR^{'ss}\cap\wt\SR'_F$
maps onto $U=\SP$. Thus $b$ is also an isomorphism.

Given a section $s$ of $H^0(\hat\SD_1(a)^{ss}\cap\hat\SD_2^{ss}
\cup\hat\SD_1(a-1)^{ss},\hat\Theta')^{inv}$, 
it extends to sections $s_1$, $s_2$ on 
$\hat\SD_1(a)\cap\hat\SD_2$ and $\hat\SD_1(a-1)$
by Lemma 4.1 since  $\hat\SD_1(a)\cap\hat\SD_2$ and 
$\hat\SD_1(a-1)$ are normal, which are equal on
$\hat\SD_1(a)^{ss}\cap\hat\SD_2^{ss}
\cap\hat\SD_1(a-1)^{ss}$. For any point $x\in\hat\SD_1(a)
\cap\hat\SD_2\cap\hat\SD_1(a-1)=\hat\SD_1(a-1)\cap\hat\SD_2$, 
we have $s_1(x)=s_2(x)$ if $x$ is semistable, and $s_1(x)=s_2(x)
=0$ if $x$ is nonsemistable (by the definition of semistability).
Thus the sections $s_1$ and $s_2$ yield a section on
$\hat\SD_1(a)\cap\hat\SD_2\cup\hat\SD_1(a-1)$, which is an 
extension of $s$. This proves that $f$ is an isomorphism. 
Hence we only need to prove that
$$H^0(\wt\SR'_F,\hat\Theta')^{inv} @>c>>H^0(\hat\SD_{F,1}(a)\cap\hat\SD_{F,2}
\cup\hat\SD_{F,1}(a-1),\hat\Theta')^{inv}$$
is surjective. Recall that $\rho:\wt\SR'_F\to\wt\SR_F$
is a grassmannian bundle over $\wt\SR_F$. 
The Lemma 4.4 implies that
$$0\to\rho_*(\hat\Theta'\otimes I_{\hat\SD_{F,1}(a)\cap\hat
\SD_{F,2}\cup\hat\SD_{F,1}(a-1)})\to\rho_*\hat\Theta'\to
\rho_*(\hat\Theta'|_{\hat\SD_{F,1}(a)\cap\hat
\SD_{F,2}\cup\hat\SD_{F,1}(a-1)})\to 0$$
is exact. In fact, we claim that the above sequence is splitting.
Note that
$$\hat\Theta'=\rho^*\Theta_{\wt\SR}\otimes(det\,\SQ)^k\otimes
(det\,\SE_y)^{-k}$$
and $\SE$ is the pullback of $\SF$ by $\rho$, we can rewrite
$$\hat\Theta'=\rho^*(\Theta_{\wt\SR}\otimes det\,\SF_y^{-k}
\otimes det\,\SF_{x_1}^{k})\otimes(det\,\SQ\otimes
det\,\SE_{x_1}^{-1})^{k}.$$
Let $X=\hat\SD_{F,1}(a)\cap\hat\SD_{F,2}\cup\hat\SD_{F,1}(a-1)$
and $\eta_{x_1}=(det\,\SQ)(det\,\SE_{x_1})^{-1}$. Then it is
enough to show that
$$0\to\rho_*(\eta_{x_1}^k\otimes I_X)\to\rho_*\eta_{x_1}^k\to
\rho_*(\eta_{x_1}^k\otimes\SO_X)\to 0\tag{4.2}$$
is splitting. The above direct image sheaves can be thought
of as vector bundles associated to representations of $GL(r)\times
GL(r)$ in (4.1) of Remark 4.1 (See Remark 5.10 of [NR]). Thus, by
the Remark 4.1, (4.2) is splitting and we proved the propostiotion.

\enddemo

In order to prove the following proposition, we need to show that
$\SW_a$ is seminormal for any $0\le a\le r$. However, we will admit
this fact, and prove it later, so that we can prove the decomposition
theorem as soon as possible.

\proclaim{Proposition 4.3} We have a (noncanonical) isomorphism 
$$H^0(\SU_X,\Theta_{\SU_X})\cong H^0(\SP,\Theta_{\SP}(-\SD_2)).$$
\endproclaim

\demo{Proof}  If we take, in the Lemma 4.3, $V=\SW_a$, $\wt V=\SD_1(a),$
$\sigma=\phi|_{\SD_1(a)}$ and $N=\Theta_{\SU_X}|_{\SW_a}$, then
we have $W=\SW_{a-1}$, $\wt W=\SD_1(a)\cap\SD_2\cup\SD_1(a-1)$ 
and $\wt N=\Theta_{\SP}|_{\SD_1(a)}$ by Proposition 2.1. Use the 
Proposition 4.1(1) and Lemma 4.3, we have
$$0\to H^0(\SD_1(a),\Theta_{\SP}\otimes I_{\SD_1(a)\cap\SD_2\cup\SD_1(a-1)})
\to H^0(\SW_a,\Theta_{\SU_X})\to H^0(\SW_{a-1},\Theta_{\SU_X})\to 0.$$
Thus we have a noncanonical isomorphism 
$$H^0(\SU_X,\Theta_{\SU_X})\cong H^0(\SW_0,\Theta_{\SU_X})\oplus
\bigoplus^r_{a=1}H^0(\SD_1(a),\Theta_{\SP}
\otimes I_{\SD_1(a)\cap\SD_2\cup\SD_1(a-1)}).\tag4.3 $$
If we define $\SD_1(-1)=\varnothing$, note that $\SD_1(0)\cong\SW_0$ and
$\SD_1(0)\cap\SD_2=\varnothing$ (By Lemma 2.4 (1) and Lemma 2.5), we can rewrite
(4.3) into
$$H^0(\SU_X,\Theta_{\SU_X})\cong 
\bigoplus^r_{a=0}H^0(\SD_1(a),\Theta_{\SP}
\otimes I_{\SD_1(a)\cap\SD_2\cup\SD_1(a-1)}).\tag4.4 $$
By Proposition 4.1, we have 
$$\aligned
& H^0(\SD_1(a),\Theta_{\SP}
\otimes I_{\SD_1(a)\cap\SD_2\cup\SD_1(a-1)})\oplus
H^0(\SD_1(a)\cap\SD_2\cup\SD_1(a-1),\Theta_{\SP})\\
&\cong H^0(\SD_1(a),\Theta_{\SP}),\endaligned$$
which and (4.4) implies that
$$\bigoplus^r_{a=0}H^0(\SD_1(a),\Theta_{\SP})\cong H^0(\SU_X,\Theta_{\SU_X})
\oplus \bigoplus^r_{a=0}H^0(\SD_1(a)\cap\SD_2\cup\SD_1(a-1),
\Theta_{\SP}).\tag4.5$$
By using the following exact sequence 
$$0\to\SO_{\SD_1(a)\cap\SD_2\cup\SD_1(a-1)}\to\SO_{\SD_1(a)\cap\SD_2}\oplus
\SO_{\SD_1(a-1)}\to\SO_{\SD_1(a-1)\cap\SD_2}\to 0,$$
and Proposition 4.1, we get
$$\aligned
&H^0(\SD_1(a)\cap\SD_2,\Theta_{\SP})\oplus H^0(\SD_1(a-1),\Theta_{\SP})\\
&\cong H^0(\SD_1(a-1)\cap\SD_2,\Theta_{\SP})\oplus
H^0(\SD_1(a)\cap\SD_2\cup\SD_1(a-1),\Theta_{\SP}).\endaligned
\tag4.6$$
By (4.5) and (4.6), we have a noncanonical isomorphism
$$H^0(\SU_X,\Theta_{\SU_X})\oplus H^0(\SD_1(r)\cap\SD_2,\Theta_{\SP})
\cong H^0(\SD_1(r),\Theta_{\SP}).$$
On the other hand, by Proposition 4.2, we have the exact sequence
$$0\to H^0(\SP,\Theta_{\SP}(-\SD_2))\to H^0(\SP,\Theta_{\SP})\to
H^0(\SD_2,\Theta_{\SP})\to 0.$$
Thus we proved our proposition if one remarks that $\SD_1(r)=\SP$.
\enddemo

We recall some facts about the representation of $GL(n)$ (See [FH]).
For any partition $\lambda=(\lambda_1\ge\cdots\ge\lambda_n\ge 0)$,
we have the so called Schur functor $\Bbb S_{\lambda}$ and Schur
polynomial $S_{\lambda}$. One get all of irreducible representations
of $GL(n)$ by applying Schur functors $\Bbb S_{\lambda}$ to the standard
representation $V$ of $GL(n)$. We denote these representations 
$\Bbb S_{\lambda}(V)$ by $R_{\lambda}:=R_{\lambda_1,\cdots,\lambda_n},$
and $D_k=(\Lambda^nV)^{\otimes k},$ then
$$R_{\lambda_1+k,\cdots,\lambda_n+k}=R_{\lambda_1,\cdots,\lambda_n}\otimes
D_k,$$
and the dual of $R_{\lambda_1,\cdots,\lambda_n}$, which is isomorphic to
$\Bbb S_{\lambda}(V^*)$, is the representation
$R_{-\lambda_n,\cdots,-\lambda_1}.$ In a more fantastic language, $R_{\lambda}$ 
is the irreducible representation with highest weight
$$(\lambda_1-\lambda_2)\omega_1+(\lambda_2-\lambda_3)\omega_2+\cdots 
+(\lambda_{n-1}-\lambda_n)\omega_{n-1}+\lambda_n\omega_n,$$
where $\omega_1,\cdots,\omega_n$ are the fundamental weights defined by
$$\omega_i(diag(a_1,\cdots,a_n))=a_1+\cdots+a_i.$$
Let $N_{\mu v\lambda}$ denote the Littlewood-Richardson number. Then we
have a general decomposition over $GL(V)\times GL(W)$
$$\Bbb S_{\lambda}(V\oplus W)=\bigoplus N_{\mu v\lambda}(\Bbb S_{\mu}V
\otimes\Bbb S_vW)\tag4.7$$
the sum over all partitions $\mu$, $v$ such that the sum of the numbers 
partitioned by $\mu$ and $v$ is the number partitioned by $\lambda$.

For $j=1,2$, let $E_j$ be $r$-dimensional vector spaces, and $Gr$ denote
the grassmannian of $r$-dimensional quotients $E_1\oplus E_2\to Q$.
We still use $E_j$ to denote the vector bundle on $Gr$ generated by
$E_j$, and use $Q$ to denote the universal quotient $E_1\oplus E_2\to Q$
on $Gr$.

\proclaim{Lemma 4.5} Let $l_j$ denote the line bundle $(det\,E_j)^{-1}
\otimes det\,Q$ on $Gr$. Then we have a natural isomorphism of
$GL(E_1)\times GL(E_2)$ modules
$$H^0(Gr,l_2^m)=\bigoplus_{\mu}\Bbb S_{\mu}(E_1)\otimes\Bbb S_{\mu}(E_2^*),$$
where  $\mu=(\mu_1,\cdots,\mu_r)$ 
runs through the integers $0\le \mu_r\le\cdots\le\mu_1\le m$.\endproclaim

\demo{Proof} It is clear that $H^0(Gr,l_2^m)=(\Lambda^rE_2)^{-m}
\otimes H^0(Gr,(det\,Q)^m),$ the space $H^0(Gr,(det\,Q)^m)$ is an irreducible 
representation of $GL(2r)$ with highest weight
$m\omega_r$
(See \S 15.4 of [FH]). Thus $$H^0(Gr,l_2^m)=(\Lambda^rE_2)^{-m}\otimes
\Bbb S_{\lambda}(E_1\oplus E_2),$$ where $\lambda=(m,\cdots,m).$ Use
(4.7), we have
$$\Bbb S_{\lambda}(E_1\oplus E_2)=\bigoplus N_{\mu v\lambda}(\Bbb S_{\mu}
E_1\otimes\Bbb S_vE_2).$$
Clearly, if $N_{\mu v\lambda}\ne 0$, $\mu=(\mu_1,\cdots,\mu_r)$
must satisfy that $0\le\mu_r\le\cdots\le\mu_1\le m$. The skew schur
function $S_{\lambda/\mu}=|H_{\lambda_i-\mu_j-i+j}|$ ($i,j=1,\cdots,r$)
can be written into
$$S_{\lambda/\mu}=\sum N_{\mu v\lambda}S_v$$
in terms of ordinary schur polynomials $S_v$ (See \S 6 of [FH]), where
$S_v=|H_{v_i+j-i}|.$ On the other hand, for a given $\mu=(\mu_1,\cdots,\mu_r)$
and $\lambda=(m,\cdots,m)$,
$$S_{\lambda/\mu}=|H_{m-\mu_j-i+j}|=|H_{m-\mu_i+i-j}|=|H_{m-\mu_{r-i+1}+j-i}|=S_v$$
where $v=(m-\mu_r,\cdots,m-\mu_1)$. Thus $N_{\mu v\lambda}=0$ when
$v\ne (m-\mu_r,\cdots,m-\mu_1)$ and $N_{\mu v\lambda}=1$  when
$v=(m-\mu_r,\cdots,m-\mu_1).$ Note that $$\Bbb 
S_{(m-\mu_r,\cdots,m-\mu_1)}(E_2)=(\Lambda^rE_2)^m\otimes\Bbb S_{\mu}(E_2^*).\tag4.8$$
We have $$\Bbb S_{\lambda}(E_1\oplus E_2)=\bigoplus_{\mu}(\Lambda^rE_2)^m
\otimes\Bbb S_{\mu}(E_1)\otimes\Bbb S_{\mu}(E^*_2),$$
which proves that 
$$H^0(Gr,l_2^m)=\bigoplus_{\mu}\Bbb S_{\mu}(E_1)\otimes\Bbb S_{\mu}(E_2^*)
\tag4.9$$
where $\mu$ runs through the integers $0\le\mu_r\le\cdots\le\mu_1\le m$.
\enddemo

Given $\mu=(\mu_1,\cdots,\mu_r)$, $\Bbb S_{\mu}(E_1)$ is the irreducible
representation of $GL(r)$ with highest weight 
$$(\mu_1-\mu_2)\omega_1+\cdots+(\mu_{r-1}-\mu_r)\omega_{r-1}+\mu_r\omega_r.$$
We can rewrite it into (Forgetting the zero terms)
$$(\mu_{r_1(x_1)}-\mu_{r_1(x_1)+1})\omega_{r_1(x_1)}+\cdots+
(\mu_{r_l(x_1)}-\mu_{r_l(x_1)+1})\omega_{r_l(x_1)}+\mu_r\omega_r.$$
Let $d_i(x_1)=\mu_{r_i(x_1)}-\mu_{r_i(x_1)+1}$ ($i=1,\cdots,l$) and
$Flag_{\vec n(x_1)}(E_1)$ the flag variety of type
$\vec n(x_1)=(n_1(x_1),\cdots,n_l(x_1))$,
where $n_i(x_1)=r_i(x_1)-r_{i-1}(x_1)$ (we set $n_1(x_1)=r_1(x_1)$).
If we denote the universal flag on $Flag_{\vec n(x_1)}(E_1)$ by 
$$E_1=F_0(E_1)\supset F_1(E_1)\supset\cdots\supset
F_{l}(E_1)\supset
F_{l+1}(E_1)=0$$
and the quotient $E_1/F_i(E_1)$ by $\SQ_{x_1,i}$, then we 
have
$$H^0(Flag_{\vec n(x_1)}(E_1),\bigotimes^l_{i=1}(det\SQ_{i,x_1})^{d_i(x_1)})\otimes
(\Lambda^rE_1)^{\mu_r}=\Bbb S_{\mu}(E_1).\tag4.10$$
Similarly, if we set $r_i(x_2)=r-r_{l-i+1}(x_1)$ and 
$d_i(x_2)=d_{l-i+1}(x_1)$, we have
$$H^0(Flag_{\vec n(x_2)}(E_2),\bigotimes^l_{i=1}(det\SQ_{i,x_2})^{d_i(x_2)})\otimes
(\Lambda^rE_2)^{-\mu_1}=\Bbb S_{\mu}(E_2^*).\tag4.11$$
We remarks that $l$ may be zero, namely, $\mu=(\mu_r,\cdots,\mu_r)$ and
$\Bbb S_{\mu}(E_1)$ is the one-dimensional irreducible representation 
$(\Lambda^rE_1)^{\mu_r}$ in this case.

Recall that $\SO^{\tilde n}\to \SF\to 0$ is the universal quotient on 
$\wt X\times \bold{\wt Q}_F$ and
$$Grass_r(\SF_{x_1}\oplus\SF_{x_2})@>f>>\bold{\wt Q}_F,\quad\wt\SR_F=
\underset{x\in I}\to{\times_{\bold{\widetilde Q}_F}}Flag_{\vec n(x)}(\Cal F_x)
@>g>>\bold{\wt Q}_F.$$
We will use $\SE$ to denote the various pullbacks of $\SF$.
Let $g_j:Flag_{\vec n(x_j)}(\SF_{x_j})\to\bold{\wt Q}_F$ ($j=1,2$)
be the relative flag scheme of type $\vec n(x_j)$ and
$$\SE_{x_j}=F_0(\SE_{x_j})\supset F_1(\SE_{x_j})\supset\cdots\supset
F_{l}(\SE_{x_j})\supset
F_{l+1}(\SE_{x_j})=0$$
the universal flag on $Flag_{\vec n(x_j)}(\SF_{x_j})$. If we set
$\SQ_{x_j,i}=\SE_{x_j}/F_i(\SE_{x_j})$ and
$$\SL_1^{\mu}=(det\,\SE_{x_1})^{\mu_r}\otimes
\bigotimes^l_{i=1}(det\SQ_{x_1,i})^{d_i(x_1)},$$ $$\SL_2^{\mu}=(det\,\SE_{x_2})^{-\mu_1}\otimes
\bigotimes^l_{i=1}(det\SQ_{x_2,i})^{d_i(x_2)}$$
where the integers $l$, $n_i(x_j)$ and $d_i(x_j)$ ($i=1,\cdots,l$) 
were defined in (4.10) and (4.11) (determined by $\mu$). 
Then we have

\proclaim{Lemma 4.6} Let $\SE_{x_1}\oplus\SE_{x_2}\to \SQ$ be
the universal $r$-quotient on $Grass(\SE_{x_1}\oplus\SE_{x_2})$, and 
$$h^{\mu}:Flag_{\vec n(x_1)}(\SF_{x_1})\times_{\bold{\wt Q}_F}
Flag_{\vec n(x_2)}(\SF_{x_2})\to\bold{\wt Q}_F.$$
Write $\eta_x:=(det\,\SQ)(det\,\SE_x)^{-1}$ for a point
$x\in\wt X$. Then we have
$$f_*(\eta_{x_2}^m)=\bigoplus_{\mu}h^{\mu}_*(\SL^{\mu}_1\otimes\SL^{\mu}_2)$$
where $\mu=(\mu_1,\cdots,\mu_r)$ runs through the integers 
$0\le\mu_r\cdots\le\mu_1\le m$.
\endproclaim 

\demo{Proof} This is the immediate corollary of Lemma 4.5 and (4.10)-(4.11).
\enddemo

For $\mu=(\mu_1,\cdots,\mu_r)$, let $\SU^{\mu}_{\wt X}$ be the moduli space of
semi-stable parabolic bundles on $\wt X$ with parabolic structures at points
$I\cup\{x_1,x_2\}$ and weights $\vec a(x)$ for $x\in I$ (See definition 1.1 in
\S 1) and for $j=1,2$
$$\vec a(x_j)=(\mu_r,\mu_r+d_1(x_j),\cdots,\mu_r+\sum^{l-1}_{i=1}d_i(x_j),
\mu_r+\sum^l_{i=1}d_i(x_j)).$$
Let $$\Theta_{\SU^{\mu}_{\wt X}}=\Theta(k,\ell,\vec a,\vec\alpha,I\cup\{x_1,x_2\})$$
be the line bundle defined in Theorem 1.2 with $\alpha_{x_1}=\mu_r$ 
and $\alpha_{x_2}=k-\mu_1$. 
Then we have the decomposition theorem

\proclaim{Theorem 4.1} There exists a (noncanonical) isomorphism
$$H^0(\SU_X,\Theta_{\SU_X})\cong\bigoplus_{\mu}H^0(\SU^{\mu}_{\wt X},
\Theta_{\SU^{\mu}_{\wt X}})$$
where $\mu=(\mu_1,\cdots,\mu_r)$ runs through the integers 
$0\le\mu_r\le\cdots\le\mu_1\le k-1.$\endproclaim

\demo{Proof} We consider the commutative diagram
$$\CD
\wt\SR'_F=Grass_r(\SF_{x_1}\oplus\SF_{x_2})\times_{\bold{\wt Q}_F}\wt\SR_F
@>\rho>>     \wt\SR_F   \\
@Vp_1VV          @VgVV  \\
Grass_r(\SF_{x_1}\oplus\SF_{x_2})     @>f>>   \bold{\wt Q}_F       
\endCD$$
and note that $\hat\Theta'=\rho^*\Theta_{\wt\SR_F}\otimes(det\SQ)^k\otimes(det\SE_y)^{-k}$.
Then 
$$\hat\Theta'\otimes\SO(-\hat\SD_2)=\rho^*\Theta_{\wt\SR_F}\otimes(det\SE_y)
^{-k}\otimes(det\SE_{x_2})^k\otimes p_1^*\eta^{k-1}_{x_2}.$$
Thus 
$$\rho_*(\hat\Theta'\otimes\SO(-\hat\SD_2))=\Theta_{\wt\SR_F}\otimes(det\SE_y)
^{-k}\otimes(det\SE_{x_2})^k\otimes g^*(f_*\eta^{k-1}_{x_2}).$$
By Lemma 4.6, we have
$$\rho_*(\hat\Theta'\otimes\SO(-\hat\SD_2))=\bigoplus_{\mu}(\Theta_{\wt\SR_F}
\otimes(det\SE_y)^{-k}\otimes(det\SE_{x_2})^k)\otimes
g^*h^{\mu}_*(\SL^{\mu}_1\otimes\SL^{\mu}_2)\tag4.12$$
where $\mu=(\mu_1,\cdots,\mu_r)$ runs through the integers $0\le\mu_1\le\cdots
\le\mu_r\le k-1.$ Let
$$\wt\SR^{\mu}:=\wt\SR\times_{\bold{\wt Q}}Flag_{\vec n(x_1)}(\SF_{x_1})\times_{\bold{\wt Q}}
Flag_{\vec n(x_2)}(\SF_{x_2})=
\underset{x\in I\cup\{x_1,x_2\}}\to{\times_{\bold{\widetilde Q}}}Flag_{\vec n(x)}(\Cal F_x)$$
and
$$\hat\Theta_{\mu}=\Theta_{\wt\SR_F}\otimes(det\SE_y)^{-k}\otimes(det\SE_{x_2})^k)
\otimes\SL^{\mu}_1\otimes\SL^{\mu}_2.$$
Recall that (See \S 1)
$$\Theta_{\wt\SR_F}
=(det\,R\pi_{\wt\SR}\Cal E)^k\otimes\bigotimes_{x\in I}\lbrace
(det\,\Cal E_x)^{\alpha_x}\otimes\bigotimes^{l_x}_{i=1}
(det\,\SQ_{x,i})^{d_i(x)}\rbrace\otimes 
(det\,\Cal E_y)^{\tilde\ell}$$
and use the definition of $\SL^{\mu}_1$ and $\SL^{\mu}_2$, one has
$$\hat\Theta_{\mu}
=(det\,R\pi_{\wt\SR^{\mu}}\Cal E)^k\otimes\bigotimes_{x\in I\cup\{x_1,x_2\}}
\lbrace(det\,\Cal E_x)^{\alpha_x}\otimes\bigotimes^{l_x}_{i=1}
(det\,\SQ_{x,i})^{d_i(x)}\rbrace\otimes 
(det\,\Cal E_y)^{\ell}$$
with $\alpha_{x_1}=\mu_r$, $\alpha_{x_2}=k-\mu_1$ and $l_{x_1}=l_{x_2}=l$.
$\hat\Theta_{\mu}$ is the restrication to $\wt\SR_F^{\mu}$ of a line bundle
linearising the $SL(\tilde n)$-action on the projective variety $\wt\SR^{\mu}$
and $\SU^{\mu}_{\wt X}$ is the GIT quotient of the semistable points 
$(\wt\SR^{\mu})^{ss}\subset\wt\SR_F^{\mu}$. Note that 
$r_i(x_2)=r-r_{l-i+1}(x_1)$ and $d_i(x_2)=d_{l-i+1}(x_1)$, we can check that 
$$ \sum_{x\in
I\cup\{x_1,x_2\}}\sum^{l_x}_{i=1}d_i(x)r_i(x)+r\sum_{x\in
I\cup\{x_1,x_2\}}
\alpha_x+r\ell=k\tilde n.$$
Thus $\hat\Theta_{\mu}$ descends to the line bundle $\Theta_{\SU^{\mu}_{\wt X}}$
on $\SU^{\mu}_{\wt X}$, and
$$H^0(\SU^{\mu}_{\wt X},\Theta_{\SU^{\mu}_{\wt X}})=H^0(\wt\SR^{\mu\,ss},
\hat\Theta_{\mu})^{inv.}=H^0(\wt\SR^{\mu}_F,\hat\Theta_{\mu})^{inv.}.$$
Let $p^{\mu}:\wt\SR^{\mu}_F\to\wt\SR_F$ be the projection, then (4.12) can be
written into
$$\rho_*(\hat\Theta'\otimes\SO(-\hat\SD_2))=\bigoplus_{\mu}
p^{\mu}_*\hat\Theta_{\mu}.\tag4.13$$
Thus 
$$H^0(\wt\SR'_F,\hat\Theta'\otimes\SO(-\hat\SD_2))^{inv.}=\bigoplus_{\mu}
H^0(\wt\SR^{\mu}_F,\hat\Theta_{\mu})^{inv.}=\bigoplus_{\mu}
H^0(\SU^{\mu}_{\wt X},\Theta_{\SU^{\mu}_{\wt X}}).$$
On the other hand, since a section of $\hat\Theta'\otimes\SO(-\hat\SD_2)$
is also a section of $\hat\Theta'$, we have
$$H^0(\wt\SR^{'ss},\hat\Theta'\otimes\SO(-\hat\SD_2))^{inv.}=
H^0(\wt\SR^{'ss}\cap\wt\SR'_F,\hat\Theta'\otimes\SO(-\hat\SD_2))^{inv.}$$
by Lemma 4.2 (See the proof of Proposition 4.2 for details), and
$$H^0(\wt\SR^{'ss}\cap\wt\SR'_F,\hat\Theta'\otimes\SO(-\hat\SD_2))^{inv.}
=H^0(\wt\SR'_F,\hat\Theta'\otimes\SO(-\hat\SD_2))^{inv.}$$
by Lemma 4.1. Thus one get a canonical
decomposition
$$H^0(\SP,\Theta_{\SP}(-\SD_2))=\bigoplus_{\mu}H^0(\SU^{\mu}_{\wt X},
\Theta_{\SU^{\mu}_{\wt X}}).\tag4.14$$
The theorem follows Proposition 4.3 and the proof is completed.
\enddemo
\remark{Remark 4.2} The proof of the above theorem gives also a
decomposition of $\rho_*(\hat\Theta')$
$$\rho_*(\hat\Theta')=\bigoplus_{\mu}
p^{\mu}_*\hat\Theta_{\mu},$$
where $\mu=(\mu_1,\cdots,\mu_r)$ runs through the integers 
$0\le\mu_1\le\cdots\le\mu_r\le k.$\endremark

Now we are in the position to deal with the seminormality of subvarieties
$\SW_a$, which was actually hidden in some literatures ([Fa], [S2] and [Tr]).
Our task here is to reveal the fact in these literatures. In order to make 
our paper self-contained, we begin with the definition of seminormality 
(See [Sw] or [NR])
and we also assume that $|I|=0$ for simplicity.

\proclaim{Definition 4.1} An extension $A\subset B$ of reduced rings is 
subintegral if\roster
\item $B$ is integral over $A$
\item $Spec(B)\to Spec(A)$ is a bijection
\item $k_{A\cap\frak p}\to k_{\frak p}$ is an isomorphism for any
$\frak p\in Spec(B)$, where $k_{\frak p}=B_{\frak p}/\frak p B_{\frak p}.$
\endroster\endproclaim

\proclaim{Definition 4.2} Let $A\subset B$ be reduced rings, we say that
$A$ is seminormal in $B$ if there is no extension $A\subset C\subset B$
with $C\neq A$ and $A\subset C$ subintegral. We say that $A$ is seminormal
if it is seminormal in its total ring of quotients. A variety $V$ is seminormal
if its local ring at any point is seminormal.\endproclaim

\proclaim{Proposition 4.4} Let $V$ be an variety and $\hat\SO_p$ denote
the completion of $\SO_p$. Let $I_1$ and $I_2$ be two radical ideals in 
a ring $A$ such that $I_1+I_2$ is radical. Then we have\roster
\item $V$ is seminormal if, for any $p\in V$, $\hat\SO_p[[u_1,\cdots,u_n]]$
is seminormal for some $n$.
\item $A/I_1\cap I_2$ is seminormal if $A/I_1$ and $A/I_2$ are seminormal.
\item A GIT quotient of a seminormal variety is seminormal.\endroster
\endproclaim

\demo{Proof} See \S 3 of [NR].\enddemo

Let $\bold Q$ be the Quot scheme of semistable torsion free sheaves of
rank $r$ and degree $d$, and $\SF$ a universal sheaf on $\bold Q\times X$.
For any $q\in\hat\SW_a\subset\bold Q,$ we will prove that 
$\hat\SO_{\hat\SW_a,q}[[u_1,\cdots,u_n]]$ is seminormal for some $n$, 
which
will imply that $\hat\SW_a$, thus $\SW_a$, is seminormal by (1) and (3) of
Proposition 4.4. Without loss of generality, we can assume that $q$ is the
point of $\hat\SW_a$ such that $$\SF_{q,x}\cong m_x^{\oplus r}.$$

To work out the local model of $\hat\SW_a$ at $q$, we have to recall
the local model of $\bold Q$ at $q$ (See Huiti{\'e}me Partie III of [S2]).
It is known that there is a subspace $W\subset H^0(X,\SF^*_q(m))$ of 
dimension $r$ such that $\SF_q(-m)\to\SO\otimes W^*$ is injective and induces the 
canonical inclusion $m_x^{\oplus r}\subset\SO_x^{\oplus r}$ for some $m$ 
(Proposition 21 of [S2]). Let $\bl$ be the category of Artinian local 
$\Bbb C$-algebras, and $X_A=X\times Spec(A)$ for any object $A$ of $\bl$.
Let $\SO_{X_A}^p\to \SF_A\to 0$ be an exact sequence, which induces 
$\SO^p_X\to \SF_q\to 0$ on $X$, $W_A\subset H^0(X_A,\SF^*_A(m))$ a free
$A$-module of rank $r$ such that $W_A\otimes_AA/m_A=W$, then 
$$\SF_A(-m)\to W^*_A=Hom_{\SO_{X_A}}(\SO_{X_A}\otimes W_A,\SO_{X_A})$$
is an injective morphism (See Lemma 19 of [S2]). Write 
$T=\SO\otimes W^*/\SF_q(-m)$ and $T_x$ the restriction
of $T$ on $\{x\}$, one has the following functors
$$F,G,H:\bl\to\bset$$
$$F(A):=\{\text{isomorphic classes of $(\SO_{X_A}^p\to\SF_A\to 0,W_A)$}\}$$
$$G(A):=\{\text{isomorphic classes of $(\SO^p_{X_A}\to\SF_A\to 
0)$}\}$$
$$H(A):=\{\text{isomorphic classes of $(\SO^r_{X_A}\to T_A\to 0)$}\},$$
where $\SF_A$ and $T_A$ are $A$-flat, $T_A$ has support $\{x\}\times Spec(A)$,
the functors satisfy
$$F(A/m_A)=\{(\SO^p_X\to\SF_q\to 0,W)\}$$
$$G(A/m_A)=\{(\SO^p_X\to\SF_q\to 0)\}$$
$$H(A/m_A)=\{(\SO^r_X\to T_x\to 0)\}.$$
We have two morphisms $f_1:F\to G$ and $f_2:F\to H$ defined by
$$f_1((\SO_{X_A}^p\to\SF_A\to 0,W_A))=(\SO_{X_A}^p\to\SF_A\to 0)$$
$$f_2((\SO_{X_A}^p\to\SF_A\to 0,W_A))=(W_A^*\to 
W^*_A/\SF_A(-m)|_{\{x\}\times Spec(A)}\to 0).$$

\proclaim{Lemma 4.7} The morphisms $f_1:F\to G$ and $f_2:F\to H$
are formally smooth.\endproclaim

\demo{Proof} See Lemma 23 and Lemma 24 of [S2].\enddemo

Suppose that $R=\hat\SO_x\cong\Bbb C[T_1,T_2]/(T_1\cdot T_2)$, and
$u=\bar T_1$, $v=\bar T_2$ the elements of $R$. Then the matrices
$$\alpha=\pmatrix
u\cdot\bold I_r&\bold 0\\
\bold 0&v\cdot\bold I_r
\endpmatrix,\quad\beta=\pmatrix
v\cdot\bold I_r&u\cdot\bold I_r
\endpmatrix$$
determine an exact sequence
$$R^{2r}@>\alpha>>R^{2r}@>\beta>>R^r@>>>\Bbb C^r@>>>0.\tag{$*$}$$
We define the functor $\Phi:\bl\to\bset$ by associating a object 
$A\in\bl$ the set of isomorphic classes
$$(R\otimes_{\Bbb C}A)^{2r}@>\alpha_A>>(R\otimes_{\Bbb C}A)^{2r}
@>\beta_A>>(R\otimes_{\Bbb C}A)^{r}@>>>T_A@>>>0$$
of deformations of $(*)$, with $T_A=Coker(\beta_A)$ $A$-flat.
One proved that $\Phi$ is isomorphic to $H$ (See Proposition 29 of [S2]).
On the other hand, we consider the variety
$$Z=\{(X,Y)\in M(r)\times M(r)|\,X\cdot Y=Y\cdot X=0\}.$$
For any $(X,Y)\in Z(A)$, where $X=(x_{ij})_{r\times r}$, 
$Y=(y_{ij})_{r\times r}$, the matrices
$$\alpha_{A(X,Y)}=\pmatrix
u\cdot\bold I_r& X\\
Y&v\cdot\bold I_r
\endpmatrix,\quad\beta_{A(X,Y)}=\pmatrix
v\cdot\bold I_r-Y&u\cdot\bold I_r-X
\endpmatrix$$
determine, if $x_{ij}\in m_A$ and $y_{ij}\in m_A$ ($i,j=1,2,\cdots,r$), a
deformation
$$(R\otimes_{\Bbb C}A)^{2r}@>\alpha_{A(X,Y)}>>
(R\otimes_{\Bbb C}A)^{2r}
@>\beta_{A(X,Y)}>>(R\otimes_{\Bbb C}A)^{r}@>>>T_A@>>>0\tag{$(*)_{(X,Y)}$}$$
of $(*)$, which gives an element of $\Phi(A)$. In fact,
$\hat\SO_{Z,(0,0)}$ represents the functor $\Phi$ (See Proposition 28 of
[S2]). Thus we get the local model $Z$ of $\bold Q$ at $q$. It is not
difficult to see that the local model of $\hat\SW_a$ at $q$ is
$$Z'=\{(X,Y)\in Z|\,rk(X)+rk(Y)\le a\}$$
if we remark that $Im(\beta_{A(X,Y)})\otimes R/m_R$ has rank
$2r-rank(X)-rank(Y)$. Namely,
$$\hat\SO_{\hat\SW_a,q}[[u_1,\cdots,u_n]]\cong \hat\SO_{Z',(0,0)}
[[v_1,\cdots,v_t]]$$
for some $n$ and $t$. To prove the seminormality of $\hat\SW_a$,
we will need

\proclaim{Lemma 4.8} Let $R$ be a ring, $X=(x_{ij})_{r\times r}$,
$Y=(y_{ij})_{r\times r}$ and 
$$W(k_1,k_2)=\{(X,Y)|XY=YX=0,rank(X)\le k_1,rank(Y)\le k_2\}$$
$$I(k_1,k_2)=(XY,YX,I_{k_1}(X),I_{k_2}(Y))R[X,Y].$$
If $B(k_1,k_2)$ is the reduced coordinate ring of $W(k_1,k_2)$, then
$$B(k_1,k_2)=R[X,Y]/I(k_1,k_2).$$
Moreover, if $R$ is Cohen-Macaulay and normal, the $W(k_1,k_2)$ are
Cohen-Macaulay and normal if $k_1+k_2\le r$.\endproclaim

\demo{Proof} See Theorem 2.9 and Theorem 2.14 of [St].\enddemo

\proclaim{Theorem 4.2} The varieties $\hat\SW_a$ ($0\le a\le r$)
are seminormal. In particular, the varieties $\SW_a$ ($0\le a\le r$)
are seminormal.\endproclaim

\demo{Proof} By Proposition 4.4, we only need to check that
$\hat\SO_{Z',(0,0)}[[v_1,\cdots,v_t]]$ is seminormal. It is clear that
$$Z'=\bigcup^a_{k=0}W(k,a-k),\quad \SO_{Z'}=\Bbb C[X,Y]/\bigcap^a_{k=0}
I(k,a-k).$$
It is easy to check that, for any $0<l\le a$,
$$\bigcap^{a-l}_{k=0}I(k,a-k)+I(a-l+1,l-1)=\bigcap^{a-l}_{k=0}I(k,l-1).
\tag4.15$$
Thus one can use (2) of Proposition 4.4 to prove that $\SO_{Z'}$ is 
seminormal since the normality of $W(k,a-k)$. But (4.15) and the normality 
are unchanged under completion, a classic fact (for example, see \S 13
 of [ZS]), we proved the seminormality of $\hat\SO_{Z',(0,0)}[[v_1,\cdots,v_t]]$
by the same reason.
\enddemo 

\heading \S 5 Codimension computations and the vanishing theorems \endheading

We are going to prove the vanishing theorems in this section. For this purpose,
we need some computations of codimensions, which may have some independent
interest. Let $V$ be a vector space of dimension $r$ and $V'\subset
V$ a $r_1$-dimensional subspace. We denote the flag variety of type 
$\vec n=(n_1,\cdots,n_{l+1})$ by $Flag_{\vec n}(V)$, and its closed point
$$(V=V_0\supset V_1\supset\cdots\supset V_l\supset V_{l+1}=0)$$
by $F(V)$. We begin the story by the following lemma.

\proclaim{Lemma 5.1} For any partition $r_1=m_1+\cdots+m_{l+1}$ with
$m_i\ge 0$, let 
$$\Omega(V')=\{F(V)\in Flag_{\vec n}(V) |dim(V'\cap V_i)\ge r_1-(m_1+\cdots+m_i)\}.$$
Then we have
$$codim(\Omega(V'))=\sum^{l+1}_{j=1}(n_j-m_j)(r_1-\sum^j_{i=1}m_i).$$
\endproclaim

\demo{Proof} The closed points of $Flag_{\vec n}(V)$ can be expressed
as the quotients
$$(V=V/V_{l+1}\to V/V_l\to\cdots\to V/V_1\to 0)$$
with $dim(V/V_i)=n_1+\cdots+n_i$ and the closed points of $\Omega(V')$
is the points of $Flag_{\vec n}(V)$ such that 
$$rank(V'\to V/V_i)\le m_1+\cdots+m_i.$$
By Proposition 9.6 and Remark 9.16 of [Fu], there exists, for any $n\ge r$,
a unique permutation $\omega\in \Bbb S_n$ such that
$$r_{\omega}(n_1+\cdots+n_i,r_1)=m_1+\cdots+m_i$$
(We take $a_i=r_1$, $1\le i\le l+1$ and $b_i=n_1+\cdots+n_{l+2-i}$ 
in Proposition 9.6 of [Fu]). 
Thus the codimension of $\Omega(V')$ is 
$\ell(\omega)$ (See Proposition 8.1 of [Fu]). By Proposition 9.6 (c) of [Fu],
we compute that 
$$\ell(\omega)=\sum^{l+1}_{j=1}(n_j-m_j)(r_1-\sum^j_{i=1}m_i),$$
which proves the lemma.\enddemo

\proclaim{Proposition 5.1} With the same notation as before, we have the
following estimations of codimensions.
\roster\item $codim(\wt\SR^{ss}\ssm\wt\SR^s)\ge (r-1)(\tilde g-1)+1$ if
$|I|>0$, and $codim(\wt\SR^{ss}\ssm\wt\SR^s)\ge (r-1)(\tilde g-1)$ if
$|I|=0$.
\item $codim(\wt\SR_F\ssm\wt\SR^{ss})\ge (r-1)(\tilde g-1)+1.$
\endroster\endproclaim

\demo{Proof}  Recall that 
$\wt\SR_F=\underset{x\in I}\to{\times_{\bold{\widetilde Q}_F}}
Flag_{\vec n(x)}(\Cal F_x)$, the tangent space of $\bold{\widetilde Q}_F$
at point $(0\to K\to\SO^{\tilde n}\to E\to 0)\in \bold{\widetilde Q}_F$ is
$H^0(\wt X,K^{\vee}\otimes E)$. Since $\Bbb C^{\tilde n}\cong H^0(E)$ and
$H^1(E)=0$ (By the definition of $ \bold{\widetilde Q}_F$), we have exact
sequence
$$0\to H^0(E^{\vee}\otimes E)\to\Bbb C^{\tilde n}\otimes\Bbb C^{\tilde n}\to
H^0(K^{\vee}\otimes E)\to H^1(E^{\vee}\otimes E)\to 0.\tag5.1$$
Riemann-Roch theorem implies that $dim\,H^0(K^{\vee}\otimes E)=r^2(\tilde g-1)+1
+dim\,PGL(\tilde n).$ Thus
$$dim\,\wt\SR_F=r^2(\tilde g-1)+1+\sum_{x\in I}dim\,Flag_{\vec n(x)}(\SF_x)+dim\,PGL(\tilde n).$$
We will deal with (1) in detail. Consider a point $E\in \wt\SR^{ss}\ssm\wt\SR^s$,
it is an extension 
$$0\to E_1\to E\to E_2\to 0$$
by two vector bundles $E_j$ of rank $r_j$ and degree $d_j$ such that
$$pardeg(E_1)=\frac{r_1}{r}pardeg(E),\tag5.2$$
where we take the induced parabolic structures on $E_j$. Let 
$$E_x=F_0(E)_x\supset F_1(E)_x\supset\cdots\supset F_{l_x}(E)_x\supset
F_{l_x+1}(E)_x=0$$
be the quasi-parabolic structure of $E$ at $x\in I$, with weight
$$0\leq a_1(x)<a_2(x)<\cdots<a_{l_x+1}(x)\le k,$$
and let $m_i(x)=dim(E_{1\,x}\cap F_{i-1}(E)_x/E_{1\,x}\cap F_i(E)_x)$, then we rewrite
(5.2) into
$$rd_1-r_1d=\frac{1}{k}\sum_{x\in I}\sum^{l_x+1}_{i=1}(r_1n_i(x)-rm_i(x))a_i(x).
\tag5.3$$
We will now describe a (countable) number of quasi-projective varieties
parametrising such parabolic bundles.

For $j=1,2$, let $d_j$, $r_j$ and $\tilde n_j$ be integers such that
$d_1+d_2=d$, $r_1+r_2=r$ and $\tilde n_1+\tilde n_2=\tilde n$. For each
$x\in I$, let $m_1(x),\cdots,m_{l_x+1}$ be non-negative integers such
that $r_1=m_1(x)+\cdots+m_{l_x+1}(x)$ and 
$$rd_1-r_1d=
\frac{1}{k}\sum_{x\in I}\sum^{l_x+1}_{i=1}(r_1n_i(x)-rm_i(x))a_i(x).$$
Let $\bold{\widetilde Q}^j$ ($j=1,2$) be the Quot scheme of rank $r_j$, degree
$d_j$ quotients $$\SO^{\tilde n_j}\to E_j\to 0$$
and $\bold{\widetilde Q}^j_F$ the open subset of locally free quotients with
vanishing $H^1(E_j)$ such that $\Bbb C^{\tilde n_j}\cong H^0(E_j)$. Let $\SF_j$
be the universal quotient on $\wt X\times\bold{\widetilde Q}^j_F$,
$V=\bold{\widetilde Q}^1_F\times\bold{\widetilde Q}^2_F$ and 
$\SF=\SF_2^{\vee}\otimes\SF_1$ on $\wt X\times V$. If we set $f:\wt X\times V\to V$
and
$$V_{h^1}=\{y\in V|h^1(f^{-1}(y),\SF|_{f^{-1}(y)})=h^1\}.$$
Then $V_{h^1}$ are locally closed subschemes (with the reduced structure) of $V$, 
and $$V=\bigcup_{h^1\ge 0}V_{h^1}$$
$R^1f_*(\SF)$ is locally free of rank $h^1$ on $V_{h^1}$. We define varieties $P_{h^1}$
as follows
\roster\item If $h^1=0$, we set $P_{h^1}=V$ and $\SF^{h^1}=\SF_1\oplus\SF_2$ on $\wt X
\times V$ \item If $h^1>0$, we define $P_{h^1}=\Bbb P((R^1f_*\SF)^{\vee})$ to be the
projective bundle on $V_{h^1}$, and $\SF^{h^1}$ to be the universal extension
$$0\to\SF_1\otimes\SO_{P_h}(1)\to\SF^{h^1}\to\SF_2\to 0$$
on $\wt X\times P_{h^1}$.\endroster 
For any $x\in I$ and $v(x)=(r_1, d_1, h^1, m_1(x),\cdots,m_{l_x+1}(x))$, we define a 
locally closed subscheme of $Flag_{\vec n(x)}(\SF^{h^1}_x)$ to be
$$ X^0_{v(x)}=\left\{\aligned &(E_x=F_0(E)_x\supset F_1(E)_x\supset\cdots
\supset F_{l_x}(E)_x\supset F_{l_x+1}(E)=0)\\
&\in Flag_{\vec n(x)}(\SF^{h^1}_x)\,|\quad dim(F_i(E)_x\cap E_1)=
r_1-\sum^i_{j=1}m_j(x)\endaligned\right\}$$ 
and let $$X_v=\underset{x\in I}\to{\times_{P_{h^1}}}X^0_{v(x)}.$$

Each $X_v$ parametrises a family of parabolic bundles $E$, which occur as 
extensions $0\to E_1\to E\to E_2\to 0$ (the extension being split if $h^1=0$),
with parabolic structures at $x\in I$ of type $\vec n(x)=(n_1(x),\cdots,n_{l_x+1}(x))$, 
whose induced parabolic structures on $E_1$
are of type $(m_1(x),\cdots,m_{l_x+1}(x))$ (we will forget $m_i(x)$ if it is zero). 
The dimension of $X_v$ are not bigger than 
$$\cases (\tilde g-1)\sum r_i^2+\sum dim\,PGL(\tilde n_i)+2+h^1-1+
\sum_{x\in I}dim\,X_{v(x)}&\text{if $h^1>0$}\\
(\tilde g-1)\sum r_i^2+\sum dim\,PGL(\tilde n_i)+2+\sum_{x\in I}
dim\,X_{v(x)}&\text{if $h^1=0$}\endcases$$
where $i=1,2$. Let $X_v^{ss}$ be the open set of semistable parabolic bundles,
and let $F(v)$ be the frame-bundle of the direct image of $\SF(v)$ (the pullback
of $\SF^{h^1}$) on $X_v^{ss}$. There is a map from each $F(v)$ to 
$\wt\SR^{ss}\ssm\wt\SR^s$, and the union of the images covers $\wt\SR^{ss}\ssm
\wt\SR^s$
$$dim(Im\,F(v))=dim(F(v))-e$$
where $e$ is the infimum of the dimensions of the irreducible components of the fibres. 
If a vector bundle $E$ is an extension
$$0\to E_1\to E\to E_2\to 0\tag5.4$$
then it is easy to see that
$$dim(Aut(E))\ge\cases 2+dim(H^0(E_2^{\vee}\otimes E_1))&\text{if (5.4) is splitting}\\
1+dim(H^0(E_2^{\vee}\otimes E_1))&\text{if (5.4) is not splitting}\endcases$$
Since the $E$ are generated by sections, any automorphism of $E$ acts nontrivially 
on the frames of $H^0(E)$, we have
\roster\item $e\ge dim(H^0(E_2^{\vee}\otimes E_1))+\tilde n_1^2+\tilde n_2^2$ if $h^1=0$, and
\item $e\ge dim(H^0(E_2^{\vee}\otimes E_1))+\tilde n_1^2+\tilde n_2^2-1$ if $h^1>0$.\endroster
Thus, by using Riemann-Roch theorem, the codimension of the images are bounded below by
$$r_1r_2(\tilde g-1)+\sum_{x\in I}codim(X^0_{v(x)})+rd_1-r_1d.$$
By Lemma 5.1 and (5.3), note that $r_1r_2=r_1(r-r_1)\ge r-1$, we have
$$codim(\wt\SR^{ss}\ssm\wt\SR^s)\ge(r-1)(\tilde g-1)+\sum_{x\in I}
\left\{\aligned &\sum^{l_x+1}_{j=1}(r_1-\sum^j_{i=1}m_i(x))(n_j(x)-m_j(x))\\
&+\sum^{l_x+1}_{j=1}(r_1n_j(x)-rm_j(x))\frac{a_j(x)}{k}\endaligned\right\}.$$
Since $r_1=\sum^{l_x+1}_{i=1}m_i(x)$ and $r=\sum^{l_x+1}_{i=1}n_i(x)$, the
first statement of the lemma follows the following Lemma 5.2.

Now we prove (2) of the lemma, the arguments is word by word as above, except
that we replace the equality (5.3) by an inequality
$$rd_1-r_1d>\frac{1}{k}\sum_{x\in I}
\sum^{l_x+1}_{i=1}(r_1n_i(x)-rm_i(x))a_i(x).$$
\enddemo

\proclaim{Lemma 5.2} For any integers $n_j>0$ and $m_j\ge 0$ ($j=1,\cdots,l+1$)
with $n_j\ge m_j$, let $0<a_1<\cdots<a_{l+1}\le 1$ be rational numbers, then
$$\sum^{l+1}_{j=1}m_j\sum^{l+1}_{j=1}(n_j-m_j)+\sum^{l+1}_{j=1}m_j
\sum^{l+1}_{j=1}n_ja_j\ge\sum^{l+1}_{j=1}(\sum^j_{i=1}m_i)(n_j-m_j)+
\sum^{l+1}_{j=1}n_j\sum^{l+1}_{j=1}m_ja_j.$$
Moreover, if $\sum_{j=1}^{l+1}n_j>\sum_{j=1}^{l+1}m_j>0$, we have the strict
inequality
$$\sum^{l+1}_{j=1}m_j\sum^{l+1}_{j=1}(n_j-m_j)+\sum^{l+1}_{j=1}m_j
\sum^{l+1}_{j=1}n_ja_j>\sum^{l+1}_{j=1}(\sum^j_{i=1}m_i)(n_j-m_j)+
\sum^{l+1}_{j=1}n_j\sum^{l+1}_{j=1}m_ja_j.$$\endproclaim

\demo{Proof}  We check it by induction for $l$, let  $(*)$ denote
the inequality, $LHS(*)$ and $RHS(*)$ denote the $\lq$left (right)-hand side of $(*)$'. 
When $l=1$, we have
$$LHS(*)-RHS(*)=m_2(n_1-m_1)+(m_1n_2-m_2n_1)(a_2-a_1),$$
which satisfies the lemma. Assume that $(*)$ is true for $l-1$, then
$$LHS(*)-RHS(*)\ge m_{l+1}\sum^l_{j=1}(n_j-m_j)-
\sum^l_{j=1}(m_{l+1}n_j-n_{l+1}m_j)(a_{l+1}-a_j),$$   
which is a strict inequality if $\sum^l_{j=1}n_j>\sum^l_{j=1}m_j>0$. When
$m_{l+1}=0$, $LHS(*)-RHS(*)\ge n_{l+1}\sum^l_{j=1}(a_{l+1}-a_j)m_j\ge 0,$
which is strict if $\sum_{j=1}^{l+1}n_j>\sum_{j=1}^{l+1}m_j>0$. 
When $m_{l+1}>0$, we have that
$$LHS(*)-RHS(*)\ge m_{l+1}\left\{\sum^l_{j=1}(n_j-m_j)-
\sum^l_{j=1}(n_j-\frac{n_{l+1}}{m_{l+1}}m_j)(a_{l+1}-a_j)\right\}\ge 0,$$
which is strict if $n_{l+1}>m_{l+1}$. The lemma is proved.
\enddemo

\proclaim{Proposition 5.2} Let $\SD_1^f=\hat\SD_1(r-1)\cup\hat\SD_1^t$
and $\SD_2^f=\hat\SD_2(r-1)\cup\hat\SD_2^t$. Then\roster
\item $Codim(\SH\ssm\wt\SR^{\prime ss})\ge (r-1)\tilde g+1.$ 
\item The complement in $\wt\SR^{\prime ss}\ssm\{\SD_1^f\cup\SD^f_2\}$ of the set
$\wt\SR^{\prime s}$ of stable points has codimension $\ge (r-1)\tilde g+1$ if $|I|>0$,
and codimension $\ge (r-1)\tilde g$ if $|I|=0.$
\endroster\endproclaim 

\demo{Proof} We will prove (1) in detail, and (2) will follow similarly.
For any $(E,Q)\in\SH\ssm\wt\SR^{\prime ss}$ with $E_{x_1}\oplus E_{x_2}@>q>>Q @>>>0$, 
there exists a nontrivial subsheaf $E_1\subset E$, of
$rank(E_1)=r_1>0$, such that $E/E_1$ is torsion free outside $\{x_1,x_2\}$ and
$$pardeg(E_1)-dim(Q^{E_1})> \frac{r_1}{r}(pardeg(E)-r).\tag5.5$$
In fact, we can choose $E_1$ such that $E/E_1=E_2$ is torsion free. If $E_2$
has torsion $_{x_1}\tau_1\,\oplus\,_{x_2}\tau_2$, and let $\wt E_1\supset E_1$
be the inverse image in $E$ of $_{x_1}\tau_1\,\oplus\,_{x_2}\tau_2$. Then
$pardeg(\wt E_1)=pardeg(E_1)+dim(\tau_1)+dim(\tau_2)$ and
$$dim(Q^{\wt E_1})-dim(Q^{E_1})\le dim(\tau_1)+dim(\tau_2)=pardeg(\wt E_1)-pardeg(E_1),$$ 
which shows that $\wt E_1$ satisfies (5.5), and we can choose $\wt E_1$ instead
of $E_1$. Thus $E$ is an extension
$$0\to E_1\to E\to E_2\to 0$$
with $E_2$ torsion free (Note that $r_2=rank(E_2) >0$) and $E_1$ satisfying (5.5).

We can write $E=E'\,\oplus\,_{x_1}\Bbb C^{s_1}\oplus\,_{x_2}\Bbb C^{s_2}$ and
$E_1=E'_1\,\oplus\,_{x_1}\Bbb C^{s_1}\oplus\,_{x_2}\Bbb C^{s_2}$ with $E'$ and
$E_1'$ torsion free. Thus $(E,Q)$ is a GPS such that $E'=E/\text{Tor($E$)}$ 
occurs as an extension $$0\to E_1'\to E'\to E'_2\to 0\quad (\text{where $E_2'=E_2$})$$
with $pardeg(E'_1)>dim(Q^{E_1})+\frac{r_1}{r}pardeg(E)-r_1-s_1-s_2.$ When
$d=deg(E)=deg(E')+s_1+s_2$ is large enough (so is $deg(E')$ since $s_1+s_2\le r$), 
we can assume that $E'_1$ and $E_2'$ are generated by global sections and $H^1(E'_1)=H^1(E'_2)=0$.

Let $d_j=deg(E_j')$, $r_j=rank(E'_j)$ ($j=1,2$) and, for any $x\in I$,
$$m_i(x)=dim(E'_{1\,x}\cap F_{i-1}(E)_x/E'_{1\,x}\cap F_i(E)_x)$$ 
where $E_x=F_0(E)_x\supset F_1(E)_x\supset\cdots\supset F_{l_x}(E)_x\supset
F_{l_x+1}(E)_x=0$
is the quasi-parabolic structure of $E$ at $x\in I$ of type $(n_1(x),\cdots,n_{l_x+1}(x))$, 
with weights
$$0\le a_1(x)<a_2(x)<\cdots<a_{l_x+1}(x)\le k.$$
Let $t=dim(Q^{E_1})$ and $s=s_1+s_2$, then $s\le t\le r$ and
$$rd_1-r_1d>r(t-s-r_1)+\frac{1}{k}\sum_{x\in I}\sum^{l_x+1}_{i=1}(r_1n_i(x)-rm_i(x))a_i(x).\tag5.6$$
Let $v=(d_1,r_1,s_1,s_2,t,\{m_1(x),\cdots,m_{l_x+1}(x)\}_{x\in I},h)$, where 
$h\ge 0$ is an integer, we will construct a variety $F(v)$ with a morphism
$F(v)\to\SH\ssm\wt\SR^{\prime ss}$
such that its image contains the point $(E,Q)$.

For $j=1,2$, let $\tilde n_j=dim\,H^0(E_j')$, $\wt{\bold Q}^j$ the Quot scheme of 
rank $r_j$, degree $d_j$ quotients $\SO^{\tilde n_j}\to E_j'\to 0$
and $\wt{\bold Q}^j_F$ the open subset of locally free quotients with vanishing 
$H^1(E'_j)$ and $E'_j$ generated by global sections. Let $\SE'_j$ be the universal
quotient on $\wt X\times\wt{\bold Q}^j_F$, $V=\wt{\bold Q}^1_F\times\wt{\bold Q}^2_F$
and $\SF=\SE_2^{\prime\vee}\otimes\SE'_1$ on $\wt X\times V$. We have 
$V=\bigcup_{h\ge 0}V_h$ and $R^1f_*(\SF)$ is locally free of rank $h$ on $V_h$
(See the proof of Proposition 5.1), where $f:\wt X\times V\to V$ is the projection. 
Let $P_h=\Bbb P((R^1f_*\SF)^{\vee})$ be the projective bundle on $V_h$ and 
$$0\to\SE'_1\otimes\SO_{P_h}(-1)\to\SE'(h)\to\SE_2'\to 0\tag5.7$$
the universal extension on $\wt X\times P_h$ (We set $P_h=V$ and $\SE'(h)=\SE'_1
\oplus\SE'_2$ if $h=0$). For $v'=(d_1,r_1,\{m_1(x),\cdots,m_{l_x+1}(x)\}_{x\in I},h)$, 
as in the proof of Proposition 5.1, we can define a 
variety $X(v')\to P_h$.
It parametrises a family of parabolic bundles $E'$, which occur as extensions
$0\to E'_1\to E'\to E'_2\to 0$ (the extension being split if $h=0$), with
parabolic structures at $x\in I$ of type 
$\vec n(x)=n_1(x),\cdots,n_{l_x+1}(x))$, whose induced parabolic structures 
on $E'_1$ are of type $(m_1(x),\cdots,m_{l_x+1}(x))$ (we will forget $m_i(x)$ if it
is zero). Let $0\to\SE'_1(-1)\to\SE'(v')\to\SE_2'\to 0$ be the pull back of
(5.7) on $\wt X\times X(v')$, and $\SE(v')=
\SE'(v')\oplus\,_{x_1}\SO^{s_1}\,\oplus\,_{x_2}\SO^{s_2}.$  We consider
$$G_{v'}=Grass_r(\SE(v')_{x_1}\oplus\SE(v')_{x_2})\to X(v')$$ and define a
subvariety of $G_{v'}$ 
$$X(v):=\left\{\aligned
&(E_{x_1}\oplus E_{x_2}@>q>>Q@>>>0)\in G_{v'})\,\text{with
$dim(ker(q)\cap(\Bbb C^{s_1}\oplus\Bbb C^{s_2}))=0$}\\ 
&\text{and $dim(ker(q)
\cap(E'_{1\,x_1}\oplus\Bbb C^{s_1}\oplus E'_{1\,x_2}\oplus\Bbb C^{s_2}))=2r_1
+s-t $}\endaligned\right\}.$$
Then $X(v)$ parametrises a family of GPS $(E=E'\oplus\,_{x_1}\Bbb C^{s_1}\,
\oplus\,_{x_2}\Bbb C^{s_2},\,Q)$, where $E'$ occurs as an extension 
$0\to E'_1\to E'\to E'_2\to 0$ (it is split if $h=0$) 
with parabolic structures at $x\in I$ of type $\vec n(x)$, 
whose induced parabolic structures on $E'_1$ are of type 
$(m_1(x),\cdots,m_{l_x+1}(x))$ (we will forget $m_i(x)$ if it is zero), 
such that $_{x_1}\Bbb C^{s_1}\,\oplus\,_{x_2}
\Bbb C^{s_2}\to Q$ is injective and
$$\text{rank$(E'_{1\,x_1}\oplus\Bbb C^{s_1}\oplus E'_{1\,x_2}\oplus\Bbb C^{s_2}\to Q)=t$}.$$
One computes $dim\,X(v)=dim\,X(v')+r(r+s)-(r-t)(2r_1+s-t)$. Let $\SE(v)$ be
the pull back of $\SE(v')$ on $\wt X\times X(v)\to\wt X\times X(v')$, and
let $F(v)$ be the frame bundle of the direct image of $\SE(v)$ on $X(v)$,
then there is a morphism $F(v)\to\SH\ssm\wt\SR^{\prime ss}$ whose image contains $(E,Q)$. 

Therefore we have a (countable) number of quasi-projective
varieties $F(v)$ and morphisms $F(v)\to\SH\ssm\wt\SR^{\prime ss}$ such that 
the union of the images covers $\SH\ssm\wt\SR^{\prime ss}$. Since the sheaf
$E'\,\oplus\,_{x_1}\Bbb C^{s_1}\oplus\,_{x_2}\Bbb C^{s_2}$ has an 
automorphism group of dimension at least $\text{dim$\,Aut(E')+rs+s^2$}$,
and the dimension of $\SH$ is 
$$r^2(\tilde g-1)+1+r^2+\sum_{x\in I}dim\,
Flag_{\vec n(x)}(\SF_x)+dim\,PGL(\tilde n),$$ we
find that the codimension of $\SH\ssm\wt\SR^{\prime ss}$ is bounded below by
$$\aligned & r_1r_2(\tilde g-1)+s^2+r_1s+(r-t)(2r_1+s-t)+\\
& rd_1-r_1d+\sum_{x\in I} \sum^{l_x+1}_{j=1}(r_1-\sum^j_{i=1}m_i(x))(n_j(x)-m_j(x)).\endaligned$$
By using (5.6), we get
$$\aligned codim(\SH\ssm\wt\SR^{\prime ss})>&r_1r_2\tilde g+(r_1-t)^2+(r_1-t+s)s\\
&+\sum_{x\in I}
\left\{\aligned &\sum^{l_x+1}_{j=1}(r_1-\sum^j_{i=1}m_i(x))(n_j(x)-m_j(x))\\
&+\sum^{l_x+1}_{j=1}(r_1n_j(x)-rm_j(x))\frac{a_j(x)}{k}\endaligned\right\}
\endaligned.\tag5.8$$
It is clear that $(r_1-t)^2+(r_1-t+s)s\ge 0$ when $t\le r_1+s$. Otherwise, if
$t>r_1+s$, we have $(r_1-t)^2+(r_1-t+s)s=s^2+(t-r_1)(t-r_1-s)>s^2$. Thus
$$codim(\SH\ssm\wt\SR^{\prime ss})\ge (r-1)\tilde g+1,$$
and we have proved (1) of the proposition.

For any $(E,Q)\in\wt\SR^{\prime ss}\ssm\{\SD_1^f\cup\SD^f_2\}\ssm\wt\SR^{\prime s}$ with $E_{x_1}\oplus
E_{x_2}@>q>>Q@>>>0$, there is a subsheaf $\wt E\subset E$ such that $E/\wt E$ is
torsion free outside $\{x_1,x_2\}$, which contracdicting the stability. One can
show that $\wt E$ has to be of rank $0<r_1<r$. Otherwise, $\wt E$ must satisfy
the exact sequence 
$$0\to\wt E\to E\to\,_{x_1}\tau_1\oplus\,_{x_2}\tau_2\to 0$$
with $dim(\tau_1\oplus\tau_2)=dim(Q/Q^{\wt E})$, and we have the diagram
$$\CD
@.\wt E_{x_1}\oplus\wt E_{x_2}@>>> E_{x_1}\oplus E_{x_2}@>>> \tau_1\oplus\tau_2@>>>  0\\
@.    @VVV
@VqVV        @|        @. \\
0@>>> Q^{\wt E}    @>>>   Q   @>>>    Q/Q^{\wt E}  @>>>0.
\endCD$$
Since $q_j:E_{x_j}\to Q$ ($j=1,2$) are isomorphisms, $\tau_j$ have to be
zero. Thus (2) now follows the same proof except that we replace the inequality (5.6) by an equality. 
\enddemo

\remark{Remark 5.1} It is not true that 
$\wt\SR^{\prime ss}\ssm\wt\SR^{\prime s}$ has codimension $>1$. Points on
$\SD_1^f=\hat\SD_1(r-1)\cup\hat\SD_1^t$
and $\SD_2^f=\hat\SD_2(r-1)\cup\hat\SD_2^t$ are never stable (See Remark 1.2).
The above codimension bound breaks down because, for 
$(E,Q)\in\SD_1^f\cup\SD^f_2$, we can not assume that the
subsheaf contradicting stability is of rank $0<r_1<r$.
\endremark

We denote the Jacobian of degree $d$ line bundles on $\wt X$ by $J^d_{\wt X}$ 
and the Poincar{\'e} line bundle on $\wt X\times J^d_{\wt X}$ by $\SL$. Let
 $$\Theta_y:=(det\,R\pi_J\SL)\otimes(det\,\SL_y)^{d+1-\tilde g}$$
and $Det: \wt\SR_F\to J^d_{\wt X}$ the morphism given by the determinant of 
the universal quotient bundle. This induces a morphism 
$\SU_{\wt X}\to J^d_{\wt X}$, which will also be denoted by $Det$. 
On $\wt\SR_F$, one sees easily that
$$(det\,R\pi_{\wt\SR_F}det\,\SE)^{-2}=
(det\,\SE_y)^{2\tilde n+2(r-1)(\tilde g-1)}\otimes Det^*\Theta_y^{-2}.\tag5.9$$ 

\proclaim{Lemma 5.3} Let $Det: \SU_{\wt X}\to J^d_{\wt X}$ be the 
induced morphism by $Det:\wt\SR_F\to J^d_{\wt X}$. Then
$$\Theta_{\SU_{\wt X}}\otimes(Det^*\Theta_y)^{-2}$$
is ample if $k>2r$.\endproclaim

\demo{Proof} Let $\SU^L_{\wt X}$ be the fibre of $Det:\SU_{\wt X}\to
J^d_{\wt X}$ at $L\in J^d_{\wt X}$. One has a $r^{2\tilde g}$-fold covering
$$f:\SU^L_{\wt X}\times J^0_{\wt X}\to \SU_{\wt X}$$
given by $f(E,L_0)=E\otimes L_0$. We will show that 
$\Theta_{\SU_{\wt X}}\otimes(Det^*\Theta_y)^{-2}$ is ample when pulled 
back to this finite cover.

One can show that $\SU^L_{\wt X}$ is unirational, which implies that 
$$Pic(\SU^L_{\wt X}\times J^0_{\wt X})=Pic(\SU^L_{\wt X})\times 
Pic(J^0_{\wt X}).$$ 
Hence it suffices to check that the restriction to each factor is ample.
The restriction to the first factor $\SU^L_{\wt X}$ is $\Theta_{\SU_{\wt X}}|
_{\SU^L_{\wt X}}$, which is clearly ample.

The restriction to the second factor is 
$f^*(\Theta_{\SU_{\wt X}})|_{J^0_{\wt X}}\otimes f^*(Det^*\Theta_y)^{-2}|_{J^0_{\wt X}}$. Write
$M_1=f^*(\Theta_{\SU_{\wt X}})|_{J^0_{\wt X}}$ and
$M_2=f^*(Det^*\Theta_y)^{-2}|_{J^0_{\wt X}}$, we are left with the task of
proving that $M_1\otimes M_2$ is ample. It is easy to
see that $M_1=f_1^*\Theta_{\SU_{\wt X}}$ and $M_2=f_2^*(\Theta_y^{-2})$, where
$f_1:J^0_{\wt X}\to\SU_{\wt X}$ and $f_2:J^0_{\wt X}\to J^d_{\wt X}$ are
given by $f_1(L_0)=E\otimes L_0$ (for a fixed $E$) and $f_2(L_0)=L_0^r\otimes L$. 
If we identify $J^0_{\wt X}$ with $J^d_{\wt X}$ by the isomorphism
$J^0_{\wt X}@>\otimes L>>J^d_{\wt X}$ and work up to algebraic equivalence, 
then $M_2=[r]^*(\Theta_y^{-2})$ is algebraically equivalent to $\Theta_y^{-2r^2}$, 
where $[r]:J^0_{\wt X}\to J^0_{\wt X}$ is the finite
cover given by $[r](L_0)=L_0^r$. To figure $M_1$ out, we consider the commutative diagram
$$\CD
\wt X\times J^0_{\wt X}  @>1\times f_1>>  \wt X\times\SU_{\wt X} \\
@V\pi_J VV                @VVV          \\
J^0_{\wt X} @>f_1>> \SU_{\wt X}.
\endCD$$ 
By the base change theorem, if $\SL$ denote a Poincar{\'e} bundle on
$\wt X\times J^0_{\wt X}$, then
$$M_1=(det\,R\pi_JE\otimes\SL)^k\otimes\bigotimes_{x\in I}
\{(det\,\SL_x)^{r\alpha_x}\otimes\bigotimes^{l_x}_{i=1}
(det\,\SL_x)^{r_i(x)d_i(x)}\}\otimes 
(det\,\SL_y)^{r\tilde\ell},$$
which is clearly algebraic equivalent to 
$$\aligned (det\,R\pi_JE\otimes\SL)^k&\otimes (det\,\SL_y)^{\sum_{x\in I}
\sum^{l_x}_{i=1}d_i(x)r_i(x)+r\sum_{x\in I}\alpha_x+r\tilde{\ell}}\\
&=(det\,R\pi_JE\otimes\SL)^k\otimes (det\,\SL_y)^{k\tilde n}.\endaligned$$
On the other hand, since $E$ is generated by sections and $det(E)=L$, we
have
$$0\to\SO_{\wt X}\otimes\Bbb C^{r-1}\to E\to L\to 0.$$
Thus $(det\,R\pi_JE\otimes\SL)=(det\,R\pi_J\SL)^{r-1}\otimes det\,R\pi_J(L\otimes \SL),$ 
and $M_1$ is algebraically equivalent to
$$\{det\,R\pi_J\SL\otimes (det\,\SL_y)^{1-\tilde g}\}^{(r-1)k}\otimes 
\{det\,R\pi_J(L\otimes \SL)\otimes (det(L\otimes\SL)_y)^{d+1-\tilde g}\}^k.$$
After identifying $J^d_{\wt X}$ with $J^0_{\wt X}$, we see that $M_1$ is
algebraically equivalent to $\Theta_y^{rk}$. Thus $M_1\otimes M_2$ is 
algebraically equivalent to $\Theta_y^{rk-2r^2}$, which is clearly ample
when $k>2r$.
\enddemo

The next lemma is a copy of Lemma 4.17 of [NR], one can see [Kn] for its
detail proof.
 
\proclaim{Lemma 5.4} Let $X$ be a normal, Cohen-Macaulay variety on which a 
reductive group $G$ acts, such that a good quotient $\pi:X\to Y$ exists.
Suppose that the action is generically free and that $dim(G)=dim(X)-dim(Y)$,
and further suppose that\roster
\item the subset where the action is not free has codimension $\ge 2$, and
\item for every prime divisor $D$ in $X$, $\pi(D)$ has codimension $\le 1$,
where $D$ need not be invariant.\endroster
Then $\omega_Y=(\pi_*\omega_X)^G$ where $\omega_X$, $\omega_Y$ are the respective 
dualising sheaves and the superscript $(\quad)^G$ denotes the
$G$-invariant direct image.\endproclaim

Fix an ample line bundle $\SO(1)$ on $\wt X$, and a set of data
$$\omega=(d,r,k,\tilde\ell,\{d_i(x)\}_{x\in I,1\le i\le l_x},\{\alpha_x\}_{
x\in I},I)$$
satisfying
$$ \sum_{x\in I}\sum^{l_x}_{i=1}d_i(x)r_i(x)+r\sum_{x\in
I}
\alpha_x+r\tilde{\ell}=k\tilde n,\tag5.10$$
$\omega$ determines a polarisation (for fixed $\SO(1)$)
$$\frac{\tilde\ell}{m}\times\prod_{x\in I}
\{\alpha_x,d_1(x),\cdots,d_{l_x}(x)\}.$$
We denote the set of semistable points for the $SL(\tilde n)$ action
under this polarisation by $\wt\SR_{\omega}^{ss}\subset\wt\SR_F$, and its
good quotient by $\SU_{\wt X,\omega}$.
$$\Theta_{\wt\SR^{ss}_{\omega}}=
(det\,R\pi_{\wt\SR^{ss}_{\omega}}\SE)^k\otimes\bigotimes_{x\in I}
\{(det\,\SE_x)^{\alpha_x}\otimes\bigotimes^{l_x}_{i=1}
(det\,\SQ_{x,i})^{d_i(x)}\}\otimes 
(det\,\SE_y)^{\tilde\ell}$$
descends to an ample line bundle $\Theta_{\SU_{\wt X,\omega}}$ on $\SU_{\wt X,\omega}$, 
and we need to prove that
$$H^1(\SU_{\wt X,\omega},\Theta_{\SU_{\wt X,\omega}})=0.$$

\proclaim{Theorem 5.1} Assume that $\tilde g\ge 2$. Then, for any set of
data $\omega$ satisfying (5.10),
$$H^1(\SU_{\wt X,\omega},\Theta_{\SU_{\wt X,\omega}})=0.$$
\endproclaim

\demo{Proof} We can assume that $r>2$ since the vanishing theorem for
$r=2$ is known (See [NR] and [Ra]). Let $\bar\omega=(d,r,\bar k,\bar\ell,
\{\bar d_i(x)\}_{x\in I,1\le i\le l_x},\{\bar\alpha_x\}_{
x\in I},I)$ 
be a new set of data with 
$\bar k=k+2r,$ $\bar\ell=2\tilde n+\tilde\ell-r|I|,$ 
$\bar d_i(x)=d_i(x)+n_i(x)+n_{i+1}(x)$ and $\bar\alpha_x=\alpha_x+n_{l_x+1}(x)$, let
$$\aligned\hat\Theta_{\bar\omega}:=&(det\,R\pi_{\wt\SR_F}\SE)^{\bar k}
\otimes\bigotimes_{x\in I}\left\{(det\,\SE_x)^{\bar\alpha_x}
\otimes\bigotimes^{l_x}_{i=1}(det\,\SQ_{x,i})
^{\bar d_i(x)}\right\}\\
&\bigotimes_{x\in I}(det\,\SE_x)^{-r}\otimes\bigotimes_q(det\,\SE_q)^{1-r}
\otimes(det\,\SE_y)^{2\tilde n+2(r-1)(\tilde g-1)+\tilde\ell}\endaligned$$
One can check that
$$ \sum_{x\in I}\sum^{l_x}_{i=1}\bar d_i(x)r_i(x)+r\sum_{x\in
I}\bar\alpha_x+r\bar{\ell}=\bar k\tilde n.\tag5.11$$
$\bar\omega$ determines a new polarisation
$$\frac{\bar\ell}{m}\times\prod_{x\in I}
\{\bar\alpha_x,\bar d_1(x),\cdots,\bar d_{l_x}(x)\}.$$
We denote the set of semistable points for the $SL(\tilde n)$ action
under the new polarisation by $\wt\SR_{\bar\omega}^{ss}\subset\wt\SR_F$, 
and its good quotient
$$\psi_{\bar\omega} :\wt\SR^{ss}_{\bar\omega}\to\SU_{\wt X,\bar\omega}.$$
$\hat\Theta_{\bar\omega}$ descends to an ample line bundle 
$\Theta_{\bar\omega}$ (See Remark 1.1 (2)). By Proposition 2.2 and (5.9), we have
$$\Theta_{\wt\SR_F}\otimes\omega^{-1}_{\wt\SR_F}=\hat\Theta_{\bar\omega}
\otimes Det^*\Theta_y^{-2}.\tag5.12$$
Since we assumed that $\tilde g\ge 2$ and $r>2$, the codimension
of $\wt\SR_F\ssm\wt\SR^{ss}_{\omega}$ for any $\omega$ is at least $3$
(See Proposition 5.1 (2)). Thus, by local cohomology theory, we have
$$H^1(\wt\SR_{\omega}^{ss},\Theta_{\wt\SR_{\omega}^{ss}})^{inv}=
H^1(\wt\SR_F,\Theta_{\wt\SR_F})^{inv}=H^1(\wt\SR^{ss}_{\bar\omega},
\Theta_{\wt\SR^{ss}_{\bar\omega}})^{inv}.\tag5.13$$
Since $codim(\wt\SR^{ss}_{\bar\omega}\ssm\wt\SR^s_{\bar\omega})\ge 2$ 
(See Proposition 5.1 (1)), by using Lemma 5.4, we have 
$$(\psi_{\bar\omega *}\omega_{\wt\SR^{ss}_{\bar\omega}})^{inv}
=\omega_{\SU_{\wt X,\bar\omega}}$$
(See Lemma 6.3 of [NR]). By (5.12), we can write 
$$\Theta_{\wt\SR^{ss}_{\bar\omega}}=\psi_{\bar\omega}^*(\Theta_
{\bar\omega}\otimes Det^*\Theta_y^{-2})
\otimes\omega_{\wt\SR^{ss}_{\bar\omega}}.\tag5.14$$
One use the fact that for good quotients the space of invariants of the
cohomology of an invariant line bundle is the same as 
the cohomology of the invariant direct image and (5.14) to prove that
$$\aligned
H^1(\wt\SR^{ss}_{\bar\omega},\Theta_{\wt\SR^{ss}_{\bar\omega}})^{inv}
&=H^1(\SU_{\wt X,\bar\omega},\Theta_{\bar\omega}\otimes Det^*\Theta_y^{-2} 
\otimes(\psi_{\bar\omega *}\omega_{\wt\SR^{ss}_{\bar\omega}})^{inv})\\
&=H^1(\SU_{\wt X,\bar\omega},\Theta_{\bar\omega}\otimes Det^*\Theta_y^{-2} 
\otimes\omega_{\SU_{\wt X,\bar\omega}}).\endaligned$$
Now since $\Theta_{\bar\omega}\otimes Det^*\Theta_y^{-2}$ is 
an ample line bundle by Lemma 5.3 (Note that $\bar k>2r$) and 
$\SU_{\wt X,\bar\omega}$ has only rational singularities, 
we can apply a Kodaira-type vanishing theorem 
(See Theorem 7.80(f) of [SS]) and conclude
that 
$$ H^1(\SU_{\wt X,\omega},\Theta_{\SU_{\wt X,\omega}})
=H^1 (\wt\SR^{ss}_{\bar\omega},\Theta_{\wt\SR^{ss}_{\bar\omega}})^{inv}=0.$$
\enddemo

\remark{Remark 5.2} (1) To check $(5.11)$, one has to show that for any $x\in I$
$$r(n_{l_x+1}(x)-r)+\sum^{l_x}_{i=1}r_i(x)(n_i(x)+n_{i+1}(x))=0.$$
Note that $r=n_{l_x+1}(x)+\sum^{l_x}_{i=1}n_i(x)$ and $n_i(x)=r_i(x)-r_{i-1}(x),$ 
we have to show that 
$$\sum^{l_x}_{i=1}\{r_{i-1}(x)(r-r_i(x))-r_i(x)(r-r_{i+1}(x))\}=0,$$
which is clearly true since $r_0(x)=0$.

(2) Since $R^i\rho_*(\hat\Theta')=0$ for all $i>0$ 
(Note that $\rho:\wt\SR'_F\to\wt\SR_F$ is a grassmannian bundle over $\wt\SR_F$), 
we have
$H^1(\wt\SR'_F,\hat\Theta')^{inv}=H^1(\wt\SR_F,\rho_*\hat\Theta')^{inv}.$
By using the canonical decomposition (See Remark 4.2) and the 
vanishing Theorem 5.1, we can show that
$H^1(\wt\SR'_F,\hat\Theta')^{inv}=0.$
\endremark

Next we will show the vanishing theorem for the moduli space
of semistable parabolic torsion free sheaves on a nodal curve $X$.

\proclaim{Theorem 5.2} Assume that $g\ge 3$. Then $H^1(\SU_X,\Theta_{\SU_X})=0.$
\endproclaim

\demo{Proof} It will be reduced to prove a vanishing theorem for $\SP$ since
the following lemma.\enddemo

\proclaim{Lemma 5.5} For $0\le a\le r$, the natural maps 
$H^1(\SW_a,\Theta_{\SU_X})\to H^1(\SD_1(a),\Theta_{\SP})$ are injective. In
particular, $H^1(\SU_X,\Theta_{\SU_X})\to H^1(\SP,\Theta_{\SP})$ is injective.
\endproclaim

\demo{Proof} It is known that $\phi_a:=\phi|_{\SD_1(a)}:\SD_1(a)\to\SW_a$ is
the normalisation of $\SW_a$ (See Proposition 2.1). If $\SW_{a-1}$ is empty,
$\phi_a$ is an isomorphism and the lemma is clear. If $\SW_{a-1}$ is not empty,
$\SW_{a-1}$ is the non-normal locus of $\SW_a$ and we are reduced to
prove that $$H^1(\SW_{a-1},\Theta_{\SU_X})\to H^1(\SD_1(a)\cap\SD_2\cup\SD_1(a-1), 
\Theta_{\SP})$$
is injective by Lemma 4.3 (2). Thus it is enough to show that
$$H^1(\SW_{a-1},\Theta_{\SU_X})\to H^1(\SD_1(a-1),\Theta_{\SP})$$ 
is injective, and we are done by induction since $\phi_0:\SD_1(0)\to\SW_0$
is always an isomorphism. \enddemo

In order to prove the vanishing theorem for
$\SP$, we have to prepare some lemmas. 

\proclaim{Lemma 5.6} Assume $\tilde g\ge 2$.
Then $(\tilde\psi'_*\omega_{\SH})^{inv}
=\omega_{\SP}$ where $\omega_{\SP}$ is the
canonical (dualising) sheaf of $\SP.$\endproclaim

\demo{Proof} We will check the conditions of
Lemma 5.4. By Proposition 5.1 (2), $\SU_{\wt X}(d-r)$ contains a stable bundle,
and thus $\SW_0$ contains a stable parabolic sheaf
by Lemma 2.8, which shows that there exist stable
parabolic bundles on $X$ since stability is an open condition. 
Thus there exist stable generalised parabolic bundles on $\wt X$ by 
Lemma 2.2 (2), and the action of $PGL(\tilde n)$
on $\SH$ is therefore generically free. We now
check conditions (1) and (2) of Lemma 5.4.\roster
\item By Proposition 5.2 (2), the nonstable locus in 
$\wt\SR^{\prime ss}\ssm\{\SD_1^f\cup\SD^f_2\}$ has codimension $\ge 2$.
We need to show that each of the $\hat\SD_j(r-1)$ and $\hat\SD_j^t$ contains
GPS with no automorphism except scales. Take $j=1$ for definiteness, let
$\wt E$ be a stable parabolic bundle on $\wt X$ of degree $d-r$, let
$E=\wt E\oplus\,_{x_2}\Bbb C^r$ and define the GPS structure on $E$ as follows.
We take $Q=\Bbb C^r$, the map $E_{x_2}\to Q$ to be the obvious projection, and
the map $E_{x_1}\to Q$ any isomorphism. This yields, after an identification
$H^0(E)\cong \Bbb C^{\tilde n}$, a point on $\hat\SD_1^t$ as required. Next
consider $E=\wt E\otimes\SO_{\wt X}(x_2)$, the GPS structure being given by
taking $Q=\wt E_{x_2}$, the map $E_{x_1}\to Q$ being zero, and the map
$E_{x_2}\to Q$ the residue $\SO_{x_2}(x_2)\cong \Bbb C$. This yields a point on
$\hat\SD_1(r-1)$ with only automorphisms by scales.
\item  If a prime divisor is not contained in the nonstable locus, its image
in $\SP$ will have codimension one. If it is contained in the nonstable locus,
then, by (2) of Proposition 5.2, it has to be one of the
$(\hat\SD_j(r-1))^{ss}$ and $(\hat\SD_j^t)^{ss}$. We have already seen that
the respective images
of these in $\SP$ are the $\SD_j$ by Proposition 3.3.
\endroster\enddemo

\proclaim{Lemma 5.7} There is a morphism $Det:\SH\to J^d_{\wt X}$ 
which extends the determinant morphism on the open set $\wt\SR'_F$. Moreover, 
it yields a
flat morphism $Det:\SP\to J^d_{\wt X}.$\endproclaim

\demo{Proof} Note that, on $\wt X\times\SH$, we have an exact squence
$$0\to\SK\to\SO^{\tilde n}\to\SE\to 0,$$
and $\SK$ is flat over $\SH$ since $\SE$
is so. One proves that $\SK$ is locally free on $\wt X\times\SH$
(By using Lemma 5.4 of [Ne]). Thus $det(\SK)^{-1}$ is a line bundle on
$\wt X\times \SH$, and gives a morphism 
$$Det:\SH\to J^d_{\wt X},$$
which is clearly an extension of the determinant morphism on the open set $\wt\SR'_F.$
Restricted to $\wt\SR^{\prime ss}$ the map $Det$ clearly factors through 
the quotient by the $SL(\tilde n)$ action and yields a morphism
$$Det:\SP\to J^d_{\wt X},$$
which we will prove to be a flat morphism. $J^0_X$ acts on $\SP$ by
$$(E,Q)\mapsto L\cdot(E,Q):=(E\otimes\pi^*L,Q\otimes L_{x_0}).$$
One checks that $Det(L\cdot(E,Q))=Det(E,Q)\otimes(\pi^*L)^r$. Note that
the pull-back map $J^0_X\to J^0_{\wt X}$ and the $r$-power map $J^0_{\wt X}\to
J^0_{\wt X}$ are surjective, and $J^0_{\wt X}$ acts transitively on 
$J^d_{\wt X},$ we can see that $Det:\SP\to J^d_{\wt X}$ is flat by generic flatness. 
\enddemo

Let $\SH^L$ denote the (reduced) fibre over $L\in J^d_{\wt X}$, and $\SP^L$
denote the (reduced) fibre of $Det$ above $L$. Clearly $\SP^L$ is the GIT 
quotient of $\SH^L$, and all of the properties of $\SH$ and $\SP$ continue to
be valid for $\SH^L$ and $\SP^L$. From the proof of above lemma, one sees that
all of the fibres of $Det:\SP\to J^d_{\wt X}$ are reduced. Thus $\SP^L$ is 
also the scheme-theoretic fibre over $L$, and we have

\proclaim{Proposition 5.3} The canonical (dualising) sheaf of $\SP^L$ is the 
restriction of $\omega_{\SP}$ to $\SP^L$.\endproclaim

\demo{Proof} The following general fact can be proved by repeated use of 
Bertini (on $U$) and the adjunction formula: Suppose $f:V\to U$ is a flat
map of varieties, with $U$ smooth and $V$ Gorenstein. Let $V_p$ be the 
scheme-theoretic fibre over $p\in U$. Then the dualising sheaf of $V_p$ is
the restriction of the dualising sheaf of $V$.\enddemo 

\proclaim{Proposition 5.4} Assume $\tilde g\ge 2$. Then 
$H^1(\SP^L,\Theta_{\SP})=0$ for any $L\in J^d_{\wt X}$.\endproclaim

\demo{Proof} Let $\omega^L_{\SH}$ denote the restriction of $\omega_{\SH}$
to $\SH^L$. Then $(\tilde\psi'_*\omega^L_{\SH})^{inv}=\omega_{\SP^L}$ by
Lemma 5.6 and Proposition 5.3. Recall that, for the polarisation
$$\frac{(\tilde{\ell}-k)}{m}\times k\times\prod_{x\in
I}\{\alpha_x,d_1(x),\cdots,d_{l_x}(x)\},$$
the line bundle $\hat\Theta'$ was defined to be
$$\aligned &(det\,R\pi_{\SH}\SE)^k\otimes\bigotimes_{x\in I}
\lbrace(det\,\SE_x)^{\alpha_x}\otimes\bigotimes^{l_x}_{i=1}
(det\,\SQ_{x,i})^{d_i(x)}\rbrace\otimes(det\,\SE_y)^{\tilde\ell}\\
&\otimes(det\,\SQ)^k\otimes(det\,\SE_y)^{-k},\endaligned$$
which descends to the ample line bundle $\Theta_{\SP}$ if the
polarisation satisfies 
$$ \sum_{x\in I}\sum^{l_x}_{i=1}d_i(x)r_i(x)+r\sum_{x\in I}
\alpha_x+r\tilde{\ell}=k\tilde n.$$
Note that $(1\times Det)^*\SL=(det\,\SK)^{-1}\otimes\pi^*_{\SH}\Cal N$
for a suitable line bundle $\Cal N$ on $\SH$, one sees that 
$\Cal N\cong det\,\SK_x$ on $\SH^L$ for any $x\in\wt X$ ($x$ may be $x_1$ and $x_2$), 
and 
$$(det\,R\pi_{\SH^L}det\,\SK^{-1})^{-2}=(det\,\SE_y)^{2\tilde n+(r-1)(2\tilde g-2)}.$$
Thus, on $\SH^L$, we have $\hat\Theta'=
\hat\Theta'_{\omega}\otimes\omega^L_{\SH}$ by Proposition 3.4, where 
$$\aligned\hat\Theta'_{\omega}=&(det\,R\pi_{\SH^L}\SE)^{\bar k}
\otimes\bigotimes_{x\in I}
\lbrace(det\,\SE_x)^{\bar\alpha_x}\otimes\bigotimes^{l_x}_{i=1}
(det\,\SQ_{x,i})^{\bar d_i(x)}\rbrace
\otimes(det\,\SE_y)^{\bar\ell}\\
&\otimes(det\,\SQ)^{\bar k}\otimes(det\,\SE_y)^{-\bar k}\endaligned$$
with $\bar k=k+2r,$ $\bar\ell=2\tilde n+\tilde\ell-r|I|,$ 
$\bar d_i(x)=d_i(x)+n_i(x)+n_{i+1}(x)$ and $\bar\alpha_x=\alpha_x+n_{l_x+1}(x)$.
One checks that
$$ \sum_{x\in I}\sum^{l_x}_{i=1}\bar d_i(x)r_i(x)+r\sum_{x\in
I}\bar\alpha_x+r\bar{\ell}=\bar k\tilde n.$$
The rest of the proof proceeds as Theorem 5.1 except that an analogue of
Lemma 5.3 is not needed. The Kodaira-type vanishing theorem and Hartogs-type
extension theorem for cohomology are applicable since $\SH^L$ and $\SP^L$
are Cohen-Macaulay and have only rational singularities.
\enddemo

\proclaim{Theorem 5.3} Assume $\tilde g\ge 2$. Then $H^1(\SP,\Theta_{\SP})=0$.
\endproclaim

\demo{Proof} We consider the flat morphism $Det:\SP\to J^d_{\wt X}$ and
try to decompose the direct image $(Det)_*\Theta_{\SP}$. One can see that
$(Det)_*\Theta_{\SP}=\{(Det_{\wt\SR^{\prime ss}})_*\hat\Theta'\}^{inv}$ and
the equalities 
$$ \{(Det_{\wt\SR^{\prime ss}})_*\hat\Theta'\}^{inv}=
\{(Det_{\SH})_*\hat\Theta'\}^{inv}=
\{(Det_{\wt\SR'_F})_*\hat\Theta'\}^{inv}$$
hold by using Lemma 4.1 and Lemma 4.2, where $Det_{\wt\SR'_F}:\wt\SR'_F\to
J^d_{\wt X}$ is clearly factorized through the projection 
$\rho:\wt\SR'_F\to\wt\SR_F$. Thus $(Det_{\wt\SR'_F})_*\hat\Theta'=
(Det_{\wt\SR_F})_*\rho_*\hat\Theta'$ and, by Remark 4.2, we have
$$(Det_{\wt\SR'_F})_*\hat\Theta'=\bigoplus_{\mu}
(Det_{\wt\SR^{\mu}_F})_*\hat\Theta_{\mu},$$
where $Det_{\wt\SR^{\mu}_F}=Det_{\wt\SR_F}\cdot p^{\mu}:\wt\SR^{\mu}_F\to
J^d_{\wt X}$, which restricting to $(\wt\SR^{\mu})^{ss}$ induces a morphism
$Det_{\mu}:\SU^{\mu}_{\wt X}\to J^d_{\wt X}.$ It is now clear that we have 
the decomposition  
$$(Det)_*\Theta_{\SP}=\bigoplus_{\mu}
(Det_{\mu})_*\Theta_{\mu},$$
which implies that $H^1(J^d_{\wt X},(Det)_*\Theta_{\SP})=0$ by Theorem 5.1
since $\tilde g\ge 2$. On the other hand, $R^1(Det)_*(\Theta_{\SP})=0$ by
Proposition 5.4. Hence we are done by using spectral sequence.\enddemo

\bigskip

\remark{Acknowledgements: Part of this work was done during my stay at the 
Mathematics Section of the Abdus Salam ICTP, Trieste, as a visiting fellow. 
I would like to express my hearty thanks to Professor M.S. Narasimhan, who
introduced me to this subject and explained patiently the ideas of [NR].
This paper would not have been possible without his help and encouragement.
It is my pleasure to thank T.R. Ramadas , 
who answered a number of my questions through e-mails and fax, 
and this was very helpful. I thank R. H{\"u}bl very much, who helped me
to understand German literature [Kn]}\endremark


\bigskip

\Refs

\widestnumber\key{B1}
\widestnumber\key{B2}
\widestnumber\key{BR}
\widestnumber\key{CS}
\widestnumber\key{DN}
\widestnumber\key{Fa}
\widestnumber\key{La}
\widestnumber\key{MS}
\widestnumber\key{Ne}
\widestnumber\key{NR}
\widestnumber\key{NS}
\widestnumber\key{Pa}
\widestnumber\key{S1}
\widestnumber\key{S2}
\widestnumber\key{Si}
\widestnumber\key{Tr}
\widestnumber\key{EGA-I}


\ref\key B1 \by U. Bhosle\paper Generalised parabolic bundles and
applications to torsionfree sheaves on nodal curves
\pages 187--215\yr1992\vol 30 \jour Arkiv f{\"o}r matematik\endref

\ref\key B2 \by U. Bhosle\paper Generalized parabolic bundles and
applications--II\pages 403--420 \yr1996\vol 106 \jour Proc. Indian
Acad. Sci. (Math. Sci.)\endref

\ref\key BR \by I. Biswas and N. Raghavendra\paper Determinants of
parabolic bundles on Riemann surfaces\pages 41--71\yr1993\vol 103
\jour Proc. Indian Acad. Sci.(Math. Sci.)\endref

\ref\key CS \by C. De Concini and E. Strickland\paper On the variety
of complexes \pages 57--77\yr1981\vol 41\jour Adv. in Math.\endref

\ref\key DN\by J.M. Drezet and M.S. Narasimhan\paper Groupe de
Picard des vari{\'e}t{\'e}s de fibr{\'e}s semistables sur les courbes
alg{\'e}briques\pages 53--94\yr1989\vol 97\jour Invent. Math.\endref

\ref\key EGA-I\by A. Grothendieck and J. Dieudonn{\'e}\book
{\'E}l{\'e}ments de G{\'e}om{\'e}trie
alg{\'e}brique I\bookinfo Grundlehren 166\publaddr
Berlin-Heidelberg-New York: Springer\yr1971\endref

\ref \key Fa \by G. Faltings\paper Moduli-stacks for bundles on semistable
curves\pages 489--515\yr1996\vol 304\jour Math. Ann.\endref

\ref\key FH \by W. Fulton and J. Harris\book Representation Theory : 
A first course\yr1991\bookinfo Graduate Texts in Mathematics
\vol 129\publaddr Springer-Verlag New York Inc.\endref

\ref\key Fu \by W. Fulton\paper Flags, Schubert polynomials, degeneracy 
loci and determinantal formula \yr1992\jour Duke Math. J.\vol 65\pages
381--420\endref

\ref \key Gi \by D. Gieseker\paper On the Moduli of vector bundles on
an algebraic surface\pages 45--60\yr1977\vol 106\jour Annals of
Math.\endref

\ref\key Ha \by R. Hartshorne\book Algebraic geometry\yr1977\publaddr 
Berlin-Heidelberg-New York. Springer\endref

\ref \key He \by W.H. Hesselink\paper Desingularizations of
varieties of nullforms\pages 141--163\yr1979\vol 55\jour Invent.
Math.\endref

\ref \key Kn \by F. Knop\paper Der kanonische Moduleines Invariantenrings \pages
40--54\yr1989\vol 127\jour Joural of algebra\endref

\ref \key La \by H. Lange\paper Universal families of extensions \pages
101--112\yr1983\vol 83\jour Joural of algebra\endref

\ref \key MRa\by V.B.Mehta and A.Ramanathan\paper Frobenius splitting and 
cohomology vanishing for Schubert varieties\pages 27--40\yr1985\vol
122\jour Annals of Math. \endref

\ref \key MR\by V.B.Mehta and T.R. Ramadas\paper Moduli of vector bundles,
Frobenius splitting, and invariant theory\pages 269--313\yr1996\vol
144\jour
Annals of Math. \endref

\ref \key MS\by V.B.Mehta and C.S.Seshadri\paper Moduli of vector bundles
on curves with parabolic structures\pages 205--239\yr1980\vol 248\jour
Math. Ann. \endref

\ref\key Ne \by P.E.Newstead\book Introduction to
moduli problems
and orbit spaces\yr1978\bookinfo TIFR lecture
notes\publaddr New
Delhi: Narosa\endref

\ref\key NR \by M.S.Narasimhan and
T.R. Ramadas\paper Factorisation of 
generalised theta functions I\pages
565--623\yr1993\vol 114 \jour
Invent. Math.\endref

\ref\key NS \by D.S.
Nagaraj and C.S. Seshadri\paper Degenerations
of the moduli spaces of
vector bundles on curves I\pages 101--137
\yr1997\vol 107 \jour Proc.
Indian Acad. Sci.(Math. Sci.)\endref

\ref\key Pa \by C. Pauly\paper
Espaces modules de fibr{\'e}s
paraboliques et blocs conformes\pages
217--235\yr1996\vol 84 \jour
Duke Math.\endref

\ref\key Ra \by T.R.
Ramadas\paper Factorisation of 
generalised theta functions II: the
Verlinde formula\pages
641--654\yr1996\vol 35
\jour
Topology\endref

\ref\key S1 \by C.S. Seshadri\paper Quotient spaces
modulo reductive
algebraic groups\pages 511-556\vol 95\jour Annals of
Math.\yr 1972\endref

\ref\key S2 \by C.S. Seshadri\paper Fibr{\'e}s
vectoriels sur les courbes
alg{\'e}briques\vol 96\jour Ast{\'e}risque\yr
1982\endref

\ref \key Si \by C. Simpson \paper Moduli of representations
of the
fundamental group of a smooth projective variety I
\pages
47--129\vol 79 \yr1994\jour I.H.E.S. Publications
Math{\'e}matiques\endref

\ref\key SS \by B. Schiffman and A.J. Sommense\book Vanishing 
theorems on complex manifolds \yr 1985
\publaddr Boston Basel Stuttgart: Birkha{\"u}ser\endref

\ref\key St \by E. Strickland\paper On the conormal bundle of
the determinantal variety \pages 523--537
\yr1982\vol 75\jour J. of Algebra\endref

\ref\key Sw \by R.G. Swan\paper On seminormality \pages 210--229
\yr1980\vol 67\jour J. of Algebra\endref

\ref \key Tr \by V. Trivedi\paper The seminormality property of
circular
complexes \pages 227--230\yr1991\vol 101\jour Proc. Indian Acad.
Sci.
(Math. Sci.) \endref 

\ref\key ZS \by O. Zariski and P. Samuel\book Commutative Algebra II
\yr1960\publaddr Princeton: Van Nostrand\endref



\endRefs
\enddocument